\documentclass{article}
\usepackage{style}
\usepackage{tikz}
\usepackage{enumitem}

\allowdisplaybreaks

\newcommand{\gl}{\mathfrak{gl}}
\newcommand{\cG}{\mathcal{G}}

\newcommand{\dev}{D}
\newcommand{\hol}{h}

\DeclareMathOperator{\diag}{diag}

\newcommand{\h}{\mathfrak{h}}

\newcommand{\bb}{\mathfrak{b}}
\newcommand{\bd}{\mathfrak{d}}
\newcommand{\bj}{\mathfrak{j}}
\newcommand{\bc}{\mathfrak{c}}
\newcommand{\w}[1]{\mathfrak{w}_{#1}}
\newcommand{\Sph}{\mathbb{S}}
\newcommand{\n}{\mathfrak{n}}

\newcommand{\Vt}{\widetilde V}
\newcommand{\Ob}{\overline\Omega}
\DeclareMathOperator{\sgn}{sgn}
\numberwithin{equation}{section}

\title{Projective structures on $\T^4$ of  rank 3 and 4}
\author{Arielle Leitner and  J.\ Maxwell Riestenberg\\ 
with an appendix by Willem de Graaf}
\date{\today}

\begin{document}

\vspace{-4cm}

\maketitle

\begin{abstract}  
Following Benoist's theory of bricks for projective manifolds with nilpotent holonomy
\cite{benoist1994nilvarietes}, \cite{Ben00}, we classify the projective and affine structures on $\mathbb T^4$ of
rank $4$, and those of rank 3 with respect to a nilpotent group without compact factors. For rank $4$, we classify the maximal abelian subgroups of
$\mathrm{SL}(5,\mathbb R)$ up to conjugacy, with an appendix by Willem de Graaf treating
the unipotent case by Galois cohomology, and determine which of them contain the group of
a left-invariant structure. For rank $3$, we introduce and determine invariant functions, which allow us to determine which bricks appear and how they can be glued together.
As a consequence of the classification, we find that, according to the holonomy, either every structure is affine or else there is an infinite family of structures with the same holonomy whose universal covers are pairwise non-isomorphic. 
We also revisit the rank $2$ structures on $\mathbb T^3$.%
\end{abstract}

\section{Introduction}
A projective (resp.\ affine) manifold is a $(G,X)$-manifold in the sense of Klein, with
$G=\mathrm{SL}(d+1,\R)$ and $X=\Sph^d$ (resp.\ $G=\mathrm{Aff}(\R^d)$ and $X=\R^d$); see \cite{thurston2014three}, \cite{goldman2022geometric}. The simplest example is the Euclidean torus $\T^2=\mathbb Z^2\backslash\R^2$,
with holonomy given by translations. Affine structures on $\T^2$ were classified by Nagano
and Yagi \cite{NY74} and by Furness and Arrowsmith \cite{furness1975locally}; for an overview of the subject and
further results, see \cite{goldman2022geometric} and, for instance, \cite{baues2008prehomogeneous}, \cite{smillie1977affinely}, \cite{choi2001decomposition}, \cite{fried1981affine}.  
 
Unlike Riemannian structures, a closed affine manifold need not be complete. The basic
example is the Hopf torus $\langle d_\lambda\rangle\backslash(\R^2\setminus\{0\})$, where $d_\lambda$
is the dilation by $\lambda>1$; its developing image is the punctured plane. 
For affine structures on $\T^d$ with $d \ge 3$, the developing map need not even be a covering space onto its image. 
To control
such phenomena for manifolds with nilpotent holonomy, Benoist \cite{benoist1994nilvarietes}, \cite{Ben00} developed a
theory of \emph{bricks}. If the holonomy group $H$ is nilpotent, a finite-index subgroup of
$H$ lies in a maximal connected nilpotent subgroup $N$ of $G$, and the developing image is a
union of the finitely many $N$-orbits. The \emph{rank} of the structure is the smallest
dimension of an $N$-orbit in the developing image (Proposition~ \ref{prop:rank}). Up to a finite cover,
the manifold is obtained by gluing finitely many bricks, quotients of closures of open
$N$-orbits, along their faces (Theorem \ref{theorem:benoistbricks} due to \cite{benoist1994nilvarietes}). A structure of rank $d$ is left invariant and
consists of a single homogeneous piece; in rank $d-1$ the bricks are glued along a
one-dimensional pattern, a ``word of bricks'', which determines the structure.

\subsection*{Rank 4} 
Full rank structures are a mild generalization of \emph{Kleinian} structures, i.e.\ the quotient of a domain in $X$ by a free proper cocompact action of a discrete group. 
The developing image of a rank $4$ structure on $\T^4$ is an open orbit of a connected $4$-dimensional abelian group $B$ in $\Sph^4$. 
Hence they can be viewed as left-invariant structures, see Definition \ref{def:rank4}. 
Every such $B$ lies in a maximal abelian subgroup of $\mathrm{SL}(5,\R)$, so
we first classify these.
 
\begin{theorem}[with W.\ de Graaf]\label{thm:intro-masa}
Up to conjugacy, the identity components of the maximal abelian subgroups of
$\mathrm{SL}(5,\R)$ are the $23$ unipotent groups of Theorem \ref{thm:unipotent-list} and the $24$ non-unipotent
groups of Table \ref{table:nonunip}.
\end{theorem}
 
The unipotent case is proved in Appendix \ref{appendix:realforms} by Willem de Graaf using Galois cohomology, starting from the
complex classification of \cite{suprunenko1968commutative}; the non-unipotent groups are sums of blocks
(Lemma \ref{lem:real-blocks} and Appendix \ref{appendix:nonunipotent}).
 
\begin{theorem}\label{thm:intro-rank4}
Every rank $4$ projective structure on $\T^4$ is the left-invariant structure defined by a marked lattice in the universal cover of a 
$4$-dimensional connected abelian subgroup $B$, with an open orbit in $\Sph^4$, of one of the
following groups: the $11$ unipotent groups of Theorem \ref{thm:unipotent-rank4}, whose structures are exactly
the complete affine ones, or the $18$ non-unipotent groups marked in Table \ref{table:nonunip}.
Conversely, every such $B$ and every marked lattice in its universal cover define a rank $4$
structure.
\end{theorem}
 
\subsection*{Rank 3} Here the structure is no longer homogeneous, and the main work is to
decide which faces can bound a brick and how bricks can be glued. For this we use
\emph{invariant functions} (Definition \ref{def:invariant-function}): a function on an open orbit whose level sets
are the orbits of the $3$-dimensional group $B$. Its behavior at the faces, which is of
pole type or of logarithmic type (Lemma \ref{lem:limits}), determines which pairs of faces bound a
brick, and whether a face keeps or changes its end when one passes to the neighbouring
copy of the open orbit (Definition \ref{def:constant-alternating}). When all faces of the second kind lie at the
same end, they can be crossed at most once, and hence never in a torus (Lemma \ref{lem:constant-faces}). Conversely, every admissible periodic word
of bricks is realized by a torus (Lemma \ref{lem:realize}). Our main result is the classification in
Theorem \ref{thm:rank3-summary}, which describes, for each of the seven maximal connected nilpotent groups $N$, the possible subgroups $B$ and
the bricks that occur. Among its consequences:
 
\begin{theorem}\label{thm:intro-rank3}
Let $V=\T^4$ carry a projective structure of rank $3$ with holonomy group almost contained in $N$ without elliptic elements (see Remark \ref{rem:compact-factors}).
\begin{enumerate}
  \item If $N$ is unipotent, there is no such structure (Proposition \ref{prop:N5}); this is the
  case $(n,p)=(4,3)$ of \cite{Ben00}, Exercise~4.6(c).
  \item Either $B$ admits a single brick, in which case the structure is affine and its words
  of bricks are the powers of this brick, or $B$ admits at least two bricks, and there are infinitely many structures whose universal covers are pairwise non-isomorphic by isomorphisms normalizing $N$ 
(Corollary~ \ref{cor:rank3-dichotomy}). 
  \item A structure with at least two bricks determines $N$ up to conjugation, while a
  single-brick structure can occur for two groups $N$, in exactly three families
  (See \Cref{sec:overlaps}).
\end{enumerate}
\end{theorem}
 
\subsection*{Rank 2 on $\T^3$} In Appendix \ref{appendix:T3} we apply the same methods to the rank $2$
structures on $\T^3$, which are the subject of \cite{Ben00}, Exercises~4.13--4.15. We give
complete solutions; for instance, for $N_{3,1}$ there are rank $2$ projective structures, none of them affine
(Proposition~ \ref{prop:B31}), and for $N_{2,1,1}$ which faces bound a common brick depends on both
generators of $\mathfrak b$ (Example \ref{ex:N211}).
 
\begin{remark}
Rank $2$ structures on $\T^4$ exist, for example products of two Hopf tori with respect to
the diagonal group. For these an invariant function takes values in $\R^2$ rather than
$\R$, and the gluing of bricks is no longer governed by a one-dimensional word; we do not
treat them here.
\end{remark}
 
\subsection*{Organization} Section \ref{section:prelim} recalls Benoist's theory and the notion of rank.
Section \ref{sec:bricks} develops invariant functions, the gluing rules and the realization of words of
bricks. Section \ref{sec:rank3} classifies the rank $3$ structures on $\T^4$, and Section \ref{sec:max_rannk} the rank $4$
ones. Appendix \ref{appendix:realforms}, by Willem de Graaf, classifies the real forms of the maximal abelian
nilpotent subalgebras of $\mathfrak{gl}(5,\mathbb C)$; Appendix \ref{appendix:T3} treats rank $2$ structures on
$\T^3$; Appendix \ref{appendix:nonunipotent} lists the non-unipotent maximal abelian subgroups of $\mathrm{SL}(5,\R)$.
 
\subsection*{Acknowledgments} We thank Anna Wienhard, Bill Goldman, Yves Benoist, John
Smillie and Charles Daly for helpful conversations. We thank Mikhail Borovoi for answering
our question on MathOverflow, pointing us to Galois cohomology, and connecting us with
Willem de Graaf. We thank the Max Planck Institute, Leipzig, for providing a paradise for
mathematical research.
 
\subsection*{AI usage statement}
The core mathematical ideas and content of this article were developed without the use of large language models. 
A complete first version of the article was written by the authors without AI assistance. 
We then used Claude to review this version, which identified a number of inefficiencies in the exposition. 
We subsequently supplied Claude with our earlier notes and computations, including Mathematica notebooks predating 2026 that contain the invariant-function computations of Section \ref{sec:rank3}, and used it to produce a reorganized and streamlined version of the article, which the authors then revised and checked in full.
Claude was also used to look up references, each of which was verified by the authors. 
The authors take full responsibility for the content of the article.


\section{Preliminaries}\label{section:prelim}

\subsection{Generalities}

Let $(G,X)$ be a geometry, i.e.\ $G$ is a Lie group acting transitively and effectively on $X$ by analytic diffeomorphisms.
Let $V$ be a closed manifold and denote by $\Gamma$ the group of deck transformations of a universal covering space $\widetilde{V} \to V$.
A (marked) %
$(G,X)$ structure on $V$ is a ``dev-hol" pair, up to equivalence, i.e.\ a local %
diffeomorphism $\dev \colon \widetilde{V} \to X$ and a homomorphism $\hol \colon \Gamma \to G$ such that $\dev( \gamma x) = \hol(\gamma) \dev(x)$ for all $\gamma \in \Gamma$ and all $x \in \widetilde{V}$, with two such pairs $(\dev,\hol),(\dev',\hol')$ considered equivalent if there is a diffeomorphism $\phi$ of $V$ isotopic to the identity with lift $\widetilde{\phi} \colon \widetilde{V} \to \widetilde{V}$ and an element $g \in G$ such that $(\dev',\hol')=(g\circ \dev \circ \widetilde{\phi},g  \hol  g^{-1})$.

\begin{example}\label{ex:T2}
  The homogeneous affine structures on $\mathbb{T}^2$ were classified independently by \cite{NY74}, \cite{furness1975locally}. 
  Here the holonomy group lifts to a lattice in the universal cover of one of the following groups: 
  \[
  \begin{array}{c} \begin{pmatrix}  1  & 0 &x \\0 &1 & y \\ 0 & 0 & 1\end{pmatrix} 
  \\
  \mathbb{R}^2 \\ 

  \end{array}
  \begin{array}{c}  \begin{pmatrix}  1&x&y\\ 0&1 &x \\ 0 & 0 & 1\end{pmatrix}\\ 
  \mathbb{R}^2 \\

 \end{array}
  \begin{array}{c} \begin{pmatrix}  1&0&x \\0 &y & 0 \\ 0 & 0 & 1\end{pmatrix} \\
  \{(x,y)\in \R^2: y>0\}  \\

 \end{array}
 \] 
  
  \[\begin{array}{c} 
  \begin{pmatrix} x &y&0 \\0 &x & 0 \\ 0 & 0 & 1 \end{pmatrix}\\ 
  \{(x,y)\in \R^2 : x>0\}\\ 

  \end{array}  
  \begin{array}{c}  \begin{pmatrix} x &0&0 \\0 &y & 0 \\ 0 & 0 & 1 \end{pmatrix}\\
  \{(x,y) \in \R^2 : x,y>0\} \\

\end{array}
  \begin{array}{c} \begin{pmatrix}  e^x \cos(y) & e^x \sin(y) &0 \\-e^x \sin(y) &e^x \cos(y) & 0 \\ 0 & 0 & 1\end{pmatrix}\\ 
  \mathbb{A}^2 \setminus\{0\} \\

  \end{array}\]

Under each group we put the developing image.
See also \cite[Example 1.6]{Ben00}.  
Notice the two half plane structures are projectively equivalent, but affinely inequivalent. 
\end{example}

\subsection{Nilpotent projective manifolds}

We follow Benoist's approach to projective manifolds with nilpotent holonomy
\cite{benoist1994nilvarietes}, \cite{Ben00}. Let $V$ be a compact projective manifold with fundamental group
$\Gamma$, developing map $D\colon\widetilde V\to\Sph^d$ and nilpotent holonomy group
$H=h(\Gamma)$. Choose a maximal connected nilpotent subgroup $N$ of
$G=\mathrm{SL}(d+1,\R)$ such that $H\cap N$ has finite index in $H$, let
$\Gamma'=\{\gamma\in\Gamma:h(\gamma)\in N\}$, and let $\pi\colon\widetilde N\to N$ be the universal
cover. The developing image is then a union of $N$-orbits, among which there are open
ones, and $V$ is cut along the lower-dimensional orbits into elementary pieces, called
\emph{bricks}: compact manifolds with corners, each the quotient of the closure of an
open $\widetilde{N}$-orbit by part of the fundamental group (to be described shortly). 
Up to a finite cover,
$V$ is obtained by gluing finitely many bricks along their faces. The precise statement
is the following.

\begin{theorem}[{Benoist's theorem of bricks \cite[Section 3]{benoist1994nilvarietes}}]
\label{theorem:benoistbricks}
Let $V,D,h,N,\pi$ as above.
   \begin{enumerate}[nosep]
   \item There exists a unique action of $\tilde N$ on $\tilde V$ such that $D(\tilde n \tilde v) = \pi ( \tilde n) D(\tilde v)$. 
   \item The group $N$ has a finite number of orbits in the model space $  \mathbb {S}^d$. 
   \item Let $ \Omega $ be an open orbit of $\tilde N$ in $\tilde V$. Then the restriction of the developing map to the closure $\bar \Omega$ of $\Omega$ is a covering space onto its image. 
   \item The group $\Gamma_0= \{ \gamma \in \Gamma' :\gamma \Omega = 
   \Omega \} $ does not depend on $\Omega$, is normal in $\Gamma' $, and the quotient $\Gamma_0 \backslash \bar \Omega$ is compact. 
   \end{enumerate} 
\end{theorem}

The quotient $\bar \omega = \Gamma_0 \backslash \bar \Omega$ is a projective (affine) manifold with corners.  The manifold $V' = \Gamma ' \backslash \tilde V$, which is a finite covering of $V$, is obtained by gluing a finite number of bricks $\bar \omega$ along their faces. 
If $N$ is simply connected (or, equivalently, has no elliptic elements), the brick $\bar \omega$ is identified with $h (\Gamma_0) \backslash D (\bar \Omega)$. 
Let us now describe the simply connected cases in more detail.

Recall $G=\SL(d+1,\R)$ acts on the projective sphere $\Sph^d=\Sph(\R^{d+1})$, and write
points of $\Sph^d$ as $[x_1:\dots:x_{d+1}]$. For a partition $i_1+\dots+i_k=d+1$, put
\[
  j_\ell=i_1+\dots+i_\ell\qquad(1\le\ell\le k),
\]
so that the $\ell$-th block consists of the coordinates $x_{j_{\ell-1}+1},\dots,x_{j_\ell}$
(with $j_0=0$), and $x_{j_\ell}$ is the last coordinate of that block. Let $U_m$ denote
the group of upper triangular unipotent $m\times m$ matrices, and define
\[
  N_{i_1,\dots,i_k}=\Bigl\{\diag(t_1u_1,\dots,t_ku_k)\ :\ t_\ell>0,\ \
  \prod_{\ell=1}^k t_\ell^{\,i_\ell}=1,\ \ u_\ell\in U_{i_\ell}\Bigr\}\ \subset G .
\]
By \cite{benoist1994nilvarietes}, the groups $N_{i_1,\dots,i_k}$ are maximal connected nilpotent subgroups of
$G$, and every maximal connected nilpotent subgroup of $G$ without elliptic elements is
conjugate to one of them; the others have the form $KAU$ with $K$ a nontrivial compact
torus.

For $1\le j\le d+1$ let $R_j=\diag(1,\dots,1,-1,1,\dots,1)$ be the reflection changing
the sign of the coordinate $x_j$. It has determinant $-1$, so it does not lie in $G$,
but it acts on $\Sph^d$ as an element of $\SL^\pm(d+1,\R)$, the group of all projective
transformations of $\Sph^d$. Let
\[
  W_{i_1,\dots,i_k}=\langle R_{j_1},\dots,R_{j_k}\rangle\cong(\Z/2)^k
\]
be the group generated by the reflections in the last coordinates of the blocks. %
The group $W_{i_1,\dots,i_k}$ permutes the open orbits of $N_{i_1,\dots,i_k}$.
For $w\in W_{i_1,\dots,i_k}$ we call $w\mathcal O$ a \emph{reflected copy} of
$\mathcal O$, and we write, for instance, $R_{j}R_{j'}\mathcal O$ as $R_{jj'} \mathcal O$.

\begin{proposition}[Benoist]\label{prop:open-orbits}
Let $N=N_{i_1,\dots,i_k}$.
\begin{enumerate}
  \item The set
  \[
    \mathcal O=\{[x_1:\dots:x_{d+1}]\in\Sph^d\ :\ x_{j_\ell}>0\ \text{for } 1\le\ell\le k\}
  \]
  is an open $N$-orbit.
  \item The open $N$-orbits are exactly the $2^k$ sets
  \[
    \mathcal O_\varepsilon=\{[x]\in\Sph^d\ :\ \varepsilon_\ell\,x_{j_\ell}>0\ \text{for }
    1\le\ell\le k\},\qquad \varepsilon\in\{\pm1\}^k ,
  \]
  and $W_{i_1,\dots,i_k}$ normalizes $N$ and acts simply transitively on the open $N$-orbits.
  \item The $N$-orbits of codimension $1$ in $\overline{\mathcal O}$ are the following. For each
  $\ell$ with $i_\ell=1$, the face
  \[
    \mathcal O_{j_\ell}=\{[x]\ :\ x_{j_\ell}=0,\ x_{j_m}>0\ \text{for } m\ne\ell\};
  \]
  and for each $\ell$ with $i_\ell\ge2$, the two faces
  \[
    \mathcal O^\pm_{j_\ell}=\{[x]\ :\ x_{j_\ell}=0,\ \pm x_{j_\ell-1}>0,\ x_{j_m}>0\ \text{for }
    m\ne\ell\}.
  \]
  The reflection $R_{j_\ell}$ fixes these faces pointwise and exchanges $\mathcal O$ with the
  reflected copy $R_{j_\ell}\mathcal O$ on the other side of them.
\end{enumerate}
\end{proposition}

\begin{proof}
This is just \cite[ \S2.3, Lemme]{benoist1994nilvarietes} and its proof, translated from lower to upper triangular
matrices. 
There the $N$-orbits are described explicitly: they are the products, over the
blocks, of the orbits of $\{tu: t>0,\ u\in U_m\}$ in $\R^m\setminus\{0\}$ (or $\{0\}$), and the
latter are the sets of vectors whose last nonzero coordinate is in a given position and
has a given sign. Statement (1) and the description of the open orbits in (2) follow
directly. Benoist shows that the group of matrices equal to $\pm\mathrm{Id}$ on each block,
which centralizes $N$, acts simply transitively on the open orbits; since $R_{j_\ell}$ and
the matrix equal to $-\mathrm{Id}$ on the $\ell$-th block have the same effect on the signs of
$x_{j_1},\dots,x_{j_k}$, the same holds for $W_{i_1,\dots,i_k}$, which normalizes $N$ because
its elements are diagonal. For (3), an orbit in $\overline{\mathcal O}$ of codimension $1$ is a
product in which exactly one factor, say the $\ell$-th, lies in an orbit of dimension
$i_\ell-1$ of its block: the vectors whose last nonzero coordinate is $x_{j_\ell-1}$ when
$i_\ell\geq2$, which gives $\mathcal O_{j_\ell}^\pm$ according to its sign, and $\{0\}$ when
$i_\ell=1$, which gives $\mathcal O_{j_\ell}$. The last statement is clear, since $R_{j_\ell}$
changes only the sign of $x_{j_\ell}$.
\end{proof}

Let $V=\T^d$ carry a projective structure, let $\Gamma\cong\Z^d$ be the group of deck
transformations of $\Vt\to V$, let $D\colon\Vt\to\Sph^d$ be the developing map, and let
$H=h(\Gamma)$ be the holonomy group, which is abelian. Fix a maximal connected
nilpotent subgroup $N=N_{i_1,\dots,i_k}$ of $G$ as above, such that
$H\cap N$ has finite index in $H$. We record the properties of $N$ that we use:
\begin{itemize}
  \item $N$ is simply connected, so $\widetilde N=N$ in Theorem \ref{theorem:benoistbricks};
  \item $N$ is the identity component of a real algebraic group, so each $N$-orbit $F$
  in $\Sph^d$ is a semialgebraic set by the Tarski-Seidenberg theorem, and $\overline F\setminus F$ is a union of
  $N$-orbits of dimension less than $\dim F$;
  \item every $N$-orbit in $\Sph^d$ is homeomorphic to a Euclidean space, hence
  contractible. Indeed, by the single-block argument in the proof of
  Proposition \ref{prop:open-orbits}, the orbit of a vector under $\{tu: t>0,\ u\in U_m\}$ is either
  $\{0\}$ or the set of vectors whose last nonzero coordinate is in the same position
  and has the same sign, which is a cell; an $N$-orbit in $\Sph^d$ is the quotient of a
  product of such cells by the free action of $\R_{>0}$.
\end{itemize}
Let $\Gamma'=\{\gamma\in\Gamma: h(\gamma)\in N\}$, a subgroup of finite index. By
Theorem \ref{theorem:benoistbricks}(1), $N$ acts on $\Vt$ with $D(n\cdot x)=n\cdot D(x)$. Let $\Omega\subset\Vt$ be
an open $N$-orbit, and let
\[
  \Gamma_0=\{\gamma\in\Gamma' : \gamma\Omega=\Omega\},
\]
which by Theorem \ref{theorem:benoistbricks}(4) does not depend on $\Omega$.
 
\begin{lemma}\label{lem:brick}
Let $\Omega$ be an open $N$-orbit in $\Vt$. Then:
\begin{enumerate}
  \item $D$ restricts to a diffeomorphism from $\Ob$ onto $D(\Ob)$, which is a union of
  $N$-orbits;
  \item $h$ is injective on $\Gamma_0$, and $h(\Gamma_0)$ is a discrete subgroup of $N$
  acting properly and cocompactly on $D(\Ob)$.
\end{enumerate}
\end{lemma}
 
\begin{proof}
(1) By Theorem \ref{theorem:benoistbricks}(3), $D|_{\Ob}$ is a covering map onto its image. The set $\Ob$ is
$N$-invariant, hence so is $D(\Ob)$. The image $D(\Omega)$ is an open $N$-orbit, which
is contractible; since $D(\Ob)$ is a manifold with corners, it deformation retracts onto
its interior $D(\Omega)$ and is therefore simply connected. So the covering is a
diffeomorphism.
 
(2) Since $D(\gamma x)=h(\gamma)D(x)$ and $D|_{\Ob}$ is injective, $D$ conjugates the
action of $\Gamma_0$ on $\Ob$ to the action of $h(\Gamma_0)$ on $D(\Ob)$. If
$h(\gamma)=1$, then $\gamma$ fixes $\Omega$ pointwise, so $\gamma=1$ because $\gamma$ is a
deck transformation. The action of $\Gamma_0$ on the closed subset $\Ob\subset\Vt$ is
proper, and it is cocompact by Theorem \ref{theorem:benoistbricks}(4); hence so is the action of $h(\Gamma_0)$ on
$D(\Ob)$. Finally, if $h(\Gamma_0)$ were not discrete in $N$, there would be distinct
$\lambda_n\in h(\Gamma_0)$ with $\lambda_n\to1$, and then $\lambda_np\to p$ for
$p\in D(\Omega)$, contradicting properness.
\end{proof}
 
Since $\Gamma\cong\Z^d$, the subgroup $\Gamma_0$ is free abelian, $\Gamma_0\cong\Z^r$;
its \emph{rank} is $r$. By Lemma~\ref{lem:brick}(2), $h(\Gamma_0)\cong\Z^r$ is an abelian
discrete subgroup of $N$.
 
Following \cite[Section 4.1]{Ben00}, let $B\subset N$ be a connected closed subgroup
containing $h(\Gamma_0)$ as a cocompact lattice (it exists by \cite[Ch. II]{raghunathan1972discrete}). As a connected closed subgroup of the
simply connected nilpotent group $N$, the group $B$ is simply connected and nilpotent;
in particular it is contractible, and its only compact subgroup is the trivial one.

Moreover, $B$ is abelian. Indeed, $B$ and $\R^r$ are simply connected nilpotent Lie
groups containing the lattices $h(\Gamma_0)$ and $\Z^r$, and these lattices are
isomorphic since $h(\Gamma_0)\cong\Gamma_0\cong\Z^r$. By Malcev's rigidity theorem
\cite{mal1949class}, see also \cite[Theorem~2.11]{raghunathan1972discrete}, an isomorphism between lattices in simply
connected nilpotent Lie groups extends to an isomorphism of the groups. Hence
$B\cong\R^r$.

\begin{proposition}\label{prop:rank}
When $N$ is simply connected, the following numbers are equal:
\begin{enumerate}
  \item the smallest dimension of an $N$-orbit contained in $D(\Vt)$;
  \item the rank of $\Gamma_0$;
  \item the dimension of $B$.
\end{enumerate}
This number is called the \emph{rank} of the projective structure with respect to $N$.
\end{proposition}
 
\begin{proof}
\emph{(2)$=$(3).}  By the discussion preceding the proposition, $B\cong\R^r$, so $\dim B=r$.

\emph{(3)$=$(1).} Fix an open $N$-orbit $\Omega\subset\Vt$ and put $P=D(\Ob)$.
 
\emph{Step 1: $B$ acts freely, properly and cocompactly on $P$.} The subgroup $B$
preserves each $N$-orbit, hence preserves $P$. Since $h(\Gamma_0)$ acts properly and
cocompactly on $P$ and $B=h(\Gamma_0)\,C$ for a compact set $C\subset B$, the action of
$B$ on $P$ is proper and cocompact. Its stabilizers are compact subgroups of the simply
connected nilpotent group $B$, hence trivial.
 
\emph{Step 2: every $N$-orbit in $P$ has dimension at least $\dim B$.} Each $N$-orbit
$F\subset P$ is $B$-invariant, and by Step~1 every $B$-orbit in $F$ has dimension
$\dim B$.
 
\emph{Step 3: a minimal $N$-orbit in $P$ has dimension exactly $\dim B$.} Let
$F\subset P$ be an $N$-orbit of smallest dimension. Since $\overline F\setminus F$
consists of orbits of smaller dimension, $F$ is closed in $P$. So $F$ is a closed
$B$-invariant subset, and $B$ acts on it freely, properly and cocompactly by Step~1.
The quotient $B\backslash F$ is then a compact manifold without boundary, and $F\to
B\backslash F$ is a principal $B$-bundle. Both $F$ and $B\cong\R^{\dim B}$ are
contractible, so by the long exact sequence of homotopy groups $B\backslash F$ is a
compact contractible manifold without boundary, that is, a point. Hence
$\dim F=\dim B$.
 
\emph{Step 4: the whole developing image.} By Theorem \ref{theorem:benoistbricks}, $V'=\Gamma'\backslash\Vt$ is
obtained by gluing bricks, so $\Vt$ is the union of the closures $\Ob$ of its open
$N$-orbits, and $D(\Vt)$ is the union of the corresponding sets $P$. Since $\Gamma_0$,
and therefore $B$, does not depend on $\Omega$, Steps~1--3 apply to every $P$ with the
same group $B$. So the smallest dimension of an $N$-orbit in $D(\Vt)$ is $\dim B$.
\end{proof}
 
\begin{remark}[Rank in the presence of elliptic elements]\label{remark:elliptics}
The standing assumption that $N$ has no elliptic elements is used in
Lemma~\ref{lem:brick}(1) and in Step~3. When $N$ has a compact factor, $D|_{\Ob}$ need
not be injective, and one must work with a lift of the holonomy to $\widetilde N$
instead of $h(\Gamma_0)$. The Hopf torus with $N=\mathrm{GL}(1,\mathbb C)$  shows that this is not a formality: there
$h(\Gamma_0) = h(\Gamma) \cong\mathbb Z$, while the structure has rank 2 in the sense that the lifted holonomy group is a lattice in $\widetilde{N}$. 
From this perspective it is more natural to view the Hopf torus as a left-invariant structure, see Definition \ref{def:rank4}.

In fact, it is still possible to show that these definitions of rank are equivalent in the presence of elliptics when the rank $r$ equals $d = \dim V$ or $d-1$. 
Here $\Omega$ and hence $\Gamma_0$ is defined by the lifted $\widetilde{N}$ action on $\widetilde{V}$, so the rank of $\Gamma_0$ can be larger than the rank of $h(\Gamma_0)$.
In this case, we note that the rank depends on the choice of $N$. 
For example, returning to the Hopf torus, if we take $N$ to be the group of diagonal matrices or the half-plane group above, the rank is 1, but if we take $N = \GL(1,\C)$, the rank is $2$. 

We always fix $N$, and all statements about rank are with respect to that choice.
\end{remark}

\begin{remark}[Compact factors]\label{rem:compact-factors}
Let $V=\T^d$ carry a projective structure whose holonomy has a finite-index subgroup in a
maximal connected nilpotent group $N=TAU$ with a nontrivial compact torus $T$
\cite[\S2.2]{benoist1994nilvarietes}. 
By \cite[\S4.2]{benoist1994nilvarietes}, after passing to a finite-index subgroup of
$\Gamma$, every holonomy element decomposes as $h(\gamma)=t_\gamma a_\gamma u_\gamma$ with
$t_\gamma\in T$, $a_\gamma\in A$ and $u_\gamma\in U$, and $\gamma\mapsto t_\gamma$ is a homomorphism.

\begin{enumerate}
  \item \emph{If $\{t_\gamma\}$ is finite}, then after passing to a finite cover of $V$ the
  holonomy lies in the connected group $AU$, all of whose elements have positive real
  eigenvalues. Decompose $\R^{d+1}$ into the generalized weight spaces $W_1,\dots,W_k$ of
  $AU$. On each $W_\ell$, the group $AU$ acts by positive scalars times unipotent matrices,
  which are simultaneously upper triangular in a suitable basis. Hence $AU$ is contained
  in a conjugate of the group $N'=N_{\dim W_1,\dots,\dim W_k}$, which is a maximal connected
  nilpotent subgroup without elliptic elements (Section \ref{section:prelim}). The structure, or a finite
  cover of it, is therefore covered by Sections \ref{sec:bricks} and \ref{sec:rank3} with respect to $N'$. Its rank
  with respect to $N'$ may differ from its rank with respect to $N$, however: for the Hopf
  torus with real dilation, $N=\mathrm{GL}(1,\mathbb C)$ gives rank $2$, while the diagonal
  group gives rank $1$ (Remark \ref{remark:elliptics}).
  \item \emph{If $\{t_\gamma\}$ is infinite}, for instance if the holonomy contains a rotation
  by an irrational angle, then the Zariski closure of the holonomy contains a nontrivial
  compact torus, and no group without elliptic elements contains a finite-index subgroup
  of the holonomy. Lemma~ \ref{lem:brick} and Proposition \ref{prop:rank} use the simple connectivity of $N$, and
  Remark \ref{remark:elliptics} shows that this is essential. Such structures of rank $d$ are covered by
  Definition~\ref{def:rank4}, and so are included in Section \ref{sec:max_rannk}. Those of rank $d-1$ are not treated in
  this paper: their study would require working with $\widetilde N$ and a lift of the
  holonomy to it, in place of $h(\Gamma_0)$.
\end{enumerate}
Accordingly, the classification of rank $3$ structures in Section \ref{sec:rank3} is with respect to
maximal connected nilpotent groups without elliptic elements. By (1), it includes every
structure whose holonomy has only finitely many elliptic components, after passing to a
finite cover.
\end{remark}

We now illustrate the combinatorics of words of bricks with two examples.

\begin{example}
    (Goldman tori: rank 1 structures on $\T^2$ with respect to the diagonal group). Let $g=\textrm{diag}(\lambda_1, \lambda_2, \lambda_3)$ be a $3 \times 3 $ diagonal matrix with $\lambda_1 > \lambda_2> \lambda_3 >0$ and eigenvectors $v_1, v_2, v_3$.  Then $g$ acts on an (open) hemisphere of $\mathbb{S}^2$ and the attracting/repelling dynamics of the action cut out a half annulus $A_g$, let $B_g$ be the corresponding half annulus for $g^{-1}$. One may glue together half annuli $A_g, B_g$ in any order to make a closed path, and at the end glue the two boundary curves of the resulting annulus together to obtain a Goldman torus $\mathbb{T}^2_{g,w,h}$ where $g$ is the diagonal matrix, $w$ is a word in half annuli, and $h$ closes up the torus. See \cite{Ben00} Example 1.10. 
\end{example}

The next example of a projective rank 1 structure on $\mathbb{T}^2$  illustrates the dynamics gluing bricks.  Only some words of bricks are permissible.  This issue will be the focus of many of our computations later on, and we resolve it using our idea of an ``invariant function".

\begin{example}[Rank $1$ structures on $\T^2$ for $N_{2,1}$]\label{ex:N21}
Let $N=N_{2,1}\subset\mathrm{SL}(3,\R)$, with open orbit $\mathcal O=\{[x:y:z]:y>0,\ z>0\}$ and faces
\[
  \mathcal O_2^\pm=\{y=0,\ \pm x>0,\ z>0\},\qquad \mathcal O_3=\{z=0,\ y>0\}.
\]
Let $B$ be the one-parameter group generated by $X=\diag(a,a,-2a)+bE_{12}$ with $a,b\neq0$;
then $B$ acts freely on $\mathcal O$ and on all three faces. The function
\[
  f=3a\,\frac xy+b\log\frac zy
\]
is constant exactly on the $B$-orbits in $\mathcal O$, so the orbit space $B\backslash \mathcal O$ is a line,
parametrized by $f$. As a point of $\mathcal O$ approaches a face, $f$ tends to one of the two ends
of this line: to $\sgn(a)\cdot\infty$ at $\mathcal O_2^+$, to $-\sgn(a)\cdot\infty$ at $\mathcal O_2^-$, and to
$-\sgn(b)\cdot\infty$ at $\mathcal O_3$.
 
A \emph{brick} is $\mathcal O$ together with two faces lying at opposite ends; $B$ acts on it
properly with a compact interval as quotient. Two faces at the same end cannot be
used together, since the orbit space would not be Hausdorff. The developing image of a
torus is obtained by gluing bricks over $\mathcal O$ and its reflected copies (changing the signs
of $y$ or $z$) along common faces, and each copy is entered through a face at one end and
left through a face at the other.
 
The two kinds of faces behave differently under this gluing. In the copy across
$\{y=0\}$ the faces $\mathcal O_2^\pm$ exchange ends, because $x/y$ changes sign, whereas the face in
$\{z=0\}$ stays at the same end in every copy, because $\log|z|$ does not change. As a
consequence, after a chain crosses a face in $\{z=0\}$, it is always entered at the end
where these faces lie, and so it can never cross such a face again. A torus has a
periodic chain, so it never crosses $\{z=0\}$ at all: the only brick that occurs is
$\mathcal O\cup \mathcal O_2^+\cup \mathcal O_2^-$, and the tori are obtained by repeatedly crossing $\{y=0\}$.
 
This example illustrates the three questions we have to answer for each group $N$: which faces can be added to $\mathcal O$ to form a brick (determined by the ends of an
invariant function), how the ends change from one copy of $\mathcal O$ to the next, and which
sequences of bricks close up to a torus.
\end{example}

\begin{remark}\label{rem:affine}
 Every affine structure is projective, by identifying $\R^n$ with a hemisphere in $\mathbb{S}^n$.  
 A projective structure is affine if the developing image is in the complement of a hyperplane, and the holonomy preserves the hyperplane.  
 So it suffices to study projective structures for the rest of the paper. 
\end{remark}

\section{Bricks and Their Gluings in Rank \texorpdfstring{$d-1$}{d-1} Structures on \texorpdfstring{$\T^d$}{Td}}\label{sec:bricks}

In this section we consider the case when $\Gamma_0$ has rank $d-1$ and so $B$ is a $(d-1)$-dimensional connected abelian subgroup of a simply connected $N$.  %
As in Section \ref{section:prelim} we assume $N$ has no elliptic elements. Then the developing image contains $N$-orbits of dimension $d-1$, namely the faces, and none of smaller dimension.

By Benoist's brick theorem, $B$ acts properly and cocompactly on a brick $\mathcal{O}'=D(\bar \Omega)$ which is $\mathcal{O}$ union two of its codimension 1 faces. 
In the next  section we first classify the tuples $(N,B,\mathcal{O}')$, %
where the action of $B$ on $\mathcal{O}'$ is proper.
Then, we use the results of this section to understand how to glue together such bordifications into a torus. %
The result is a classification of rank $3$ projective structures on $\mathbb{T}^4$ in Section \ref{sec:rank3}, and clear answers to Benoist's exercises, \cite[Section 4.6]{Ben00}, classifying rank 2 structures on $\T^3$ in Appendix \ref{appendix:T3}.

We sketch the approach common to all cases:
\begin{enumerate}
    \item Enumerate the connected $(d-1)$-dimensional (abelian) subgroups $B$ of $N$ which act \emph{freely and properly} on $\mathcal{O}$ and at least two of its faces.
    We do this by finding a ``normal form" for the Lie algebra $\mathfrak{b}$ of $B$, i.e.\ a preferred basis where it is easy to read off the properness of the action on each face.
    \item We then find an \emph{invariant function} (see Definition \ref{def:invariant-function}),  %
    $f \colon \mathcal{O} \to \mathbb{R}$. 
    By taking limits $f(x_n)$ for $x_n$ limiting to a face of $\mathcal{O}$, we can determine when the action of $B$ on the union of $\mathcal{O}$ and two of its codimension $1$ faces is proper.
    The tuples $(B,\mathcal{O},F_1,F_2)$ where $F_1,F_2$ are codimension $1$ faces of $\mathcal{O}$ and $B$ acts properly on $\mathcal{O} \cup F_1 \cup F_2$, where  the two faces must lie at opposite ends of an invariant function. These $\mathcal{O} \cup F_1 \cup F_2$ are Benoist's ``bricks" which make up the building blocks of a projective rank $d-1$ structure (see Lemma \ref{lem:bricks}).
    \item Using the results of this section (Lemmas \ref{lem:bricks}, \ref{lem:constant-faces}), we determine how the bricks can be glued to form the whole structure. 
    We call this data a ``word of bricks." 
    The possible words of bricks are determined by how the dynamics of the action of $B$ on $\mathcal{O}$ differs from the dynamics of $B$ on an adjacent, reflected copy of $\mathcal{O}$. 

\end{enumerate}

\subsection{Words of bricks}

Throughout this section $N$ is fixed and has no elliptic elements (Section \ref{section:prelim}), and we
consider projective structures on $\T^d$ of rank $d-1$ with respect to $N$. With the
notation of Section \ref{section:prelim}, let $S\subset\Vt$ be the set of points that $D$ maps into an
$N$-orbit of dimension $d-1$; its components are the \emph{faces} in $\Vt$. Label each $\Omega_i$ and each $F_i$ by the $N$-orbit in $\Sph^d$ that contains its
image under $D$. The \emph{label sequence} of the structure is the resulting
bi-infinite sequence, well defined up to shift and reversal; its terms
$(D(\Omega_i),D(F_i),D(F_{i+1}))$ are the bricks.

\begin{lemma}\label{lem:chain}
Each closure $\overline\Omega$ of an open $N$-orbit $\Omega\subset\Vt$ contains exactly
two faces. The open $N$-orbits and the faces of $\Vt$ form a bi-infinite chain
\[
  \cdots\ \Omega_{-1}\ \ F_0\ \ \Omega_0\ \ F_1\ \ \Omega_1\ \ F_2\ \cdots,
  \qquad \overline{\Omega_i}=F_i\cup\Omega_i\cup F_{i+1},
\]
which is unique up to shift and reversal of the indexing.
\end{lemma}

\begin{proof}
By Proposition \ref{prop:rank}, $B$ has dimension $d-1$, and by Step~1 of its proof $B$ acts freely,
properly and cocompactly on $\overline\Omega\cong D(\overline\Omega)$. Since the rank is $d-1$, the developing image contains no $N$-orbits of dimension less than $d-1$ (Proposition \ref{prop:rank}). So $\overline\Omega$ consists of $\Omega$ together with faces, and is a manifold with boundary.
The orbit space
$B\backslash\overline\Omega$ is therefore a compact connected $1$-manifold with nonempty boundary
whose interior $B\backslash\Omega$ is an open interval, so it is a closed interval. Since $B$ acts freely and properly on the face, its orbits are open and closed in it, and the face is connected.
Each
face in $\overline\Omega$ has dimension $d-1$ and is a single closed $B$-orbit, so the
faces in $\overline\Omega$ are exactly the two endpoints.

By Theorem \ref{theorem:benoistbricks}, $\Vt$ is the union of the sets $\overline\Omega$. Each face is a closed
embedded hypersurface of the simply connected manifold $\Vt$, so it separates $\Vt$.
Hence the graph with one vertex for each open $N$-orbit and one edge for each face,
joining the two orbits it bounds, is a tree, and every vertex has degree $2$. It is
infinite because $\Gamma'/\Gamma_0$ acts on it freely, so it is a bi-infinite line.
\end{proof}

Indeed, since $h(\Gamma')\subset N$ preserves every $N$-orbit, while adjacent open orbits carry different labels (they differ by a reflection $R_j$), no element of $\Gamma'$ can reverse the chain. The group
$\Gamma'$ acts on the chain by shifts with kernel $\Gamma_0$, so the label sequence is periodic of some length $k$. 
Moreover $\Gamma'/\Gamma_0$ acts freely by translations, and is therefore isomorphic to $\mathbb Z$.
A \emph{word of bricks} of the structure is the sequence of $k$
consecutive bricks forming one period; it is defined up to cyclic rotation and reversal. A torus with word $w^k$ can be taken to be a $k$-fold cover of one with word $w$; this is why the next proposition is
stated for primitive roots.

\begin{proposition}\label{prop:words-distinct}
Let $\mathcal N$ be the set of $N$-orbits in $\Sph^d$, and let $\Pi_N$ be the group of
permutations of $\mathcal N$ induced by the normalizer of $N$ in $\mathrm{SL}^\pm(d+1,\R)$; it
is finite.
\begin{enumerate}
  \item Let $(\Vt,D)$ and $(\Vt',D')$ be the universal covers of two projective structures
  on $\T^d$ of rank $d-1$ with respect to $N$, with label sequences $L$ and $L'$. Suppose
  that $D'\circ\varphi=g\circ D$ for a diffeomorphism $\varphi\colon\Vt\to\Vt'$ and an element
  $g$ normalizing $N$. Then $L'$ is obtained from $L$ by a shift, possibly a reversal, and
  the permutation of labels induced by $g$.
  \item Consequently, if the structures have words $w$ and $w'$, then their universal
  covers are isomorphic by an isomorphism normalizing $N$ only if the primitive roots of
  $w$ and $w'$ are related by a cyclic rotation, possibly a reversal, and an element of
  $\Pi_N$. Each such class contains only finitely many primitive words. In particular, if
  there are infinitely many primitive admissible words, there are infinitely many
  structures whose universal covers are pairwise non-isomorphic by isomorphisms
  normalizing $N$.
\end{enumerate}
\end{proposition}

\begin{proof}
$\Pi_N$ is a group of permutations of the finite set $\mathcal N$ (Theorem~ \ref{theorem:benoistbricks}(2)), so it is
finite.

(1) Since $g$ normalizes $N$, it maps $N$-orbits to $N$-orbits, preserving their
dimensions. Since $D'\circ\varphi=g\circ D$, a point $x\in\Vt$ develops into an $N$-orbit $X$
if and only if $\varphi(x)$ develops into $gX$. Hence $\varphi$ maps the set $S$ of points
developing into $(d-1)$-dimensional orbits onto $S'$, maps each open $N$-orbit $\Omega_i$
onto some $\Omega'_j$, and maps each face $F_i$ onto some $F'_j$, with the labels changed by
$g$. As a homeomorphism, $\varphi$ preserves adjacency, so it induces an isomorphism between
the chains of Lemma~\ref{lem:chain}. An isomorphism of bi-infinite chains is a shift, possibly
composed with a reversal.

(2) By (1), $\cdots www\cdots$ and $\cdots w'w'w'\cdots$ agree up to shift, reversal and an
element of $\Pi_N$. Since each of these operations preserves the length of a period, the
primitive roots have the same length and are related by a cyclic rotation, possibly a
reversal, and an element of $\Pi_N$. A class therefore consists of words of a single
length and contains at most $2\,\ell\,|\Pi_N|$ primitive words of length $\ell$. Hence
infinitely many primitive words give infinitely many classes.
\end{proof}

\begin{remark}
The normalizer of $N=N_{i_1,\dots,i_k}$ in $\mathrm{SL}^\pm(d+1,\R)$ contains the reflections
$R_j$ and the permutations of blocks of equal size (Lemma~\ref{lem:conventions}(1)), which act nontrivially
on $\mathcal N$. So $\Pi_N$ is in general nontrivial, and some distinct primitive words do give
isomorphic universal covers. We do not know whether every isomorphism of universal covers
of rank $d-1$ structures normalizes $N$. 

\end{remark}

\subsection{Invariant functions, ends and bricks}\label{sec:ends}

This section explains rules for how bricks can be glued together, using the idea of \emph{invariant functions}. Throughout, $\mathcal O\subset\Sph^d$ is an open orbit of $N=N_{i_1,\dots,i_k}$, and $B\subset N$ is a connected
subgroup of dimension $d-1$ acting locally freely on $\mathcal O$, so that every $B$-orbit in $\mathcal O$
is an immersed hypersurface. For $Z\in\mathfrak b$ we write $Z^*$ for the fundamental
vector field on $\mathcal O$, $Z^*_p=\frac{d}{dt}\big|_{t=0}\exp(tZ)\cdot p$.

We represent points of
$\Sph^d$ by unit vectors, so that each coordinate $x_i$ is a continuous function on
$\Sph^d$, and a homogeneous polynomial $Q$ of degree $k$ defines the continuous function
$p\mapsto Q(p)$ on unit vectors. Let $F$ be a codimension-one face of $\mathcal O$, contained in
the hyperplane $\{x_j=0\}$ with $j=j_\ell$ (Proposition \ref{prop:open-orbits}(3)); the coordinate $x_j$
has a constant sign $\epsilon\in\{\pm1\}$ on $\mathcal O$.

\begin{definition}\label{def:invariant-function}
An \emph{invariant function} for $B$ on $\mathcal O$ is a smooth submersion
$f\colon \mathcal O\to\R$ whose fibers $f^{-1}(c)$, for $c\in f(\mathcal O)$, are exactly the $B$-orbits
in $\mathcal O$.
\end{definition}

Equivalently, $f$ is constant on each $B$-orbit, and distinct $B$-orbits lie in distinct
fibers. The image $f(\mathcal O)$ is then an open interval $I\subset\R$, and $f$ induces a
homeomorphism from the orbit space $B\backslash \mathcal O$ onto $I$. Post-composing $f$ with a
diffeomorphism of $I$ onto another interval gives another invariant function, and every
invariant function arises from a given one in this way. In particular, when
$I\neq\R$ we may reparametrize so that $f(\mathcal O)=\R$.

\begin{lemma}\label{lem:invariant-function}
Let $Z_1,\dots,Z_{d-1}$ be a basis of $\mathfrak b$, and let $f\colon \mathcal O\to\R$ be a
smooth submersion with connected fibers such that $Z_i^*(f)=0$ for all $i$. Then $f$ is
an invariant function for $B$.
\end{lemma}

\begin{proof}
Since $Z_i^*(f)=0$ for a basis of $\mathfrak b$, $f$ is constant along every curve
$t\mapsto\exp(tZ)\cdot p$; as $B$ is connected, it is generated by such one-parameter
subgroups, so $f$ is constant on $B$-orbits. Hence each fiber $F=f^{-1}(c)$ is a union of
$B$-orbits. Since $f$ is a submersion, $F$ is an embedded hypersurface, and each
$B$-orbit in $F$ has dimension $d-1=\dim F$, so it is open in $F$. The orbits
therefore partition $F$ into disjoint open subsets, and since $F$ is connected, $F$ is a
single orbit.
\end{proof}

\begin{remark}\label{rem:limits}
Let $F$ be a codimension-one face of $\mathcal O$ (an $N$-orbit in $\partial \mathcal{O}$). We say that $f$
\emph{tends to $+\infty$ (resp.\ $-\infty$) at $F$} if $f(p_n)\to+\infty$ (resp.\
$-\infty$) for every sequence $p_n\in \mathcal O$ converging to a point of $F$. Since $f$ is
$B$-invariant, it suffices to check this along one sequence converging to each point of
$F$, and when $B$ acts on $F$ with dense orbits, or transitively, along a single
sequence. These limits determine the bricks (and how they can be glued), see  Example \ref{ex:N21}, Section \ref{sec:rank3}, and Appendix \ref{appendix:T3}.
\end{remark}

\begin{lemma}[Limits at a face]\label{lem:limits}
\begin{enumerate}
  \item \emph{(Pole type.)} Suppose that near $F$ we can write
  \[
    f=\frac{Q}{x_j^{\,k}}+R ,
  \]
  where $k\ge1$, $Q$ is a homogeneous polynomial of degree $k$ which does not vanish on
  $F$, and $x_j^{\,k}R\to0$ at every point of $F$. Then $f$ tends to
  $\epsilon^k\operatorname{sgn}(Q|_F)\cdot\infty$ at $F$.
  \item \emph{(Logarithmic type.)} Suppose that near $F$ we can write
  $f=c\log|x_j|+R$ with $c\neq0$ and $R$ bounded near each point of $F$. Then $f$ tends
  to $-\operatorname{sgn}(c)\cdot\infty$ at $F$.
\end{enumerate}
\end{lemma}

\begin{proof}
The face $F$ is connected (it is a cell), so $Q$ has a constant sign on it. Let
$p_n\in \mathcal O$ converge to $q\in F$. Then $x_j(p_n)\to0$ with sign $\epsilon$ and
$Q(p_n)\to Q(q)\neq0$, so $Q(p_n)/x_j(p_n)^k\to\epsilon^k\operatorname{sgn}(Q(q))\cdot\infty$,
while $R(p_n)$ is negligible compared with this term.  This proves (1); (2) follows in
the same way from $\log|x_j(p_n)|\to-\infty$, since $R$ absorbs the other logarithms which are bounded near $F$. 
\end{proof}

\begin{corollary}[Parity criterion]\label{cor:parity}
Let $i_\ell\ge2$, $j=j_\ell$, and let $\mathcal O^\pm_j=\{x_j=0,\ \pm x_{j-1}>0,\dots\}$ be the two
faces of $\mathcal O$ in $\{x_j=0\}$. Suppose $f$ is of pole type at both faces with the same
$Q$ and $k$, and that
\[
  Q|_{\{x_j=0\}}=x_{j-1}^{\,k}\,u
\]
for a function $u$ that is positive on $\mathcal O^+_j\cup \mathcal O^-_j$. Then $\mathcal O^+_j$ and $\mathcal O^-_j$ lie at
the same end of $\R$ if $k$ is even, and at opposite ends if $k$ is odd.
\end{corollary}

\begin{proof}
By Lemma~\ref{lem:limits}(1), $f$ tends to $\epsilon^k(\pm1)^k\cdot\infty$ at $\mathcal O^\pm_j$.
\end{proof}

 The bricks and the way they
are glued are governed by the following rules, which we use throughout this section.

\begin{lemma}[Bricks and gluing]\label{lem:bricks}
Let $B$ act freely on $\mathcal{O}$ with surjective invariant function $f \colon \mathcal{O} \to \R$.
\begin{enumerate}

   \item Let $F_1,F_2$ be faces of $\mathcal O$ on which $B$ acts freely. Then $B$ acts properly on
  $P=\mathcal O\cup F_1\cup F_2$ if and only if $f$ tends to opposite ends of $\R$ at $F_1$ and
  $F_2$. In that case $B\backslash P$ is a closed interval with endpoints $F_1$ and $F_2$, and
  $P$ is a brick for $B$.
  \item In the chain of Lemma \ref{lem:chain}, each $\Omega_i$ is entered through $F_i$ and left
  through $F_{i+1}$, and these two faces lie at opposite ends for an invariant function
  on the copy $D(\Omega_i)$.
  \item \emph{(Periodicity.)}  Let $F$ be a face label (an $N$-orbit in $\mathbb S^d$). Suppose that $F$ is crossed at most once in every chain satisfying (2) in which the copy across a face in $\{x_j=0\}$ is $R_j$ applied to the current copy. Then $F$ is not crossed in the chain of any torus. 
\end{enumerate}
\end{lemma}

\begin{proof}

(1) Recall that $B$ acts simply transitively on each $F_i$ (Remark \ref{rem:free-proper-faces}), and that near a
point $q$ of a face $F\subset\{x_j=0\}$ the set $\mathcal O\cup F$ is a closed half-ball: the signs
of the coordinates that determine $\mathcal O$ are constant near $q$, except for that of $x_j$,
which vanishes on $F$ (Proposition~\ref{prop:open-orbits}). In particular $P$ is a manifold with boundary
$F_1\cup F_2$.

\emph{Step 1: local charts.} Let $q\in P$. Choose a smooth curve $\gamma$ through
$q=\gamma(0)$, defined on $J=(-\varepsilon,\varepsilon)$ if $q\in\mathcal O$ and on $J=[0,\varepsilon)$ if
$q\in F_i$, as follows. If $q\in\mathcal O$, take $\gamma$ with $(f\circ\gamma)'(0)\neq0$, which is
possible since $f$ is a submersion. If $q\in F_i$, take $\gamma$ transverse to the hyperplane
of $F_i$, with $\gamma(s)\in\mathcal O$ for $s>0$. Define
\[
  \Psi\colon B\times J\to P,\qquad \Psi(b,s)=b\,\gamma(s).
\]
At $(e,0)$ the differential of $\Psi$ is invertible. On $B$ it is the differential of the
orbit map, which is injective with image the tangent space of the orbit $B\cdot q$: this
orbit is a level set of $f$ if $q\in\mathcal O$, and equals $F_i$ if $q\in F_i$. On $J$ it is
$\gamma'(0)$, which is transverse to that orbit. By the inverse function theorem (for
manifolds with boundary when $q\in F_i$), after shrinking $\varepsilon$ there is an open
neighborhood $B_0$ of $e$ such that $\Psi$ restricts to a diffeomorphism from
$B_0\times J$ onto an open neighborhood of $q$ in $P$.

Let $\varphi(s)=f(\gamma(s))$ for $s\in J$ with $\gamma(s)\in\mathcal O$. Since $f$ is
$B$-invariant, $f(\Psi(b,s))=\varphi(s)$. Since $f$ is a submersion on $\mathcal O$ and $\Psi$ is a
local diffeomorphism, $d(f\circ\Psi)$ never vanishes there, and as $f\circ\Psi$ does not
depend on $b$, we get $\varphi'(s)\neq0$. Hence $\varphi$ is strictly monotone: on $J$ if
$q\in\mathcal O$, and on $(0,\varepsilon)$ if $q\in F_i$, where moreover $\varphi(s)\to\pm\infty$ as
$s\to0$ by assumption.

\emph{Step 2: the extended function.} Suppose that $f$ tends to $-\infty$ at $F_1$ and to
$+\infty$ at $F_2$. Define $\bar f\colon P\to[-\infty,+\infty]$ by $\bar f=f$ on $\mathcal O$,
$\bar f=-\infty$ on $F_1$ and $\bar f=+\infty$ on $F_2$. By the definition of ``tends to
$\pm\infty$ at $F_i$'' (Remark \ref{rem:limits}), $\bar f$ is continuous. It is $B$-invariant, and its
fibers are exactly the $B$-orbits in $P$: in $\mathcal O$ by Definition \ref{def:invariant-function}, and on the faces
because each $F_i$ is a single orbit. Since $f$ is unbounded in both directions, $f(\mathcal O)=\R$
and $\bar f$ is onto. The induced map
\[
  \hat f\colon B\backslash P\to[-\infty,+\infty]
\]
is therefore a continuous bijection. It is also open: by Step~1, near every point $\bar f$
is given in the chart $\Psi$ by the strictly monotone function $\varphi$ (extended by
$\varphi(0)=\pm\infty$ at a face), which is an open map onto its image. Hence $\hat f$ is a
homeomorphism, so $B\backslash P$ is Hausdorff and is a closed interval with endpoints $F_1$
and $F_2$.

\emph{Step 3: properness.} First we show that every $q\in P$ has a neighborhood $U$ such
that $\{b\in B:bU\cap U\neq\emptyset\}$ is relatively compact. Take $\Psi$, $B_0$ and $J$ as in
Step~1, choose a compact neighborhood $K\subset B_0$ of $e$, and put $U=\Psi(\mathrm{int}K\times J)$.
If $b\,b_1\gamma(s_1)=b_2\gamma(s_2)$ with $b_1,b_2\in K$, then applying $\bar f$ gives
$\varphi(s_1)=\varphi(s_2)$, with $\varphi(0)=\pm\infty$ in the case of a face. Since $\varphi$ is
injective, $s_1=s_2$, so $b_2^{-1}bb_1$ fixes $\gamma(s_1)$. Since $B$ acts freely on $P$, this
gives $b=b_2b_1^{-1}\in KK^{-1}$, which is compact.

Now let $p_n\to p$ and $b_np_n\to p'$ in $P$. The images of $p_n$ and $b_np_n$ in
$B\backslash P$ coincide and converge to the images of $p$ and $p'$; since $B\backslash P$ is
Hausdorff, $p'=bp$ for some $b\in B$. Then $b^{-1}b_np_n\to p$ and $p_n\to p$, so for a
neighborhood $U$ of $p$ as above, eventually $p_n\in U$ and $b^{-1}b_np_n\in U$. Hence
$b^{-1}b_n$ lies in a compact set, and $(b_n)$ has a convergent subsequence. So the action is
proper, and $P$ is a brick.

\emph{Step 4: the converse.} Suppose that $f$ tends to $+\infty$ at both $F_1$ and $F_2$
(the case $-\infty$ is the same). Let $q_i\in F_i$, and let $V_1,V_2$ be $B$-invariant
neighborhoods of $q_1,q_2$ in $P$. By Step~1, $V_i$ contains a set $\Psi_i(B_0\times[0,\varepsilon))$,
on which $\bar f$ takes all values in $(c_i,+\infty]$ for some $c_i$. For
$c>\max(c_1,c_2)$, the orbit $f^{-1}(c)$ therefore meets both $V_1$ and $V_2$. So the
distinct orbits $F_1$ and $F_2$ cannot be separated by $B$-invariant open sets, and
$B\backslash P$ is not Hausdorff. A proper action has Hausdorff orbit space, so the action is
not proper.

(2) is Lemma \ref{lem:chain}, since each
$\overline{\Omega_i}$ is a brick. (3) By Lemma~3.1 the label sequence of a torus is
periodic, so a face crossed once is crossed infinitely often.
\end{proof}

\begin{remark}\label{rem:free-proper-faces}
If $B$ acts freely on a face $F$, it acts properly on it. Indeed, $\dim B=d-1=\dim F$, so
each $B$-orbit in $F$ is open; since $F$ is connected, $B$ acts simply transitively on
$F$, so the orbit map $B\to F$ is a diffeomorphism and the action is proper.
\end{remark}

\begin{remark}\label{rem:copies}
On a reflected copy $w \mathcal O$, an invariant function for $B$ is given by the same formula as on
$\mathcal O$, with each $\log x_i$ replaced by $\log|x_i|$: the equations $Z^*(f)=0$ are local, and
the formula solves them wherever it is defined. Note that $f\circ w^{-1}$ is in general
\emph{not} an invariant function for $B$ on $w \mathcal O$, since it is invariant under $wBw^{-1}$,
which may differ from $B$. Since only the relative signs within one copy matter
(Lemma~\ref{lem:bricks}(2)), Lemma~\ref{lem:limits} applies to each copy separately, with
$\epsilon$ the sign of $x_j$ on that copy.
\end{remark}

\begin{definition}\label{def:constant-alternating}
Fix, on every copy $w\mathcal O$, the invariant function $f_w$ given by 
Remark~\ref{rem:copies}. A face $F'$ is shared by exactly two copies $w \mathcal O$ and
$R_jw \mathcal O$. We call $F'$ \emph{constant} if $f_w$ and $f_{R_jw}$ tend to the same end
at $F'$, and \emph{alternating} if they tend to opposite ends.
\end{definition}

 By Lemma \ref{lem:limits} and Corollary \ref{cor:parity}, logarithmic faces and even-order pole pairs are constant, and odd-order pole pairs are alternating.

\begin{lemma}[Constant faces]\label{lem:constant-faces}
Suppose that, with respect to the functions $f_w$ of Definition \ref{def:constant-alternating}, all constant faces on which $B$ acts freely and properly lie at the same end $c$ in every copy. Then no constant face is crossed in the chain of a torus.
\end{lemma}

\begin{proof}
For each copy $\Omega_i$ in the chain, let $s_i$ be the end at which its entry face
$F_i$ lies. By Lemma~\ref{lem:bricks}(2), the exit face $F_{i+1}$ lies at the end $-s_i$.
If $F_{i+1}$ is alternating, its end in the next copy is $-(-s_i)=s_i$, so
$s_{i+1}=s_i$. If $F_{i+1}$ is constant, its end in the next copy is again $-s_i$, so
$s_{i+1}=-s_i$; moreover $-s_i=c$, since all constant faces lie at $c$. Thus crossing a
constant face forces $s_{i+1}=c$. After that the entry end stays $c$ under alternating
crossings, while crossing a constant face would require an entry end $-c$. So a
constant face is crossed at most once, and by Lemma~\ref{lem:bricks}(3) not at all in a
torus.
\end{proof}

\begin{remark}\label{rem:const-opp} If there are two constant faces which lie at opposite ends, then we can alternate between them, and still use them in words of bricks.   For example this happens in $N_{2,1,1,1}$ (Proposition \ref{prop:N2111}) and $N_{3,1,1}$ (Proposition \ref{prop:N311}). 
\end{remark}

\subsection{Realizing words of bricks}\label{sec:realizing}
 
As in the rest of this section, $N=N_{i_1,\dots,i_k}$ is fixed, and
$B\subset N$ is a connected abelian subgroup of dimension $d-1$.
 
\begin{definition}\label{def:admissible-word}
An \emph{admissible periodic word} for $B$ is a bi-infinite sequence
$(C_i,F_i)_{i\in\Z}$, where each $C_i$ is an open $N$-orbit (a copy $w\mathcal O$) and each $F_i$ is
a codimension-one face of both $C_{i-1}$ and $C_i$, such that:
\begin{enumerate}
  \item $B$ acts freely on each $C_i$ and on each $F_i$;
  \item for each $i$, there is an invariant function for $B$ on $C_i$ that tends to opposite ends of
  $\R$ at $F_i$ and at $F_{i+1}$;
  \item if $F_{i+1}\subset\{x_j=0\}$, then $C_{i+1}=R_jC_i$;
  \item there is $k\ge1$ with $(C_{i+k},F_{i+k})=(C_i,F_i)$ for all $i$.
\end{enumerate}
The sets $P_i=F_i\cup C_i\cup F_{i+1}\subset\Sph^d$ are its \emph{bricks}; by (1), (2),
Lemma \ref{lem:bricks}(1), and Remark \ref{rem:free-proper-faces}, $B$ acts properly on each $P_i$, and $B\backslash P_i$ is a closed interval
with endpoints $F_i$ and $F_{i+1}$.
\end{definition}
 
By Lemma \ref{lem:chain} and Lemma \ref{lem:bricks}(2), the chain of a projective structure on $\T^d$ of rank
$d-1$ with respect to $N$, with group $B$, is an admissible periodic word. The following
lemma gives the converse.
 
\begin{lemma}\label{lem:realize}
Let $(C_i,F_i)_{i\in\Z}$ be an admissible periodic word for $B$ with period $k$, let
$\Gamma_0\subset B$ be a lattice, and let $g$ be an element of the centralizer of $B$ in
$N$. Then there is a projective structure on $V=\T^d$ of rank $d-1$ with respect to $N$,
with group $B$, whose chain is the given word, and whose holonomy group is the group
generated by $\Gamma_0$ and $g$.
\end{lemma}
 
\begin{proof}
\emph{The universal cover.} Let $\Vt$ be the space obtained from the disjoint union of the
bricks $P_i$ by identifying, for every $i$, the face $F_{i+1}\subset P_i$ with the face
$F_{i+1}\subset P_{i+1}$ by the identity of $\Sph^d$. Let $D\colon\Vt\to\Sph^d$ be the map
given by the inclusion on each $P_i$.
 
\emph{$\Vt$ is a manifold and $D$ is a local diffeomorphism.} Let $q\in F_{i+1}$, and let
$F_{i+1}\subset\{x_j=0\}$. An open $N$-orbit is determined by the signs of the
coordinates $x_{j_1},\dots,x_{j_k}$ (Proposition \ref{prop:open-orbits}), and these signs are constant near
$q$ except for $x_j$, which vanishes at $q$. If the block of $x_j$ has size at least two,
the sign of $x_{j-1}$ is also constant near $q$. So a small neighborhood of $q$ in $\Sph^d$
consists of points of $C_i$ (one sign of $x_j$), of $C_{i+1}=R_jC_i$ (the other sign) and
of $F_{i+1}$ (where $x_j=0$). Hence $D$ maps a neighborhood of $q$ in $\Vt$
homeomorphically onto a neighborhood of $q$ in $\Sph^d$. Away from the faces this is
clear, so $\Vt$ is a manifold without boundary and $D$ is a local diffeomorphism.
 
\emph{The action of $B$.} Since $B$ preserves each $C_i$ and each $F_i$, it acts on $\Vt$,
and $D$ is $B$-equivariant. The action is free by (1). It is proper: the bricks form a
locally finite cover of $\Vt$ by closed $B$-invariant subsets, so a compact
$K\subset\Vt$ meets finitely many of them, and if $bK\cap K\neq\emptyset$, then
$b(K\cap P_i)\cap(K\cap P_i)\neq\emptyset$ for one of these $i$, which confines $b$ to a
compact set by properness on $P_i$. The orbit space $B\backslash\Vt$ is the union of the
closed intervals $B\backslash P_i$ glued end to end, so it is homeomorphic to $\R$.
Hence $\Vt\to\R$ is a principal $B$-bundle, which is trivial since $B\cong\R^{d-1}$ is
contractible, and $\Vt\cong B\times\R$ is simply connected.
 
\emph{The deck group.} By (4), $P_{i+k}=P_i$ as subsets of $\Sph^d$, so there is a
diffeomorphism $\tau$ of $\Vt$ mapping $P_i$ to $P_{i+k}$ by the identity of $\Sph^d$. The
element $g$ commutes with $B$ and preserves every $N$-orbit, so it preserves each $P_i$
and acts on $\Vt$, commuting with $B$ and with $\tau$. Put $\tau_g=\tau\circ g$. Then
$D\circ\tau_g=g\circ D$, and $\tau_g$ commutes with $B$. On $B\backslash\Vt\cong\R$,
the map induced by $g$ preserves each interval $B\backslash P_i$ and its endpoints, so
the map induced by $\tau_g$ sends each interval to the one $k$ steps further. It
therefore acts freely, properly and cocompactly on $\R$, with quotient a circle.
 
Let $\Gamma=\Gamma_0\times\langle\tau_g\rangle\cong\Z^d$, where $\Gamma_0$ acts through
$B$. The quotient $\langle\tau_g\rangle\backslash\Vt$ is a principal $B$-bundle over a
circle, hence trivial, so $\Gamma\backslash\Vt\cong(\Gamma_0\backslash B)\times S^1$ is a
torus $\T^d$, and $\Gamma$ acts freely and properly on $\Vt$.
 
\emph{The structure.} Define $h\colon\Gamma\to N$ by $h|_{\Gamma_0}=\mathrm{id}$ and
$h(\tau_g)=g$. Then $D(\gamma x)=h(\gamma)D(x)$ for all $\gamma\in\Gamma$, so $(D,h)$ is
a developing pair for a projective structure on the torus $V=\Gamma\backslash\Vt$, whose
universal cover is $\Vt$, with holonomy
group $\langle\Gamma_0,g\rangle$. Its developing image is the union of the $C_i$ and the
$F_i$, so the smallest dimension of an $N$-orbit in it is $d-1$, and the structure has
rank $d-1$ with respect to $N$. By construction, its chain is the given word.
\end{proof}
 
\begin{remark}
Taking $g=1$ gives holonomy group $\Gamma_0\cong\Z^{d-1}$; the holonomy is then not
injective, and the developing image need not be simply connected. Taking $g\in B$ only
changes the identification of the deck group. The words $w$ and $w^m$ give structures
with the same universal cover; the torus defined by $w^m$ with $g^m$ covers the torus
defined by $w$ with $g$.
\end{remark}

\section{Rank 3 Structures on \texorpdfstring{$\T ^4$}{T4}}\label{sec:rank3}

In this section we classify the projective structures on $\T^4$ of rank $3$. Let $V=\T^4$
carry such a structure, with respect to a maximal connected nilpotent subgroup $N$ of
$\mathrm{SL}(5,\R)$ without elliptic elements (Section \ref{section:prelim};  see Remark \ref{rem:compact-factors} for the case of compact factors), and let $B\subset N$ be its
$3$-dimensional abelian group, with Lie algebra $\mathfrak b$ (Proposition \ref{prop:rank}). Up to
conjugation, $N$ is one of the seven groups $N_{i_1,\dots,i_k}$ with $i_1\ge\dots\ge i_k$ and
$i_1+\dots+i_k=5$ \cite{Ben00} Section~4.1. The results are summarized as follows.
 
\begin{theorem}\label{thm:rank3-summary}
Let $V=\T^4$ carry a projective structure of rank $3$ with respect to $N$. Then, up to
conjugation, $N$ and $\mathfrak b$ are as in one of the propositions listed below, and the
bricks occurring in the chain of $\Vt$ are as indicated. (Other bricks may exist over $\mathcal O$,
but by Lemma \ref{lem:constant-faces} they are never used.)
\begin{center}
\renewcommand{\arraystretch}{1.25}
\begin{tabular}{lll}
$N$ & reference & bricks occurring in the chain of a torus \\ \hline
$N_5$ & Prop.~\ref{prop:N5} & no structures \\
$N_{4,1}$ & Prop.~\ref{prop:N41} & one brick $\mathcal O\cup \mathcal O_4^\pm$ \\
$N_{3,2}$ & Prop.~\ref{prop:N32} & one brick $\mathcal O\cup \mathcal O_5^\pm$ \\
$N_{3,1,1}$ & Prop.~\ref{prop:N311} & case (1a): $2$, $3$ or $4$ bricks; cases (1b), (1c), (2): one brick \\
$N_{2,2,1}$ & Prop.~\ref{prop:N221} & case (1): $4$ bricks; case (2): one brick \\
$N_{2,1,1,1}$ & Prop.~\ref{prop:N2111} & case (1): $1$, $4$ or $6$ bricks; case (2a): $1$ or $2$ bricks \\
$N_{1,1,1,1,1}$ & Prop.~\ref{prop:diagonal} & $n_+n_-$ bricks, $n_\pm\ge1$, $n_++n_-\le5$
\end{tabular}
\end{center}
Conversely, every Lie algebra $\mathfrak b$ and every admissible periodic word occurring in
these propositions, together with a lattice $\Gamma_0\subset B$ and an element $g$ of the
centralizer of $B$ in $N$, defines a projective structure on $\T^4$ of rank $3$, with
holonomy group $\langle\Gamma_0,g\rangle$.
\end{theorem}
 
\begin{corollary}\label{cor:rank3-dichotomy}
In the situation of Theorem~\ref{thm:rank3-summary}, exactly one of the following holds.
\begin{enumerate}
  \item There is a single brick. Then the words of bricks are its powers, the chain of
  $\Vt$ has period $2$ or $4$, and the structure is affine.
  \item There are at least two bricks. Then there are 
  infinitely many structures whose universal covers are pairwise non-isomorphic by isomorphisms normalizing $N$ (Proposition \ref{prop:words-distinct}).
\end{enumerate}
\end{corollary}
 
\begin{proof}
This is read off from the propositions cited in the table: in every case with a single
brick the structure is shown to be affine, since one hyperplane is never crossed, so Lemma \ref{lem:conventions}(5) applies. And in every case with at least two bricks the
chain can choose between at least two faces at some end infinitely often.
\end{proof}

\subsection*{The general outline} \label{sec:conventions}
 
The subsections below treat the seven groups in turn. To avoid repetition, we collect
here the steps that are common to all of them. Points of $\Sph^4$ are written
$[v:w:x:y:z]$, $E_{ij}$ denotes an elementary matrix, $\mathcal O$ is the open orbit of
Proposition \ref{prop:open-orbits}, and its faces are labelled by the vanishing coordinate as in
Proposition \ref{prop:open-orbits}(3).
 
\begin{lemma}\label{lem:conventions}
\begin{enumerate}
  \item \emph{(Position.)} After conjugating the structure by an element of
  $W_{i_1,\dots,i_k}$, which normalizes $N$, we may assume that $\mathcal O$ is one of the copies in
  the chain of $\Vt$. Conjugating by a permutation of blocks of equal size, which also
  normalizes $N$, we may permute those blocks.
  \item \emph{(Freeness.)} For $p\in\R^5\setminus\{0\}$, the group $B$ acts freely at $[p]$ if
  and only if $b_1p,b_2p,b_3p,p$ are linearly independent for a basis $b_1,b_2,b_3$ of
  $\mathfrak b$, that is, if and only if some $4\times4$ minor of the matrix
  $(b_1p,b_2p,b_3p,p)$ is nonzero. A nilpotent $X\in\mathfrak b$ fixes $[p]$ if and only if
  $Xp=0$. If $B$ acts freely on a face, it acts properly on it (Remark \ref{rem:free-proper-faces}).
  \item \emph{(Invariant functions.)} Let $f$ be a function on $\mathcal O$ annihilated by the
  fundamental vector fields of a basis of $\mathfrak b$. If $f=\alpha\,x_i/x_m+g$
  with $\alpha\neq0$ constant and $g$ independent of $x_i$ %
  or if $f$ is a nonzero linear combination of
  logarithms of coordinate ratios, then the fibers of $f$ are graphs over the remaining variables, hence connected, and
  $f$ is an invariant function (Lemma \ref{lem:invariant-function}).
  \item \emph{(Ends and gluing.)} At a face in $\{x_j=0\}$, an invariant function has
  pole type or logarithmic type (Lemma \ref{lem:limits}), and Corollary \ref{cor:parity} and Definition \ref{def:constant-alternating}
  determine the ends and whether the face is alternating (odd pole) or constant (even
  pole or logarithm). Along the chain of a torus, each copy is entered at one end and
  left at the other; crossing an alternating face preserves the end at which the next
  copy is entered, and crossing a constant face reverses it (proof of Lemma \ref{lem:constant-faces}). If
  all constant faces on which $B$ acts freely lie at one end, they are never crossed
  (Lemma \ref{lem:constant-faces}).
  \item \emph{(Affine structures.)} Let $j$ be the last coordinate of a block. If the
  word of a structure never crosses a face in $\{x_j=0\}$, then the developing image lies
  in $\{x_j>0\}$, and the structure is affine.
  \item \emph{(Existence.)} Every admissible periodic word of bricks, together with a
  lattice $\Gamma_0\subset B$ and an element $g$ of the centralizer of $B$ in $N$, is
  realized by a structure with holonomy group $\langle\Gamma_0,g\rangle$ (Lemma \ref{lem:realize}).
\end{enumerate}
\end{lemma}
 
\begin{proof}
(1) The reflections $R_j$ and the permutations of equal blocks are diagonal or permutation
matrices, so conjugation by them preserves $N$, and $W_{i_1,\dots,i_k}$ acts transitively on
the open orbits (Proposition \ref{prop:open-orbits}(2)).
 
(2) The differential at the identity of the orbit map $h\mapsto h\cdot[p]$ is
$X\mapsto Xp\bmod\R p$, so the stabilizer of $[p]$ in $B$ has trivial Lie algebra exactly
under the stated condition. It is then trivial: if $\exp(X)p=\lambda p$ with
$X\in\mathfrak b$, write $X=S+U$ with $S$ diagonalizable over $\R$ and $U$ nilpotent,
commuting. Then $p$ is an eigenvector of $\exp S$, hence of $S$, with eigenvalue
$\mu=\log\lambda$ (the eigenvalues of elements of $N$ are positive), and $\exp(U)p=p$, so
$Up=0$; hence $Xp=\mu p$, and $X$ lies in the Lie algebra of the stabilizer. For nilpotent
$X$, the only eigenvalue is $0$, so $Xp\in\R p$ forces $Xp=0$.
 
(3) In the first case, solving $f=c$ for $x_i$ expresses each fiber as a graph over the
other coordinates of $\mathcal O$; in the second, each fiber is the preimage of an affine
hyperplane under a submersion with connected fibers onto a vector space of
logarithms. Since $\mathcal O$ is a cell, the fibers are connected.
 
(4) is a restatement of Lemma \ref{lem:limits}, Corollary \ref{cor:parity}, Definition \ref{def:constant-alternating} and Lemma \ref{lem:constant-faces}, together
with Lemma \ref{lem:bricks}(2).
 
(5) Crossing a face in $\{x_{j'}=0\}$ replaces the current copy by its image under $R_{j'}$,
which changes only the sign of $x_{j'}$. If $j'\ne j$ throughout, every copy, and every face
between them, lies in $\{x_j>0\}$. The hyperplane $\{x_j=0\}$ is $N$-invariant, and the
holonomy $\langle\Gamma_0,g\rangle$ lies in $N$, so the structure is affine by
Remark \ref{rem:affine}.
 
(6) is Lemma \ref{lem:realize}.
\end{proof}
 
In each subsection we therefore only state: the normal forms of $\mathfrak b$, the
freeness conditions (read off from the minors, using (2)), an invariant function (with
the reason for (3)), the ends of the faces (using (4)), and the resulting bricks and words.
Existence is then given by (6), and affineness, where it holds, by (5); we do not repeat
these. A SymPy script verifying all symbolic computations in Sections \ref{sec:rank3}, \ref{sec:max_rannk} and Appendices \ref{appendix:T3} and \ref{appendix:nonunipotent} is available as supplementary material.

\subsection{The unipotent case, \texorpdfstring{$N_5$}{N5}}\label{section:N5}

By \cite{Ben00} Exercise 4.6(c),(d), a compact projective 4-manifold with unipotent holonomy that is not complete in the sense of \cite{Ben00} has rank $p=4$; in particular there is no rank 3 structure with respect to $N_5$.  Complete means $p=0$, which is impossible for a torus, since then $\widetilde V$ would be $\mathbb S^4$. Since this case of the exercise is left to the reader in \cite{Ben00}, we include a proof.

\begin{proposition}\label{prop:N5}
There is no projective structure on $\T^4$ of rank $3$ with respect to $N_5$.
\end{proposition}

\begin{proof}
Since $t^5=1$ and $t>0$ force $t=1$, the group $N_5$ is the group $U_5$ of upper
triangular unipotent matrices, with Lie algebra $\mathfrak u_5$. By
Proposition \ref{prop:open-orbits}, applied with a single block of size $5$, the open $N_5$-orbits are
$\mathcal O_\pm=\{\pm z>0\}$, and the codimension-one faces of $\mathcal O= \mathcal O_+$ are
\[
  \mathcal O_5^\pm=\{[v:w:x:y:0]\ :\ \pm y>0\}.
\]

Suppose there is a rank $3$ structure, and let $B\subset N_5$ be the group of
Proposition \ref{prop:rank}, which is $3$-dimensional and abelian, with Lie algebra
$\mathfrak b\subset\mathfrak u_5$. %

By Lemma \ref{lem:chain}, the closure of an open $N_5$-orbit in $\widetilde V$ develops onto $D(\Omega)\cup F_1\cup F_2$, where $D(\Omega)$ is $\mathcal O_+$ or $\mathcal O_-$. Conjugating by $R_5$ if necessary (it normalizes $N_5$ and exchanges $\mathcal O_\pm$), we may assume $D(\Omega)=\mathcal O$,
so $\{F_1,F_2\}=\{\mathcal O_5^+,\mathcal O_5^-\}$. The group $B$ acts freely on $\mathcal O$, $\mathcal O_5^+$ and
$\mathcal O_5^-$, and an invariant function for $B$ on $\mathcal O$ tends to opposite ends of $\R$ at
$\mathcal O_5^+$ and $\mathcal O_5^-$. We show that these requirements are incompatible.

Since every $X\in\mathfrak b$ is nilpotent, $Xp\in\R p$ forces $Xp=0$. So freeness at
$[p]$ means that $X\mapsto Xp$ is injective on $\mathfrak b$.

\emph{Step 1: normal form.} Freeness at $[e_4]\in \mathcal O_5^+$ shows that
$X\mapsto Xe_4=X_{14}e_1+X_{24}e_2+X_{34}e_3$ is injective on $\mathfrak b$, hence an
isomorphism $\mathfrak b\to\langle e_1,e_2,e_3\rangle$. Choose the basis
$X_1,X_2,X_3$ of $\mathfrak b$ with $X_ie_4=e_i$, and write
\[
  X_i=E_{i4}+A_i+\textstyle\sum_{k=1}^3 d_{ki}E_{k5}+h_iE_{45},
\]
where $A_i$ is strictly upper triangular and supported in the upper left $3\times3$
block, and $(A_i)_{jk}$ denotes its entries. A direct computation of the brackets
$[X_i,X_j]$ shows that $\mathfrak b$ is abelian if and only if the following hold:
$A_1=0$; $(A_2)_{12}(A_3)_{23}=0$; $(A_2)_{23}=0$; $(A_2)_{13}=(A_3)_{12}$;
$h_1=h_2=0$; $h_3=(A_3)_{23}\,d_{32}$; $(A_3)_{23}\,d_{31}=0$;
$h_3=(A_3)_{12}d_{21}+(A_3)_{13}d_{31}$; $(A_2)_{12}\,d_{21}+(A_2)_{13}\,d_{31}=0$; and
$(A_3)_{12}d_{22}-(A_2)_{13}d_{33}+(A_3)_{13}d_{32}=(A_2)_{12}d_{23}$.

\emph{Step 2: the case $h_3=0$ is impossible.} Then $Xe_5\in\langle e_1,e_2,e_3\rangle$
for every $X\in\mathfrak b$, and there is a linear map $T$ of
$\langle e_1,e_2,e_3\rangle$ with $Xe_5=T(Xe_4)$ for all $X\in\mathfrak b$. Let
$u\neq0$ be a real eigenvector of $T$, say $Tu=\mu u$, and let $X\in\mathfrak b$ be the
element with $Xe_4=u$. Then $X(e_5-\mu e_4)=Tu-\mu u=0$, so $X$ fixes the point
$[e_5-\mu e_4]\in \mathcal O$. This contradicts freeness on $\mathcal O$.

\emph{Step 3: the case $h_3\neq0$.} Put $a=(A_3)_{23}$. The relations of Step~1 give
$a\neq0$, $d_{32}\neq0$, $d_{31}=0$ and $(A_2)_{12}=0$. Then $h_3=(A_3)_{12}d_{21}$, so
$d_{21}\ne0$, and moreover
\[
  (A_2)_{13}=(A_3)_{12}=\frac{a\,d_{32}}{d_{21}},\qquad
  (A_3)_{13}=\frac{a\,(d_{33}-d_{22})}{d_{21}},\qquad h_3=a\,d_{32}.
\]
Let $\Delta=(d_{ki})$, a $3\times3$ matrix with $d_{31}=0$, and let $\sigma_1,\sigma_2,
\sigma_3$ be its trace, the sum of its principal $2\times2$ minors, and its determinant.
A direct computation shows that $f=P/z^4$ is $B$-invariant on $\mathcal O$, where
\[
\begin{aligned}
  P={}&3y^4+z\bigl(4\sigma_1y^3-12a\,d_{32}\,xy^2\bigr)\\
   &+z^2\bigl(6a^2d_{32}^2x^2-12a\,d_{32}(d_{11}+d_{22})\,xy+12a\,d_{32}^2\,wy
     +6\sigma_2y^2\bigr)\\
   &+z^3\bigl(12a\,d_{32}^2(d_{11}w-d_{21}v)+12a\,d_{32}(d_{12}d_{21}-d_{11}d_{22})\,x
     +12\sigma_3y\bigr).
\end{aligned}
\]
The coefficient of $v$ in $f$ is $-12a\,d_{32}^2d_{21}/z\neq0$. So $f$ is a submersion
whose fibers are graphs over the remaining coordinates, hence connected, and by
Lemma \ref{lem:invariant-function} it is an invariant function.

At every point of $\mathcal O_5^+\cup \mathcal O_5^-$ we have $z=0$ and $y\neq0$, so $P=3y^4>0$ there. By
continuity, $P$ is bounded below by a positive constant near each such point, while
$z^4\to0^+$. Hence $f\to+\infty$ at both $\mathcal O_5^+$ and $\mathcal O_5^-$, so the two faces lie at
the same end  (Corollary \ref{cor:parity} with $k=4$), contradicting Lemma \ref{lem:bricks}(1).
\end{proof}

\begin{remark}
There are examples where $B \subset N_5$ acts properly on $\mathcal{O}$ and two of its codimension 1 faces.
For instance,
\[
  \mathfrak b=\langle E_{14}+E_{25},\ E_{24}+E_{13}+E_{35},\
  E_{34}+E_{12}+E_{23}+E_{45}\rangle
\]
is abelian and $B = \exp(B)$ acts freely on $\mathcal O$, $\mathcal O_5^+$ and $\mathcal O_5^-$: the minors of the matrix
$(X_1p,X_2p,X_3p)$ include $z^3$ and $y^3$. 
Nevertheless, the action of $B$ on the union $\mathcal O_5^- \cup \mathcal O_5 \cup \mathcal O_5^+$ is not proper:
the invariant function has even degree $4$ in $y$ at infinity, so both faces lie at the
end $+\infty$.
\end{remark}

\subsection{\texorpdfstring{$N_{4,1}$}{N41} }\label{section:N41}

Here $\mathfrak n=\R D\oplus\mathfrak u$, with $D=\diag(1,1,1,1,-4)$ and $\mathfrak u$ the
strictly upper triangular matrices in the upper left $4\times4$ block. The open orbit
$\mathcal O=\{y>0,\ z>0\}$ has the faces
\[
  \mathcal O_4^\pm=\{y=0,\ \pm x>0,\ z>0\},\qquad \mathcal O_5=\{z=0,\ y>0\}.
\]
 
\emph{Normal forms.} Suppose first that $\mathfrak b\subset\mathfrak u$. Every
$X\in\mathfrak u$ satisfies $Xe_5=0$ and $Xe_4\in\langle e_1,e_2,e_3\rangle$, so freeness at
$[e_4+e_5]\in\mathcal O$ makes $X\mapsto Xe_4$ an isomorphism
$\mathfrak b\to\langle e_1,e_2,e_3\rangle$. Hence $Xe_3=T(Xe_4)$ for a linear map $T$ with
image in $\langle e_1,e_2\rangle$. If $Tu=0$ with $u\neq0$, the element $X$ with $Xe_4=u$
fixes $[\pm e_3+e_5]\in\mathcal O_4^\pm$. So only $\mathcal O_5$ is available, and there is no
brick. Hence $\mathfrak b=\langle D+u_1,b_2,b_3\rangle$ with $u_1,b_2,b_3\in\mathfrak u$.
 
A brick needs two faces, so $B$ must be free on $\mathcal O_4^\pm$. At
$p=(v,w,x,0,z)$, the nilpotent elements act by
$Xp=(X_{12}w+X_{13}x)e_1+X_{23}x\,e_2$, so $B$ is free on $\mathcal O_4^\pm$ if and only if every
nonzero $X\in\langle b_2,b_3\rangle$ has $X_{23}\ne0$, or $X_{12}=0\ne X_{13}$. Normalizing, and
imposing that $\mathfrak b$ is abelian, we obtain
\[
  b_1=D+dE_{14}+eE_{24},\qquad b_2=gE_{12}+E_{23}+jE_{14}+kE_{24}+gq\,E_{34},\qquad
  b_3=E_{13}+pE_{14}+qE_{24},
\]
with $eg=0$. This gives two cases: (a) $g\neq0$ and $e=0$, and (b) $g=0$.
 
\emph{Freeness on $\mathcal O$.} In case (a), $B$ is free on $\mathcal O$ if and only if
$q\neq0$: if $q=0$, then $b_3$ fixes the points with $x=-py$, and if $q\ne0$, the minor
$-5gq^2y^3z$ is nonzero on $\mathcal O$. In case (b), a nonzero
$\beta b_2+\gamma b_3$ fixes a point of $\mathcal O$ if and only if
$s^2+(p+k)s+pk-jq$ has a real root $s=x/y$, so $B$ is free on $\mathcal O$ if and only if
$(p-k)^2+4jq<0$.
 
\emph{Invariant functions.} In case (b) the function $f=x/y$ is invariant. In case (a),
\[
  f=-10\frac{x^3}{y^3}+30gq\frac{wx}{y^2}-15(k+p)\frac{x^2}{y^2}-30gq^2\frac vy
    +30gpq\frac wy+30(jq-kp)\frac xy+6dgq^2\log\frac yz
\]
is invariant. Both are invariant functions by Lemma~\ref{lem:conventions}(3), with
$x_i=x$ in case (b) and $x_i=v$ in case (a).
 
\emph{Faces.} At $\mathcal O_4^\pm$ the function $f$ has pole type with $k=1$, $Q=x$ in
case (b), and with $k=3$, $Q=-10x^3$ in case (a). So $\mathcal O_4^+$ and $\mathcal O_4^-$ lie at
opposite ends and are alternating. At $\mathcal O_5$, either $f$ has a finite limit (case (b),
or case (a) with $d=0$), so $\mathcal O_5$ is attached to no end and lies in no brick, or $f$
has logarithmic type, so $\mathcal O_5$ is constant. In the latter case it is the only
constant face, so by Lemma~\ref{lem:conventions}(4) it is never crossed.
 
\begin{proposition}\label{prop:N41}
Let $V=\T^4$ carry a projective structure of rank $3$ with respect to $N_{4,1}$. Then, up
to conjugation, $\mathfrak b$ is as above, in case (a) with $q\neq0$ or in case (b) with
$(p-k)^2+4jq<0$. The only brick is $\mathcal O\cup\mathcal O_4^+\cup\mathcal O_4^-$, the chain of
$\Vt$ alternates between $\mathcal O$ and $R_4\mathcal O$, and the structure is affine.
\end{proposition}
 
\begin{proof}
The chain never crosses $\{z=0\}$, so the structure is affine by
Lemma~\ref{lem:conventions}(5).
\end{proof}

\subsection{\texorpdfstring{$N_{3,2}$}{N32}}

Here $\mathfrak n=\R D\oplus\mathfrak u$, with $D=\diag(2,2,2,-3,-3)$ central and
$\mathfrak u=\langle E_{12},E_{13},E_{23},E_{45}\rangle$, so that
\begin{equation}\label{eq:N32-bracket}
  [X,Y]=(X_{12}Y_{23}-X_{23}Y_{12})\,E_{13}\qquad(X,Y\in\mathfrak n).
\end{equation}
The open orbit $\mathcal O=\{x>0,\ z>0\}$ has the faces
\[
  \mathcal O_3^\pm=\{x=0,\ \pm w>0,\ z>0\},\qquad \mathcal O_5^\pm=\{z=0,\ \pm y>0,\ x>0\}.
\]
 
\emph{Freeness.} An element with a nonzero $D$-coefficient has no eigenray on
$\mathcal O$ or its faces: there one of $x,w$ and one of $y,z$ is nonzero, so an eigenvalue
would equal both $2$ and $-3$ times that coefficient. 
For 
$X \in W=\mathfrak b\cap\mathfrak u$ we may write
$Xp=(X_{12}w+X_{13}x)e_1+X_{23}x\,e_2+X_{45}z\,e_4$ for $X\in\mathfrak u$, and observe that the group $B$ acts
freely
\begin{itemize}
  \item on $\mathcal O_3^\pm$ iff no nonzero $X\in W$ has $X_{12}=X_{45}=0$;
  \item on $\mathcal O_5^\pm$ iff every nonzero $X\in W$ with $X_{23}=0$ has $X_{12}=0\neq X_{13}$;
  \item on $\mathcal O$ iff no nonzero $X\in W$ has $X_{23}=X_{45}=0\neq X_{12}$.
\end{itemize}
If $\mathfrak b\subset\mathfrak u$, then $\dim W=3$, and nonzero elements with
$X_{12}=X_{45}=0$, respectively $X_{13}=X_{23}=0$, show that $B$ would not act freely on any face. 
Hence $\dim W=2$.
 
\emph{Normal forms.} By \eqref{eq:N32-bracket}, $\mathfrak b$ is abelian iff the vectors
$(X_{12},X_{23})$, $X\in\mathfrak b$, span at most a line.
\begin{itemize}
  \item[(A)] If $B$ is free on $\mathcal O_3^\pm$, then $X\mapsto(X_{12},X_{45})$ is an
  isomorphism $W\to\R^2$, and
  \[
    b_1=D+cE_{13},\qquad b_2=E_{12}+gE_{13}+eE_{23},\qquad b_3=kE_{13}+E_{45}.
  \]
  $B$ is free on $\mathcal O$ iff $e\neq0$, and on $\mathcal O_5^\pm$ iff $k\neq0$.
  \item[(B)] If $B$ is not free on $\mathcal O_3^\pm$, it must be free on $\mathcal O$ and
  $\mathcal O_5^\pm$. Then the vectors $(X_{12},X_{23})$ span a line, the kernel of
  $X\mapsto(X_{12},X_{23})$ on $W$ is spanned by an element with nonzero $E_{13}$-entry,
  and
  \[
    b_1=D+dE_{45},\qquad b_2=fE_{12}+E_{23}+hE_{45},\qquad b_3=E_{13}+lE_{45},
  \]
  with $fl=0$. This $B$ is always free on $\mathcal O$ and $\mathcal O_5^\pm$, and free on
  $\mathcal O_3^\pm$ iff $fl\neq0$.
\end{itemize}
 
\emph{Invariant functions.} The functions
\[
  f_A=5\frac vx-\frac{5g}{e}\frac wx-\frac{5}{2e}\frac{w^2}{x^2}-5k\frac yz+c\log\frac zx,
  \qquad
  f_B=5\frac yz-5l\frac vx-5h\frac wx+d\log\frac zx
\]
are invariant functions, by Lemma~\ref{lem:conventions}(3) with $x_i=v$ for $f_A$ and
$x_i=y$ for $f_B$.
 
\emph{Faces.} At $\mathcal O_5^\pm$ both functions have pole type with $k=1$ ($Q=-5ky$,
respectively $Q=5y$), so $\mathcal O_5^\pm$ lie at opposite ends and are alternating. In
case (A), at $\mathcal O_3^\pm$ the function $f_A$ has pole type with $k=2$ and
$Q|_{\{x=0\}}=-\tfrac5{2e}w^2$, so $\mathcal O_3^\pm$ lie at the same end and are constant. If
$k=0$, these are the only faces, so there is no brick. If $k\neq0$, they are the only
constant faces, so by Lemma~\ref{lem:conventions}(4) they are never crossed.
 
\begin{proposition}\label{prop:N32}
Let $V=\T^4$ carry a projective structure of rank $3$ with respect to $N_{3,2}$. Then, up
to conjugation, $\mathfrak b$ is as in case (A) with $ek\neq0$, or as in case (B) with
$fl=0$. The only brick is $\mathcal O\cup\mathcal O_5^+\cup\mathcal O_5^-$, the chain of $\Vt$
alternates between $\mathcal O$ and $R_5\mathcal O$, and the structure is affine.
\end{proposition}
 
\begin{proof}
The chain never crosses $\{x=0\}$, so the structure is affine by
Lemma~\ref{lem:conventions}(5).
\end{proof}

\subsection{\texorpdfstring{$N_{3,1,1}$}{N31}}

Here $\mathfrak n=\mathfrak a\oplus\mathfrak u$, with $\mathfrak a=\langle D_1,D_2\rangle$ central,
$D_1=\diag(2,2,2,-3,-3)$, $D_2=\diag(0,0,0,1,-1)$, and
$\mathfrak u=\langle E_{12},E_{13},E_{23}\rangle$, so that $[X,Y]=(X_{12}Y_{23}-X_{23}Y_{12})E_{13}$.
The open orbit $\mathcal O=\{x>0,\ y>0,\ z>0\}$ has the faces
\[
  \mathcal O_3^\pm=\{x=0,\ \pm w>0,\ y,z>0\},\qquad
  \mathcal O_4=\{y=0,\ x,z>0\},\qquad \mathcal O_5=\{z=0,\ x,y>0\}.
\]
 
\emph{Normal forms.} Since $\mathfrak u$ is not abelian, $\mathfrak b\to\mathfrak a$ is nonzero,
and $\mathfrak b$ is abelian iff the vectors $(X_{12},X_{23})$, $X\in\mathfrak b$, span at most a
line. If $\mathfrak b\to\mathfrak a$ is onto, then $\mathfrak b\cap\mathfrak u$ is spanned by one
element $iE_{12}+hE_{13}+gE_{23}$, which fixes a point of $\mathcal O$ iff $g=0$ and ($i\ne0$ or
$h=0$). Conjugating by a diagonal matrix and subtracting, we obtain three subcases. If
$\mathfrak b\to\mathfrak a$ has rank $1$, then $\mathfrak b\cap\mathfrak u$ is a $2$-dimensional abelian
subalgebra of $\mathfrak u$, so it contains $E_{13}$. Altogether:
\[
\begin{array}{lll}
  \text{(1a)} & b_1=D_1+bE_{13},\ \ b_2=D_2+eE_{13}, & b_3=E_{12}+hE_{13}+E_{23};\\[2pt]
  \text{(1b)} & b_1=D_1+bE_{13},\ \ b_2=D_2+eE_{13}, & b_3=hE_{13}+E_{23};\\[2pt]
  \text{(1c)} & b_1=D_1+sU,\ \ b_2=D_2+tU, & b_3=E_{13},\qquad U=pE_{12}+qE_{23};\\[2pt]
  \text{(2)} & b_1=aD_1+bD_2,\ \ b_2=dE_{12}+cE_{23}, & b_3=E_{13}.
\end{array}
\]
 
\emph{Freeness and invariant functions.} By a computation of minors, $B$ is free on the
orbits as follows:
\begin{center}
\renewcommand{\arraystretch}{1.2}
\begin{tabular}{lllll}
 & $\mathcal O$ & $\mathcal O_3^\pm$ & $\mathcal O_4$ & $\mathcal O_5$ \\ \hline
(1a) & always & always & $b\neq5e$ & $b\neq-5e$ \\
(1b) & always & never & $b\neq5e$ & $b\neq-5e$ \\
(1c) & always & never & $q(s-5t)\neq0$ & $q(s+5t)\neq0$ \\
(2) & $c\neq0$ & never & $c(5a+b)\neq0$ & $c(5a-b)\neq0$
\end{tabular}
\end{center}
The following are invariant functions, by Lemma~\ref{lem:conventions}(3) with $x_i=v$ in
(1a) and (1b), $x_i=w$ in (1c), and logarithms in (2):
\[
\begin{aligned}
  \text{(1a)}\quad & -5\,\frac{2hxw+w^2-2vx}{x^2}+(b-5e)\log\frac yx+(b+5e)\log\frac zx,\\
  \text{(1b)}\quad & 10\,\frac{v-hw}{x}+(b-5e)\log\frac yx+(b+5e)\log\frac zx,\\
  \text{(1c)}\quad & \frac{10}{q}\,\frac wx+(s-5t)\log\frac yx+(s+5t)\log\frac zx,\\
  \text{(2)}\quad & (5a-b)\log\frac zx-(5a+b)\log\frac yx .
\end{aligned}
\]
 
\emph{Faces.} Every face is constant: $\mathcal O_3^\pm$ has pole type with $k=2$ and
$Q|_{\{x=0\}}=-5w^2$, so both lie at $-\infty$; $\mathcal O_4$ and $\mathcal O_5$ have
logarithmic type and lie at the following ends:
\begin{center}
\renewcommand{\arraystretch}{1.2}
\begin{tabular}{lll}
 & $\mathcal O_4$ & $\mathcal O_5$ \\ \hline
(1a), (1b) & $-\sgn(b-5e)$ & $-\sgn(b+5e)$ \\
(1c) & $-\sgn(s-5t)$ & $-\sgn(s+5t)$ \\
(2) & $\sgn(5a+b)$ & $-\sgn(5a-b)$
\end{tabular}
\end{center}
By Lemma~\ref{lem:conventions}(4), the entry end changes at every crossing, so a chain
alternates between faces at $+\infty$ and faces at $-\infty$; the bricks are the unions of
$\mathcal O$ with one free face from each end, and every periodic sequence alternating in
this way is admissible.
 
\begin{proposition}\label{prop:N311}
Let $V=\T^4$ carry a projective structure of rank $3$ with respect to $N_{3,1,1}$. Then, up
to conjugation, $\mathfrak b$ is one of the following.
\begin{enumerate}
  \item[(1a)] The faces $\mathcal O_3^\pm$ lie at $-\infty$, and:
  \begin{itemize}
    \item if $b<-5|e|$, then $\mathcal O_4,\mathcal O_5$ lie at $+\infty$, and there are $4$ bricks;
    \item if $|b|<5|e|$, then exactly one of $\mathcal O_4,\mathcal O_5$ lies at $+\infty$, and
    there are $3$ bricks, each containing it;
    \item if $b=5e<0$, then only $\mathcal O_5$ lies at $+\infty$ and $\mathcal O_4$ is not
    free, and there are $2$ bricks; similarly if $b=-5e<0$, with $\mathcal O_4$ and
    $\mathcal O_5$ exchanged;
    \item in all other cases there is no brick.
  \end{itemize}
  There are infinitely many primitive words, hence infinitely many structures with
  pairwise non-isomorphic universal covers by isomorphisms normalizing $N$ (Proposition \ref{prop:words-distinct}). A structure whose word
  avoids $\mathcal O_4$ (respectively $\mathcal O_5$) is affine.
  \item[(1b)] With $|b|<5|e|$, (1c) with $q\neq0$ and $|s|<5|t|$, and (2) with $c\neq0$
  and $|b|<5|a|$: the only brick is $\mathcal O\cup\mathcal O_4\cup\mathcal O_5$, the chain has
  period $4$, through $\mathcal O$, $R_4\mathcal O$, $R_{45}\mathcal O$, $R_5\mathcal O$, and the
  structure is affine.
\end{enumerate}
\end{proposition}
 
\begin{proof}
The brick counts follow from the tables: for instance, $b<-5|e|$ means $b-5e<0$ and
$b+5e<0$. In (1b), (1c) and (2), $\mathcal O_4$ and $\mathcal O_5$ lie at opposite ends
exactly under the stated inequalities. In (1a), both faces $\mathcal O_3^\pm$ lie at $-\infty$, so the chain has a free
choice between them at every crossing at $-\infty$, which gives infinitely many primitive
words. The affineness statements follow from Lemma~\ref{lem:conventions}(5):
in (1b), (1c) and (2) the chain never crosses $\{x=0\}$.
\end{proof}
 
\begin{remark}
Case (1a) is the example of Remark \ref{rem:const-opp} in which constant faces lie at both ends. In
(1b), (1c) and (2), the element $b_3$ fixes $\mathcal O_3^\pm$ pointwise, and only the
logarithmic faces $\mathcal O_4,\mathcal O_5$ are available.
\end{remark}

\subsection{\texorpdfstring{$N_{2,2,1}$}{N221}}

Here $\mathfrak n=\mathfrak a\oplus\mathfrak u$ is abelian, with $\mathfrak a=\langle D_1,D_2\rangle$,
$D_1=\diag(1,1,0,0,-2)$, $D_2=\diag(0,0,1,1,-2)$, and $\mathfrak u=\langle E_{12},E_{34}\rangle$.
The open orbit $\mathcal O=\{w>0,\ y>0,\ z>0\}$ has the faces
\[
  \mathcal O_2^\pm=\{w=0,\ \pm v>0,\ y,z>0\},\qquad
  \mathcal O_4^\pm=\{y=0,\ \pm x>0,\ w,z>0\},
\]
and $\mathcal O_5=\{z=0,\ w,y>0\}$.
Exchanging the two blocks of size $2$ (Lemma~\ref{lem:conventions}(1)) exchanges
$D_1\leftrightarrow D_2$, $E_{12}\leftrightarrow E_{34}$ and $\mathcal O_2^\pm\leftrightarrow\mathcal O_4^\pm$.
 
\emph{Normal forms.} If $\mathfrak b\to\mathfrak a$ has rank $1$, then $\mathfrak b\supset\mathfrak u$,
and $E_{12}$, $E_{34}$ fix $\mathcal O_2^\pm$, respectively $\mathcal O_4^\pm$, pointwise, so only
$\mathcal O_5$ is left and there is no brick. Hence $\mathfrak b\cap\mathfrak u$ is spanned by one
element $eE_{12}+fE_{34}$, and after subtracting, rescaling and possibly exchanging the
blocks:
\[
\begin{array}{lll}
  \text{(1)}\ e,f\neq0: & b_1=D_1+aE_{12},\quad b_2=D_2+dE_{34}, & b_3=E_{12}+E_{34};\\[2pt]
  \text{(2)}\ e=0: & b_1=D_1+aE_{12},\quad b_2=D_2+cE_{12}, & b_3=E_{34}.
\end{array}
\]
 
\emph{Freeness and invariant functions.} By a computation of minors, in case (1) $B$ is
free on $\mathcal O$, $\mathcal O_2^\pm$ and $\mathcal O_4^\pm$; in case (2) it is free on $\mathcal O$
and $\mathcal O_2^\pm$, while $b_3=E_{34}$ fixes $\mathcal O_4^\pm$ pointwise. (On $\mathcal O_5$, $B$ is
free iff $a\neq d$, respectively $a+c\neq0$; this plays no role.) The functions
\[
\begin{aligned}
  \text{(1)}\quad & -\frac vw+\frac xy-d\log\frac yw+\frac{d-a}{5}\log\frac{y^2z}{w^3},\\
  \text{(2)}\quad & 5\frac vw+(2c-3a)\log w+(2a-3c)\log y+(a+c)\log z
\end{aligned}
\]
are invariant functions, by Lemma~\ref{lem:conventions}(3) with $x_i=v$.
 
\emph{Faces.} At $\mathcal O_2^\pm$ the functions have pole type with $k=1$ ($Q=-v$,
respectively $Q=5v$), and in case (1) at $\mathcal O_4^\pm$ pole type with $k=1$ and $Q=x$.
So each pair lies at opposite ends, and these faces are alternating. The face
$\mathcal O_5$ has logarithmic type, so it is constant; being the only constant face, it is
never crossed (Lemma~\ref{lem:conventions}(4)). Since all faces in use are alternating,
the chain enters every copy at the same end.
 
\begin{proposition}\label{prop:N221}
Let $V=\T^4$ carry a projective structure of rank $3$ with respect to $N_{2,2,1}$. Then, up
to conjugation, $\mathfrak b$ is as in case (1) or (2).
\begin{itemize}
  \item In case (1) there are four bricks,
  $\mathcal O\cup\mathcal O_2^+\cup\mathcal O_2^-$, $\mathcal O\cup\mathcal O_2^-\cup\mathcal O_4^-$,
  $\mathcal O\cup\mathcal O_2^+\cup\mathcal O_4^+$ and $\mathcal O\cup\mathcal O_4^+\cup\mathcal O_4^-$, and at
  every step the chain may cross either $\{w=0\}$ or $\{y=0\}$. There are infinitely many
  primitive words, hence infinitely many structures with pairwise non-isomorphic
  universal covers by isomorphisms normalizing $N$ (Proposition \ref{prop:words-distinct}). A structure whose word never crosses $\{y=0\}$
  (respectively $\{w=0\}$) is affine.
  \item In case (2) the only brick is $\mathcal O\cup\mathcal O_2^+\cup\mathcal O_2^-$, the chain
  alternates between $\mathcal O$ and $R_2\mathcal O$, and the structure is affine.
\end{itemize}
\end{proposition}
 
\begin{proof}
In case (1), in $\mathcal O$ the faces $\mathcal O_2^-$ and $\mathcal O_4^+$ lie at one end and
$\mathcal O_2^+$, $\mathcal O_4^-$ at the other, and since all faces are alternating, every copy has
one face in $\{w=0\}$ and one in $\{y=0\}$ at its exit end. This gives the bricks and the
free choice at every step. The affineness statements follow from
Lemma~\ref{lem:conventions}(5).
\end{proof}

\subsection{\texorpdfstring{$N_{2,1,1,1}$}{N2111}}

Here $\mathfrak n=\mathfrak a\oplus\mathfrak u$ is abelian, with $\mathfrak a=\langle D_1,D_2,D_3\rangle$,
$D_1=\diag(1,1,0,0,-2)$, $D_2=\diag(0,0,1,0,-1)$, $D_3=\diag(0,0,0,1,-1)$, and
$\mathfrak u=\langle E_{12}\rangle$. The open orbit $\mathcal O=\{w,x,y,z>0\}$ has the faces
\[
  \mathcal O_2^\pm=\{w=0,\ \pm v>0,\ x,y,z>0\}
\]
and $\mathcal O_3=\{x=0\}$, $\mathcal O_4=\{y=0\}$, $\mathcal O_5=\{z=0\}$ (with the other three of
$w,x,y,z$ positive). Exchanging $x$ and $y$ (Lemma~\ref{lem:conventions}(1)) exchanges
$D_2\leftrightarrow D_3$ and $\mathcal O_3\leftrightarrow\mathcal O_4$.
 
\emph{Normal forms.} Either $\mathfrak b\cap\mathfrak u=0$, or $E_{12}\in\mathfrak b$ and
$\mathfrak b=\mathfrak a'\oplus\langle E_{12}\rangle$ for a plane $\mathfrak a'\subset\mathfrak a$. In the
second case, if some element of $\mathfrak a'$ has a nonzero $D_1$-coefficient, then, after
exchanging $x$ and $y$ if necessary, $\mathfrak a'$ projects isomorphically onto
$\langle D_1,D_2\rangle$ along $D_3$; otherwise $\mathfrak a'=\langle D_2,D_3\rangle$. This gives:
\[
\begin{array}{lll}
  \text{(1)} & b_1=D_1+aE_{12},\quad b_2=D_2+bE_{12}, & b_3=D_3+cE_{12};\\[2pt]
  \text{(2a)} & b_1=D_1+aD_3,\quad b_2=D_2+bD_3, & b_3=E_{12};\\[2pt]
  \text{(2b)} & b_1=D_2,\quad b_2=D_3, & b_3=E_{12}.
\end{array}
\]
 
\emph{Freeness and invariant functions.} Put
\[
\begin{aligned}
  &\alpha_3=a-4b+c,\quad \alpha_4=a+b-4c,\quad \alpha_5=a+b+c,\\
  &\beta_3=a-4b-1,\quad \beta_4=a+b+4,\quad \beta_5=a+b-1.
\end{aligned}
\]
By a computation of minors, $B$ is free on $\mathcal O$ in all cases (in (2a) because
$\beta_4-\beta_5=5$), and on the faces as follows:
\begin{center}
\renewcommand{\arraystretch}{1.2}
\begin{tabular}{lllll}
 & $\mathcal O_2^\pm$ & $\mathcal O_3$ & $\mathcal O_4$ & $\mathcal O_5$ \\ \hline
(1) & always & $\alpha_3\neq0$ & $\alpha_4\neq0$ & $\alpha_5\neq0$ \\
(2a) & never & $\beta_3\neq0$ & $\beta_4\neq0$ & $\beta_5\neq0$ \\
(2b) & never & always & always & always
\end{tabular}
\end{center}
The following are invariant functions, by Lemma~\ref{lem:conventions}(3) with $x_i=v$ in
(1) and logarithms in (2a), (2b):
\[
\begin{aligned}
  \text{(1)}\quad & 5\frac vw+\alpha_3\log x+\alpha_4\log y+\alpha_5\log z+(2b+2c-3a)\log w,\\
  \text{(2a)}\quad & \beta_3\log\frac xw+\beta_4\log\frac yw+\beta_5\log\frac zw,\qquad
  \text{(2b)}\quad \log\frac{xyz}{w^3}.
\end{aligned}
\]
 
\emph{Faces.} In (1), $\mathcal O_2^\pm$ have pole type with $k=1$ and $Q=5v$, so they are
alternating, with $\mathcal O_2^+$ at $+\infty$. The faces $\mathcal O_3,\mathcal O_4,\mathcal O_5$ have
logarithmic type and are constant; $\mathcal O_j$ lies at $-\sgn(\alpha_j)$ in (1),
$-\sgn(\beta_j)$ in (2a), and at $-\infty$ in (2b). Let $S_\pm$ be the sets of free constant
faces at $\pm\infty$.
 
\begin{proposition}\label{prop:N2111}
Let $V=\T^4$ carry a projective structure of rank $3$ with holonomy group almost contained in $N_{2,1,1,1}$. Then, up
to conjugation, $\mathfrak b$ is as in (1) or (2a), and:
\begin{itemize}
  \item in case (1) with $S_+$ or $S_-$ empty (the nonzero $\alpha_j$ all have the same
  sign, or all vanish), the only brick is $\mathcal O\cup\mathcal O_2^+\cup\mathcal O_2^-$, the chain
  alternates between $\mathcal O$ and $R_2\mathcal O$, and the structure is affine;
  \item in case (1) with $S_+$ and $S_-$ nonempty, the bricks are the unions of $\mathcal O$
  with one face from each end among $\mathcal O_2^\pm$ and $S_\pm$, so there are
  $(1+|S_+|)(1+|S_-|)$ of them, and there are infinitely many primitive words; a
  structure whose word avoids $\mathcal O_j$ ($j=3,4,5$) is affine;
  \item in case (2a), bricks exist iff $S_+$ and $S_-$ are nonempty, i.e.\ two of the
  $\beta_j$ are nonzero with opposite signs; the bricks are $\mathcal O\cup F_+\cup F_-$ with
  $F_\pm\in S_\pm$, so there are $|S_+|\,|S_-|\in\{1,2\}$ of them, and the structure is
  affine. With one brick the chain has period $4$; with two there are infinitely many
  primitive words.
\end{itemize}
In case (2b) all faces lie at $-\infty$, and there is no brick.
\end{proposition}
 
\begin{proof}
In case (1) with $S_+$ or $S_-$ empty, the free constant faces lie at one end and are
never crossed (Lemma~\ref{lem:conventions}(4)). Otherwise constant faces lie at both ends;
crossing them reverses the entry end and crossing $\mathcal O_2^\pm$ does not, and every
periodic sequence consistent with this is admissible, which gives the bricks and words.
In (2a) all faces in use are constant, so the chain alternates between $S_+$ and $S_-$.
Affineness follows from Lemma~\ref{lem:conventions}(5), since in case (2a) the chain never
crosses $\{w=0\}$.
\end{proof}

\subsection{The diagonal case, \texorpdfstring{$N_{1,1,1,1,1}$}{N11111}}\label{sec:diag}

Here $\mathfrak n$ consists of the traceless diagonal matrices. We write points as
$[x_1:\dots:x_5]$; the open orbit $\mathcal O=\{x_1,\dots,x_5>0\}$ has the faces
$\mathcal O_i=\{x_i=0,\ x_j>0\ \text{for}\ j\neq i\}$, and permutations of the coordinates may
be used as in Lemma~\ref{lem:conventions}(1). Since $\mathfrak b$ has codimension one in
$\mathfrak n$,
\[
  \mathfrak b=\Bigl\{\diag(X_1,\dots,X_5)\ :\ \textstyle\sum_iX_i=0,\ \sum_i\kappa_iX_i=0\Bigr\}
\]
for a nonzero $\kappa\in\R^5$ with $\sum_i\kappa_i=0$, unique up to scale. (If $\mathfrak b$ is
given modulo scalars as $\{\diag(a,b,c,d,0):\alpha a+\beta b+\gamma c+\delta d=0\}$, then
$\kappa=(\alpha,\beta,\gamma,\delta,-\alpha-\beta-\gamma-\delta)$.)
 
\emph{Freeness and invariant function.} A diagonal matrix fixes $[p]$ iff it is scalar on
the coordinates where $p\neq0$. So $B$ is free on $\mathcal O$, and on $\mathcal O_i$ it can only
fix points through a multiple of $X^{(i)}=\diag(1,\dots,-4,\dots,1)$ ($-4$ in position
$i$); since $\sum_j\kappa_jX^{(i)}_j=-5\kappa_i$, $B$ is free on $\mathcal O_i$ iff
$\kappa_i\neq0$. The function $f=\sum_i\kappa_i\log x_i$ is an invariant function
(Lemma~\ref{lem:conventions}(3)).
 
\emph{Faces.} Each free face $\mathcal O_i$ has logarithmic type, so it is constant and lies at
$-\sgn(\kappa_i)\cdot\infty$. Let $S_+=\{\mathcal O_i:\kappa_i<0\}$, $S_-=\{\mathcal O_i:\kappa_i>0\}$ and
$n_\pm=|S_\pm|$. Since $\sum_i\kappa_i=0$ and $\kappa\neq0$, we have $n_\pm\ge1$ and
$n_++n_-\le5$. By Lemma~\ref{lem:conventions}(4), a chain alternates between $S_+$ and
$S_-$, and every periodic sequence alternating in this way is admissible.
 
\begin{proposition}\label{prop:diagonal}
Let $V=\T^4$ carry a projective structure of rank $3$ with respect to $N_{1,1,1,1,1}$. Then
the bricks are the sets $\mathcal O\cup F_+\cup F_-$ with $F_\pm\in S_\pm$, so there are $n_+n_-$
of them. If $n_+=n_-=1$, the words are the powers of the single brick, the chain has
period $4$, and the structure is affine. Otherwise there are infinitely many primitive
words, and a structure whose word avoids $\mathcal O_i$ is affine. Since $n_\pm\ge1$ for every
$\kappa$, every $3$-dimensional subgroup of $N$ is the group of a rank $3$ structure.
\end{proposition}
 
\begin{proof}
If $n_+=n_-=1$, say $S_\pm=\{\mathcal O_i\},\{\mathcal O_j\}$, the chain passes through $\mathcal O$,
$R_i\mathcal O$, $R_{ij}\mathcal O$, $R_j\mathcal O$; otherwise there is a choice of face at every
other step. Affineness follows from Lemma~\ref{lem:conventions}(5), and existence from
Lemma~\ref{lem:conventions}(6).
\end{proof}

\subsection{When the supergroup \texorpdfstring{$N$}{N} is not unique}\label{sec:overlaps}

Theorem~\ref{thm:rank3-summary} classifies structures with respect to a chosen $N$. 
When the holonomy group of a projective structure has finite index subgroup in $N \cap N'$, the structure is counted twice in our list.

The developing image is
then saturated by $N'$-orbits (Theorem \ref{theorem:benoistbricks}), and it contains $3$-dimensional $N'$-orbits,
so the structure also has rank $3$ with respect to $N'$.
 
To find all such cases, note that a conjugate of $N'=N_{i_1,\dots,i_k}$ is determined by a
decomposition $\R^5=W_1\oplus\dots\oplus W_k$ with $\dim W_\ell=i_\ell$. Hence
$\mathfrak b\subset\mathfrak n'$, up to conjugation, if and only if $\R^5$ splits into
$\mathfrak b$-invariant subspaces of dimensions $i_1,\dots,i_k$, each contained in a single
generalized weight space of $\mathfrak b$. Starting from the blocks of $N$, a different
partition arises in one of two ways: by merging two blocks on which $\mathfrak b$ has the
same weight, or by splitting a block on which $\mathfrak b$ acts decomposably. Checking
every case of Propositions~\ref{prop:N41}--\ref{prop:diagonal}, this happens, among
parameters for which structures exist, only in the following three families:
\begin{enumerate}
  \item[(I)] $N_{2,1,1,1}$, case (1) with $a=b=c=0$, and $N_{1,1,1,1,1}$ with
  $\kappa$ proportional to $(1,-1,0,0,0)$. Here
  $\mathfrak b=\{\diag(X_1,X_1,X_3,X_4,X_5)\}$ is the diagonal part of $\mathfrak n_{2,1,1,1}$,
  and the developing image is $\{x,y,z>0,\ (v,w)\neq(0,0)\}$.
  \item[(II)] $N_{3,1,1}$, case (2) with $b=0$, and $N_{3,2}$, case (B) with $d=h=l=0$.
  Here $\mathfrak b=\langle\diag(2,2,2,-3,-3),\ dE_{12}+cE_{23},\ E_{13}\rangle$, and the
  developing image is $\{x>0,\ (y,z)\neq(0,0)\}$.
  \item[(III)] $N_{2,2,1}$, case (2) with $a=c=0$, and $N_{2,1,1,1}$, case (2a) with
  $(a,b)$ equal to $(0,1)$, $(-2,-2)$ or $(-1,-\tfrac12)$, which are related by permuting
  the coordinates $x,y,z$. Here $\mathfrak b=\langle\diag(1,1,0,0,-2),\ \diag(0,0,1,1,-2),\
  E_{34}\rangle$, and the developing image is $\{y>0,\ z>0,\ (v,w)\neq(0,0)\}$.
\end{enumerate}
In each family the two descriptions cut the same developing image into bricks
differently: the chain has period $2$ for the coarser group and period $4$ for the finer
one. All other cases are excluded either because the weights of $\mathfrak b$ on the
blocks of $N$ are distinct and $\mathfrak b$ acts indecomposably on each block (in
$N_{4,1}$ case (b), indecomposability is equivalent to the freeness condition
$(p-k)^2+4jq<0$), or because the existence conditions fail at the coincidences, as for
$N_{3,1,1}$ case (2) with $b=\pm5a$.

Note that the holonomy representation may determine $N$ even when $B$ does not: an element $\gamma$ of $\Gamma' \setminus \Gamma_0$ may determine $N$ uniquely.
In particular the gluing element $g$ of Lemma \ref{lem:realize} may lie in the
centralizer of $B$ in $N$ but not in $N'$. For example, in (I) the element $\exp(E_{12})$
centralizes $B$ but is not diagonal. Such structures occur only for $N$.

It turns out that each of these three families admit only one brick for their respective group $B$.

Hence, by
Corollary~\ref{cor:rank3-dichotomy}, a rank $3$ structure on $\T^4$ with at least two
bricks determines $N$ uniquely up to conjugation, whereas a single-brick (affine)
structure may occur for two groups, exactly in the families (I)--(III).

\section{Maximal Rank}\label{sec:max_rannk}

In this section, with computations in the appendices, we prove Theorems~ \ref{thm:intro-masa} and \ref{thm:intro-rank4} from the
introduction.
 
\subsection{Reduction to maximal abelian groups}

If a projective structure on $\T^4$ has rank $4$ with group $B$, then $B$ is connected
and abelian, so it lies in the identity component of a maximal abelian subgroup of
$\SL(5,\R)$. To classify rank $4$ structures it therefore suffices to classify these
identity components and to decide which of them contain a $4$-dimensional connected
subgroup with an open orbit. Lemma~\ref{lem:reduction} shows that this is equivalent to
the group itself having an open orbit, and Lemma~\ref{lem:C-criterion} gives a criterion to test this. 

Fix the universal covering space $\R^4 \to \T^4$ with deck transformation group $\Z^4$. 

\begin{definition}\label{def:rank4}
A (marked) \emph{left-invariant projective structure} on the torus is a tuple $(B,h,[p])$
where: 
\begin{enumerate}
    \item $B\subset\SL(5,\R)$ is a connected abelian subgroup of dimension $4$ with universal cover $\pi \colon \widetilde{B} \to B$.
    \item $h \colon \Z^4 \to \widetilde{B}$ is an isomorphism onto a lattice of $\widetilde{B}$
    \item The $B$-orbit at $[p]$, that is $\Omega = B \cdot [p] \subset \Sph^4$, is open.
\end{enumerate}
\end{definition}

To be explicit, the associated dev-hol pair consists of the compositions
    \[ D \colon \R^4 \to \widetilde{B} \to \Omega, \quad D(v) = \pi(\bar{h}(v))[p], \quad hol \colon \Z^4 \to \widetilde{B} \to B, \quad hol(\gamma) = \pi(h(\gamma)) \]
where $\bar{h} \colon \R^4 \to \widetilde{B}$ is the unique Lie group isomorphism extending $h$. 
The fact that $D$ is a local homeomorphism follows from the fact that the action of $\widetilde{B}$ is locally free, and equivariance is clear since $\bar{h}$ restricts to $h$.

The developing image is $\Omega$, which turns out to always avoid a hyperplane, hence induces at least one affine structure.
Let us emphasize that the pullback of such an affine structure of $\Omega$ is almost always different from the standard affine structure on $\R^4$.

Since the exponential map identifies the Lie algebra $\mathfrak{b}$ of $B$ with $\widetilde{B}$, the map $h$ can be written as a map to $\mathfrak{b}$. 
The point $[p]$ can be varied within $\Omega$ by elements of $B$, yielding equivalent structures. 
The upshot is that the map $h \colon \Z^4 \to \mathfrak{b}$ determines the left-invariant structure up to equivalence, and two such structures are equivalent when they are conjugated by an element $g$ of $\SL^\pm(5,\R)$: more precisely, when $\Ad(g) \circ h=h'$.

The equivalent descriptions of rank from \Cref{prop:rank} extend to non-simply connected $N$ in the full rank case, with the caveat that the rank of $h(\Gamma_0)$ in $\widetilde{B}$ can differ from the rank of $\pi(h(\Gamma_0))$ in $B$, as in the case of the Hopf torus. 
This follows from the fact that the developing image is already a $B$-orbit of full dimension. 
We conclude (c.f.\ \cite[Section 4.3]{benoist1994nilvarietes})

\begin{proposition}
    A projective structure on $\T^4$ has rank $4$ if and only if it is equivalent to a left-invariant structure.
\end{proposition}

Let $H\subset\SL(5,\R)$ be a connected abelian subgroup with Lie algebra $\h$, and let
$X_1,\dots,X_k$ be a basis of $\h$. For $p\in\R^5\setminus\{0\}$ let
\[
  r(p)=\operatorname{rank}\bigl(X_1p,\ \dots,\ X_kp,\ p\bigr),
\]
the rank of a $5\times(k+1)$ matrix. When $k=4$ this matrix is square, and we put
\[
  \delta(p)=\det\bigl(X_1p,\ X_2p,\ X_3p,\ X_4p,\ p\bigr),
\]
a homogeneous polynomial of degree $5$; then $r(p)=5$ if and only if $\delta(p)\neq0$.
Changing the basis multiplies $\delta$ by a nonzero constant.
 
\begin{lemma}\label{lem:C-criterion}
The orbit $H\cdot[p]\subset\Sph^4$ is open if and only if $r(p)=5$. Consequently:
\begin{enumerate}
  \item $H$ has an open orbit if and only if $r(p)=5$ for some $p$; when $\dim H=4$,
  if and only if $\delta\not\equiv0$;
  \item when $\dim H=4$ and $\delta\not\equiv0$, the open orbits of $H$ are exactly the
  connected components of $\{[p]\in\Sph^4:\delta(p)\neq0\}$, and $H$ acts on each of them
  with discrete stabilizers;
  \item if $H$ has an open orbit, so does some $4$-dimensional connected subgroup of $H$.
\end{enumerate}
\end{lemma}
 
\begin{proof}
The tangent space of $\Sph^4$ at $[p]$ is $\R^5/\R p$, and the orbit map
$o\colon H\to\Sph^4$, $h\mapsto h\cdot[p]$, has differential at the identity
$\mathrm{d}o_e(X)=Xp \bmod \R p$. Hence $\mathrm{d}o_e$ is onto if and only if $X_1p,\dots,X_kp,p$ span
$\R^5$, that is, if and only if $r(p)=5$. Since $o(gh)=g\cdot o(h)$, the rank of $\mathrm{d}o$ is
the same at every point of $H$.
 
If $r(p)=5$, then $o$ is a submersion, so its image $H\cdot[p]$ is open. If $r(p)<5$,
then $\mathrm{d}o$ has rank at most $3$ everywhere, so every point of $H\cdot[p]$ is a critical
value of $o$; by Sard's theorem $H\cdot[p]$ has measure zero and is not open. This
proves the first statement and (1).
 
For (2): every point of $U=\{\delta\neq0\}$ has an open orbit, and orbits are disjoint,
so each connected component of $U$ is a disjoint union of open orbits, hence a single
orbit. Conversely, an open orbit lies in $U$ by the first part, hence equals the
component of $U$ containing it. At a point with $\delta(p)\ne0$ the map $\mathrm{d}o_e$ is an
isomorphism $\h\to\R^5/\R p$, so the stabilizer has trivial Lie algebra and is discrete.
 
For (3): if $r(p)=5$, choose $Y_1,\dots,Y_4\in\h$ whose images $Y_ip\bmod\R p$ form a
basis of $\R^5/\R p$. Since $\h$ is abelian, their span is a subalgebra, and the
corresponding connected subgroup has an open orbit at $[p]$.
\end{proof}

\begin{lemma}\label{lem:reduction}
Let $G\subset\SL(5,\R)$ be a connected abelian subgroup. The following are equivalent:
\begin{enumerate}
  \item $G$ contains a $4$-dimensional connected subgroup $B$ with an open orbit in
  $\Sph^4$;
  \item $G$ has an open orbit in $\Sph^4$;
  \item there is a projective structure on $\T^4$ of rank $4$ with a group $B\subset G$
  (Definition~\ref{def:rank4}).
\end{enumerate}
\end{lemma}

\begin{proof}
(1)$\Leftrightarrow$(2) is Lemma~\ref{lem:C-criterion}(3).

(3)$\Rightarrow$(1) is immediate from Definition~\ref{def:rank4}.

(1)$\Rightarrow$(3): let $\Omega=B\cdot[p]$ be an open orbit. By
Lemma~\ref{lem:C-criterion}(2) the stabilizer $\Delta$ of $[p]$ in $B$ is discrete, so
$\Omega=B/\Delta$. Let $\pi\colon\widetilde B\to B$ be the universal cover; since $B$ is
connected and abelian, $\widetilde B\cong\R^4$. The map
$D(\tilde b)=\pi(\tilde b)\cdot[p]$ is the composition of the covering maps
$\widetilde B\to B\to B/\Delta=\Omega$, hence a covering map, and it satisfies
$D(\tilde g\tilde b)=\pi(\tilde g)D(\tilde b)$. For any isomorphism $h$ of $\Z^4$ onto a lattice of
$\widetilde B$, the triple $(B,h,[p])$ therefore
defines a left-invariant projective structure on the torus which has rank $4$ with group $B$.
\end{proof}

Maximal abelian subgroups of $\SL(5,\R)$ need not be connected; in what follows we always
list their identity components, equivalently their Lie algebras.
 
\subsection{Maximal abelian subalgebras over \texorpdfstring{$\C$}{C} and over \texorpdfstring{$\R$}{R}}
 
Denote real and complex centralizers by $C_{\R}$ and $C_{\C}$.

\begin{lemma}\label{lem:complexification}
A real subalgebra $A_{\R}\subset\gl(5,\R)$ is maximal abelian if and only if its
complexification $A_{\C}=A_{\R}\otimes\C$ is maximal abelian in $\gl(5,\C)$.
\end{lemma}

\begin{proof}
An abelian subalgebra is maximal abelian if and only if it equals its centralizer. The
centralizer of $A_{\R}$ is the solution space of a system of real linear equations, and
the centralizer of $A_{\C}$ is the solution space of the same system over $\C$. Hence
$C_{\C}(A_{\C})=C_{\R}(A_{\R})\otimes\C$, and $C_{\R}(A_{\R})=A_{\R}$ if and only if
$C_{\C}(A_{\C})=A_{\C}$.
\end{proof}

So the real maximal abelian subalgebras are exactly the real forms of the complex ones.
Over $\C$ the structure is described by \cite{suprunenko1968commutative}: by \cite{suprunenko1968commutative} Theorem 2, every maximal
abelian subalgebra of $\gl(n,\C)$ is, up to conjugacy, a direct sum
\begin{equation}\label{eq:complex-blocks}
  A_{\C}=\bigoplus_{i=1}^k\bigl(\C I_{n_i}\oplus\n_i\bigr),\qquad n_1+\dots+n_k=n,
\end{equation}
where $\C^n=\bigoplus_i\C^{n_i}$ is the decomposition into generalized weight spaces of
$A_{\C}$, and $\n_i$ is a maximal abelian nilpotent subalgebra of $\gl(n_i,\C)$. The
latter are classified in \cite{suprunenko1968commutative} Chapter 3, Section 3 for $n_i\le6$.

A real form of \eqref{eq:complex-blocks} is determined by a complex conjugation
preserving $A_{\C}$. Such a conjugation permutes the weight spaces, so it fixes some
blocks and exchanges others in pairs of equal size. This gives two kinds of real blocks.

\begin{lemma}\label{lem:real-blocks}
Every maximal abelian subalgebra of $\gl(5,\R)$ is conjugate to a direct sum of real
blocks of the following two kinds:
\begin{enumerate}
  \item \emph{split blocks} $\R I_m\oplus\n$, where $\n$ is a real form of a maximal
  abelian nilpotent subalgebra of $\gl(m,\C)$; these come from blocks of
  \eqref{eq:complex-blocks} fixed by the conjugation;
  \item \emph{complex blocks}, the realifications of $\C I_m\oplus\n$ acting on
  $\C^m\cong\R^{2m}$; these come from pairs of blocks exchanged by the conjugation.
\end{enumerate}
Since $2m\le5$, a complex block has $m=1$ or $m=2$. For $m=1$ it is $\R I_2\oplus\bc$,
where $\bc$ is the rotation block of Appendix \ref{appendix:nonunipotent}. For $m=2$, the only maximal abelian
nilpotent subalgebra of $\gl(2,\C)$ is spanned by a single Jordan block, and the complex
block is $\R I_4\oplus\w8$, the realification of a $2\times2$ complex Jordan block.
\end{lemma}

In the language of \cite{WZ84}, the complex blocks are the ``indecomposable but not
absolutely indecomposable'' MASAs; $\w8$ is item~8 of \cite{WZ84}, in the corrected form of
\cite{diatta2022classification} Section 5.2. Lemma~\ref{lem:real-blocks} reduces the classification to two
tasks: find the real forms of the nilpotent algebras $\n$ of size $m\le5$, which is
done in Appendix \ref{appendix:realforms} using Galois cohomology; and assemble the blocks, which is done for
the non-nilpotent algebras in Appendix \ref{appendix:nonunipotent}.

\begin{remark}
Galois cohomology classifies the real forms up to $\GL(5,\R)$-conjugacy. Since $5$ is
odd, $-I$ has determinant $-1$ and $\GL(5,\R)=\R^\times\cdot\SL(5,\R)$, so this is the
same as $\SL(5,\R)$-conjugacy.
\end{remark}
 
\subsection{The unipotent groups}

The case of a single block of size $5$ with $\n$ nilpotent gives the unipotent groups.
We say that an abelian algebra $\n$ of nilpotent matrices has \emph{nilpotence class}
$l$ if $a^l=0$ for all $a\in\n$ and $a^{l-1}\neq0$ for some $a\in\n$. We denote the
unipotent groups by $N^k_{li}$, where $l$ is the nilpotence class, $i$ is the position
in the list of \cite{suprunenko1968commutative} (reordered slightly to put transposes next to each other, and
conjugated to be upper triangular), and $k=1,2$ labels the real forms when there is
more than one.

\begin{theorem}\label{thm:unipotent-list}
Up to conjugacy, there are $23$ maximal abelian unipotent subgroups of
$\SL(5,\R)$. They are represented by the following groups, with $a,\dots,f\in\R$:
\begin{align*}
\begin{array}{cccc}
  \makebox[3.3cm]{$N_{21}$} & \makebox[3.3cm]{$N_{22}$} & \makebox[3.3cm]{$N_{23}$} & \makebox[3.3cm]{$N_{24}$} \\[4pt]
  \makebox[3.3cm]{$\begin{pmatrix}1&0&0&0&d\\0&1&0&0&c\\0&0&1&0&b\\0&0&0&1&a\\0&0&0&0&1\end{pmatrix}$} &
  \makebox[3.3cm]{$\begin{pmatrix}1&a&b&c&d\\0&1&0&0&0\\0&0&1&0&0\\0&0&0&1&0\\0&0&0&0&1\end{pmatrix}$} &
  \makebox[3.3cm]{$\begin{pmatrix}1&0&0&a&b\\0&1&0&c&d\\0&0&1&e&f\\0&0&0&1&0\\0&0&0&0&1\end{pmatrix}$} &
  \makebox[3.3cm]{$\begin{pmatrix}1&0&a&b&c\\0&1&d&e&f\\0&0&1&0&0\\0&0&0&1&0\\0&0&0&0&1\end{pmatrix}$}
\end{array} \\[14pt]
\begin{array}{cccc}
  \makebox[3.3cm]{$N_{35}$} & \makebox[3.3cm]{$N_{36}$} & \makebox[3.3cm]{$N_{37}$} & \makebox[3.3cm]{$N^1_{38}$} \\[4pt]
  \makebox[3.3cm]{$\begin{pmatrix}1&a&0&b&c\\0&1&0&a&0\\0&0&1&d&e\\0&0&0&1&0\\0&0&0&0&1\end{pmatrix}$} &
  \makebox[3.3cm]{$\begin{pmatrix}1&0&0&a&d\\0&1&a&b&c\\0&0&1&0&b\\0&0&0&1&a\\0&0&0&0&1\end{pmatrix}$} &
  \makebox[3.3cm]{$\begin{pmatrix}1&a&b&c&d\\0&1&0&b&a\\0&0&1&a&0\\0&0&0&1&0\\0&0&0&0&1\end{pmatrix}$} &
  \makebox[3.3cm]{$\begin{pmatrix}1&0&b&a&d\\0&1&a&b&c\\0&0&1&0&b\\0&0&0&1&a\\0&0&0&0&1\end{pmatrix}$}
\end{array} \\[14pt]
\begin{array}{cccc}
  \makebox[3.3cm]{$N^2_{38}$} & \makebox[3.3cm]{$N^1_{39}$} & \makebox[3.3cm]{$N^2_{39}$} & \makebox[3.3cm]{$N_{310}$} \\[4pt]
  \makebox[3.3cm]{$\begin{pmatrix}1&0&b&-a&d\\0&1&a&b&c\\0&0&1&0&b\\0&0&0&1&a\\0&0&0&0&1\end{pmatrix}$} &
  \makebox[3.3cm]{$\begin{pmatrix}1&a&b&c&d\\0&1&0&b&a\\0&0&1&a&b\\0&0&0&1&0\\0&0&0&0&1\end{pmatrix}$} &
  \makebox[3.3cm]{$\begin{pmatrix}1&a&b&c&d\\0&1&0&b&-a\\0&0&1&a&b\\0&0&0&1&0\\0&0&0&0&1\end{pmatrix}$} &
  \makebox[3.3cm]{$\begin{pmatrix}1&0&0&0&d\\0&1&0&0&c\\0&0&1&a&b\\0&0&0&1&a\\0&0&0&0&1\end{pmatrix}$}
\end{array} \\[14pt]
\begin{array}{cccc}
  \makebox[3.3cm]{$N_{311}$} & \makebox[3.3cm]{$N^1_{312}$} & \makebox[3.3cm]{$N^2_{312}$} & \makebox[3.3cm]{$N^1_{313}$} \\[4pt]
  \makebox[3.3cm]{$\begin{pmatrix}1&a&b&c&d\\0&1&a&0&0\\0&0&1&0&0\\0&0&0&1&0\\0&0&0&0&1\end{pmatrix}$} &
  \makebox[3.3cm]{$\begin{pmatrix}1&0&0&0&d\\0&1&a&b&c\\0&0&1&0&a\\0&0&0&1&b\\0&0&0&0&1\end{pmatrix}$} &
  \makebox[3.3cm]{$\begin{pmatrix}1&0&0&0&d\\0&1&a&b&c\\0&0&1&0&a\\0&0&0&1&-b\\0&0&0&0&1\end{pmatrix}$} &
  \makebox[3.3cm]{$\begin{pmatrix}1&a&b&c&d\\0&1&0&a&0\\0&0&1&b&0\\0&0&0&1&0\\0&0&0&0&1\end{pmatrix}$}
\end{array} \\[14pt]
\begin{array}{cccc}
  \makebox[3.3cm]{$N^2_{313}$} & \makebox[3.3cm]{$N^1_{314}$} & \makebox[3.3cm]{$N^2_{314}$} & \makebox[3.3cm]{$N_{415}$} \\[4pt]
  \makebox[3.3cm]{$\begin{pmatrix}1&a&b&c&d\\0&1&0&a&0\\0&0&1&-b&0\\0&0&0&1&0\\0&0&0&0&1\end{pmatrix}$} &
  \makebox[3.3cm]{$\begin{pmatrix}1&a&b&c&d\\0&1&0&0&a\\0&0&1&0&b\\0&0&0&1&c\\0&0&0&0&1\end{pmatrix}$} &
  \makebox[3.3cm]{$\begin{pmatrix}1&a&b&c&d\\0&1&0&0&a\\0&0&1&0&-b\\0&0&0&1&-c\\0&0&0&0&1\end{pmatrix}$} &
  \makebox[3.3cm]{$\begin{pmatrix}1&a&b&c&d\\0&1&a&b&0\\0&0&1&a&0\\0&0&0&1&0\\0&0&0&0&1\end{pmatrix}$}
\end{array} \\[14pt]
\begin{array}{ccc}
  \makebox[3.3cm]{$N_{416}$} & \makebox[3.3cm]{$N_{417}$} & \makebox[3.3cm]{$N_{518}$} \\[4pt]
  \makebox[3.3cm]{$\begin{pmatrix}1&b&a&c&d\\0&1&0&0&0\\0&0&1&a&c\\0&0&0&1&a\\0&0&0&0&1\end{pmatrix}$} &
  \makebox[3.3cm]{$\begin{pmatrix}1&a&c&b&d\\0&1&a&0&c\\0&0&1&0&a\\0&0&0&1&b\\0&0&0&0&1\end{pmatrix}$} &
  \makebox[3.3cm]{$\begin{pmatrix}1&a&b&c&d\\0&1&a&b&c\\0&0&1&a&b\\0&0&0&1&a\\0&0&0&0&1\end{pmatrix}$} 
\end{array}
\end{align*}
\end{theorem}

The proof uses Galois cohomology and is in Appendix \ref{appendix:realforms} by Willem de Graaf. The groups
$N_{23}$ and $N_{24}$ have dimension $6$, the group $N_{35}$ has dimension $5$, and all
the others have dimension $4$. For each group the centralizer of its Lie algebra $\n$ in
$\gl(5,\R)$ is $\R I\oplus\n$, so these groups are indeed maximal abelian in
$\SL(5,\R)$.

\begin{remark}\label{rem:real-forms}
The pairs of real forms can be told apart by a quadratic form, as in the case of
$\GL(4,\R)$ \cite{WZ84, diatta2022classification}. For a general element $X$ of the Lie algebra, $X^2$ has a
single nonzero entry in the following cases (see Appendix \ref{appendix:nonunipotent} for $\w6^\pm$):
\begin{center}
\renewcommand{\arraystretch}{1.2}
\begin{tabular}{lll}
algebra & nonzero entry of $X^2$ & signatures of the two real forms \\ \hline
$\w6^{\pm}\subset\gl(4,\R)$ & $a^2\pm b^2$ & $(2,0)$ and $(1,1)$ \\
$N_{312}^{1,2}$ & $a^2\pm b^2$ & $(2,0)$ and $(1,1)$ \\
$N_{313}^{1,2}$ & $a^2\pm b^2$ & $(2,0)$ and $(1,1)$ \\
$N_{314}^{1,2}$ & $a^2+b^2+c^2$, resp.\ $a^2-b^2-c^2$ & $(3,0)$ and $(1,2)$
\end{tabular}
\end{center}
The form is only determined up to a nonzero scalar, so for $N_{314}$ the signatures
$(1,2)$ and $(2,1)$ give the same class; this is why $N_{314}$ has exactly two real
forms rather than three. For $N_{38}$ and $N_{39}$ the real forms differ in the same way, by the sign of $a^2$
in the entry $a^2\pm b^2$ of $X^2$; that there are exactly two real forms in each of
these cases is proved in Appendix \ref{appendix:realforms} by Willem de Graaf.
\end{remark}

We now determine which unipotent groups contain the group $B$ of a rank $4$ projective
structure. For unipotent groups the open-orbit test becomes a condition on a single
column. We use the following standard fact.
 
\begin{lemma}[{Kostant--Rosenlicht: \cite[ Theorem 2]{rosenlicht1961quotient}, \cite[Proposition 4.10]{borel2012linear} }]\label{lem:KR}
Let $U$ be a unipotent group acting on $\R^n$ by affine transformations whose linear
parts are unipotent. Then every $U$-orbit is closed.
\end{lemma}

\begin{proof} First let $U$ act over $\mathbb C$.
Let $O=U\cdot q$. The Zariski closure $\overline O$ is irreducible and $O$ is open in it,
so $Z=\overline O\setminus O$ is Zariski closed and $U$-invariant. Suppose
$Z\neq\emptyset$ and pick a nonzero polynomial $f$ on $\overline O$ vanishing on $Z$. The
$U$-translates of $f$ span a finite-dimensional $U$-invariant space of such functions,
and since $U$ is unipotent it contains a nonzero $U$-fixed vector $F$. Then $F$ is
constant on $O$, hence on $\overline O$, and vanishes on $Z$, so $F\equiv0$, a
contradiction. So $O$ is Zariski closed, hence closed. For the real group, the stabilizer $U_q$ is unipotent, hence has trivial Galois cohomology, so the real orbit $U(\mathbb R)\cdot q$ is the set of real points of the complex orbit $U(\mathbb C)\cdot q$. The complex orbit is closed, hence so is its intersection with $\mathbb R^n$.
\end{proof}

\begin{theorem}\label{thm:unipotent-rank4}
Let $N$ be one of the groups of Theorem~\ref{thm:unipotent-list}, with Lie algebra $\n$.
The following are equivalent:
\begin{enumerate}
  \item $N$ contains the group $B$ of a rank 4 projective structure on $\T^4$;
  \item $N$ contains the holonomy of a complete affine structure on $\T^4$;
  \item the vectors $Xe_5$, $X\in\n$, span $\R e_1\oplus\dots\oplus\R e_4$; that is,
  the last column of the general element of $N$ contains four independent parameters.
\end{enumerate}
These conditions hold exactly for: $
  N_{21},\ N_{36},\ N^1_{38},\ N^2_{38},\ N_{310},\ N^1_{312},\ N^2_{312},\
  N^1_{314},\ N^2_{314},\ N_{417},\ N_{518}.
$
\end{theorem}

\begin{proof}
Every $X\in\n$ is strictly upper triangular, so $N$ preserves the hyperplane $\{z=0\}$
and each of the two open hemispheres $U_\pm=\{\pm z>0\}$. In the affine chart $z=1$,
which identifies $U_+$ with $\R^4$, the group $N$ acts by affine transformations with
unipotent linear parts.
 
\emph{Step 1: $N$ has an open orbit if and only if (3) holds, and then the open orbits
are $U_+$ and $U_-$.} If (3) holds, then $\n e_5+\R e_5=\R^5$, so the orbit of $[e_5]$
is open by Lemma~\ref{lem:C-criterion}. Conversely, let $O$ be an open orbit. It meets $U_+$ or $U_-$,
say $U_+$ (the other case is symmetric), and since both sets are $N$-invariant,
$O\subset U_+$. By Lemma~\ref{lem:KR}, $O$ is closed in $U_+$, and $U_+$ is connected,
so $O=U_+$. In particular $[e_5]\in O$, so $\n e_5+\R e_5=\R^5$ by Lemma~\ref{lem:C-criterion}; since
every $Xe_5$ lies in $\R e_1\oplus\dots\oplus\R e_4$, this is (3). The same argument
shows that any open orbit of any connected subgroup of $N$ is $U_+$ or $U_-$.
 
\emph{(1)$\Leftrightarrow$(3).} By Lemma~\ref{lem:reduction}, (1) holds if and only if
$N$ has an open orbit, which by Step~1 is equivalent to (3).
 
\emph{(3)$\Rightarrow$(2).} By Lemma~\ref{lem:C-criterion}(3) some $4$-dimensional subgroup $B\subset N$
has an open orbit, which by Step~1 is $U_+$ (up to sign). Now $B\cong\R^4$ and
$U_+\cong\R^4$ are simply connected, so the orbit map $B\to U_+$, a covering map, is a
diffeomorphism: $B$ acts simply transitively on the affine chart $U_+$. For a lattice
$\Gamma\cong\Z^4$ in $B$, the quotient $\Gamma\backslash U_+$ is a complete affine torus
with holonomy $\Gamma\subset N$.
 
\emph{(2)$\Rightarrow$(3).} Suppose a complete affine structure on $T^4$ has holonomy
$\Gamma\cong\Z^4$ contained in $N$. Its developing map is a diffeomorphism onto an
affine chart $A\subset\Sph^4$, and $\Gamma$ acts on $A$ freely, properly and
cocompactly. Let $B=\exp(\operatorname{span}\log\Gamma)\subset N$, a connected
unipotent group of dimension at most $4$ containing $\Gamma$. It lies in the Zariski
closure of $\Gamma$, so it preserves the hyperplane $\Sph^4\setminus A$, and it acts on
$A\cong\R^4$ by affine maps with unipotent linear parts. Fix $q\in A$. The orbit
$B\cdot q$ is closed in $A$ by Lemma~\ref{lem:KR} and $\Gamma$-invariant, so $\Gamma$
acts on it properly and cocompactly. The stabilizer $B_q$ is a connected algebraic
subgroup, so $B\cdot q\cong B/B_q\cong\R^k$ with $k\le4$, and $\Gamma$ acts on it by
translations. A free, proper action of $\Z^4$ by translations on $\R^k$ is an
injective homomorphism onto a discrete subgroup, which forces $k\ge4$. Hence $k=4$,
the orbit $B\cdot q$ is open, and $N$ has an open orbit. By Step~1, (3) holds.

Finally, (3) can be read off the matrices in Theorem \ref{thm:unipotent-list} by reading off when the last column has 4 independent parameters.

\end{proof}

\begin{remark}
Equivalently, for each of the eleven groups in Theorem~\ref{thm:unipotent-rank4} the
determinant $\delta$ of Lemma~\ref{lem:C-criterion} equals $\pm z^5$, and its open orbits are the two
hemispheres $U_\pm$; for the other twelve groups the rank $r(p)$ is at most $4$ for
every $p$. This is consistent with the theorem of Fried--Goldman--Hirsch
\cite{fried1981affine} Theorem A that a compact affine manifold with unipotent holonomy is complete.
\end{remark}

\subsection{The non-unipotent groups}
 
The remaining maximal abelian subgroups have more than one block (possibly complex). They are listed in Appendix \ref{appendix:nonunipotent}: by Lemma~\ref{lem:real-blocks} there are $24$ of
them, namely nine of block type $(4,1)$, including $A(\w8,\bd)$ with the complex block
$\w8$, and fifteen with smaller blocks. By Lemma~\ref{lem:reduction} and
Proposition~\ref{prop:C-blocks}, such a group contains the group $B$ of a rank 4 projective structure if and only if each of its blocks is good.
Exactly $18$ of the $24$ groups have this property; they are marked in
Table \ref{table:nonunip}. As in the two half plane structures in Example \ref{ex:T2}, some of them can be dehomogenized in different ways,
giving several affine structures.

\bibliography{bibliography.bib}

@article{baues2008prehomogeneous,
  title={Prehomogeneous {A}ffine {R}epresentations and {F}lat {P}seudo-{R}iemannian {M}anifolds},
  author={Baues, Oliver},
  journal={arXiv preprint arXiv:0809.0824},
  year={2008}
}

@article{benoist1994nilvarietes,
  title={Nilvari{\'e}t{\'e}s {P}rojectives},
  author={Benoist, Yves},
  journal={Commentarii Mathematici Helvetici},
  volume={69},
  number={3},
  pages={447--473},
  year={1994},
  publisher={Zurich: Orell Fussli Verlag, 1929-}
}

@article{Ben00,
  title={Tores {A}ffines},
  author={Benoist, Yves},
  journal={Contemporary Mathematics},
  volume={262},
  pages={1--38},
  year={2000},
  publisher={Providence, RI: American Mathematical Society}
}

@book{choi2001decomposition,
  title={The {D}ecomposition and {C}lassification of {R}adiant {A}ffine 3-{M}anifolds},
  author={Choi, Suhyoung},
  volume={730},
  year={2001},
  series={Mem. Amer. Math. Soc. 730},
  publisher={American Mathematical Soc.}
}

@article{diatta2022classification,
  title={Classification of 2-{S}olvable {F}robenius Lie {A}lgebras {G}iven by {N}onderogatory {M}atrices},
  author={Diatta, Andre and Manga, Bakary and Mbaye, Ameth},
  journal={arXiv preprint arXiv:2209.05595},
  year={2022}
}

@article{fried1981affine,
  title={Affine {M}anifolds with {N}ilpotent {H}olonomy},
  author={Fried, David and Goldman, William and Hirsch, Morris W},
  journal={Commentarii Mathematici Helvetici},
  volume={56},
  pages={487--523},
  year={1981},
  publisher={Birkh{\"a}user-Verlag}
}

@article{furness1975locally,
  title={Locally {S}ymmetric {S}paces},
  author={Furness, PMD and Arrowsmith, DK},
  journal={Journal of the London Mathematical Society},
  volume={2},
  number={4},
  pages={487--499},
  year={1975},
  publisher={Oxford University Press}
}

@book{goldman2022geometric,
  title={Geometric {S}tructures on {M}anifolds},
  author={Goldman, William M},
  volume={227},
  year={2022},
  publisher={American Mathematical Society}
}

@article{NY74,
  title={The {A}ffine {S}tructures on the {R}eal Two-{T}orus. {I}},
  author={Nagano, Tadashi and Yagi, Katsumi},
journal={Osaka J. Math},
  year={1974}}

@book {suprunenko1968commutative,
    AUTHOR = {Suprunenko, D. A. and Tyshkevich, R. I.},
     TITLE = {Commutative {M}atrices ({P}erestanovochnye {M}atritsy)},
   EDITION = {Second},
 PUBLISHER = {\`Editorial URSS, Moscow},
      YEAR = {2003},
     PAGES = {104},
      ISBN = {5-354-00387-3},
   MRCLASS = {16S50 (15A27 16N40)},
  MRNUMBER = {2118458},
MRREVIEWER = {Mykola\ Ya.\ Komarnitski\u i},
}

@article{WZ84,
  title={The {S}tructure of {M}aximal {A}belian {S}ubalgebras of {C}lassical {L}ie and {J}ordan algebras, dans  XIIIth International Colloquium on Group Theoretical Methods in Physics (College Park, Md., 1984) },
  author={Winternitz, P and Zassenhaus, H},
  journal={World Sci, Publishing},
  year={1984}
}

@book{serre_gal,
    AUTHOR = {Serre, Jean-Pierre},
     TITLE = {Galois {C}ohomology},
    SERIES = {Springer Monographs in Mathematics},
   EDITION = {English},
      NOTE = {Translated from the French by Patrick Ion and revised by the
              author},
 PUBLISHER = {Springer-Verlag, Berlin},
      YEAR = {2002}
}

@book{berhuy,
    AUTHOR = {Berhuy, Gr\'egory},
     TITLE = {An Introduction to {G}alois {C}ohomology and its {A}pplications},
    SERIES = {London Mathematical Society Lecture Note Series},
    VOLUME = {377},
      NOTE = {With a foreword by Jean-Pierre Tignol},
 PUBLISHER = {Cambridge University Press, Cambridge},
      YEAR = {2010}
}

@article {borwdg,
    AUTHOR = {Borovoi, Mikhail and de Graaf, Willem A.},
     TITLE = {Computing {G}alois {C}ohomology of a {R}eal {L}inear {A}lgebraic
              {G}roup},
   JOURNAL = {J. Lond. Math. Soc. (2)},
  FJOURNAL = {Journal of the London Mathematical Society. Second Series},
    VOLUME = {109},
      YEAR = {2024},
    NUMBER = {5},
     PAGES = {Paper No. e12906, 53}
}

@incollection{thurston2014three,
  title={Three-{D}imensional {G}eometry and {T}opology, volume 1:(pms-35)},
  author={Thurston, William P},
  booktitle={Three-Dimensional Geometry and Topology, Volume 1},
  year={2014},
  publisher={Princeton university press}
}

@phdthesis{smillie1977affinely,
  title={Affinely {F}lat {M}anifolds.},
  author={Smillie, John David},
  year={1977},
  school={The University of Chicago}
}

@article{rosenlicht1961quotient,
  title={On {Q}uotient {V}arieties and the {A}ffine {E}mbedding of {C}ertain {H}omogeneous {S}paces},
  author={Rosenlicht, Maxwell},
  journal={Transactions of the American Mathematical Society},
  volume={101},
  number={2},
  pages={211--223},
  year={1961}
}

@book{borel2012linear,
  title={Linear {A}lgebraic {G}roups},
  author={Borel, Armand},
  year={2012},
  publisher={GTM 126}
}

@article{mal1949class,
  title={On a {C}lass of {H}omogeneous {S}paces},
  author={Mal'tsev, Anatolii Ivanovich},
  journal={Izvestiya Akademii Nauk SSSR},
  volume={13},
  number={1},
  pages={9--32},
  year={1949},
  publisher={Russian Academy of Sciences, Steklov Mathematical Institute of Russian~…}
}

@book{raghunathan1972discrete,
  title={ Discrete {S}ubgroups of {L}ie {G}roups},
  author={Raghunathan, Madabusi Santanam},
  publisher={Springer},
  Series={Ergebnisse der Mathematik und ihrer Grenzgebiete 68},
  year={1972}
}
\bibliographystyle{alpha}

\appendix

\section{ The real forms of the maximal abelian nilpotent matrix algebras, by Willem de Graaf}\label{appendix:realforms}

In this appendix we consider maximal abelian nilpotent Lie subalgebras (MANSs) of $\gl(n,F)$ consisting of nilpotent matrices. Here we are interested in the fields $F=\R$ and $F=\C$. 
Secondly, we want the subalgebras up to the conjugation action by the group $\GL(n,F)$. Thirdly, we are mainly concerned with $n=5$. 

A first remark is that such a subalgebra can be conjugated into the strictly upper triangular matrices. For $n\leq 6$ the MANSs in $\gl(n,\C)$ have been classified in 
\cite{suprunenko1968commutative}.  The paper \cite{WZ84} contains a classification of the MANSs in $\gl(4,\R)$. Here we focus on the MANSs in $\gl(5,\R)$.
The result is listed in Theorem \ref{thm:unipotent-list}. We have obtained it using Galois cohomology. 

Now we have several subsections. The first one gives a very short introduction into Galois cohomology. The second subsection outlines our method to classify the real forms.
The third subsection gives details on the classification of the real MANSs in $\gl(5,\R)$.

\subsection{Galois cohomology}

Here we deal with Galois cohomology relative to the Galois group $\Gamma = \mathrm{Gal}(\C/\R)$. 
For a general treatment we refer to the books by Berhuy and Serre (\cite{berhuy,serre_gal}).

Write $\Gamma = \{1,\gamma\}$. Let $\cG$ be a group on which $\Gamma$ acts by automorphisms. 
Let $X$ be a set on which $\Gamma$ acts. Let $\sigma : \cG\to\cG$, $\sigma : X\to X$
be the maps corresponding to the action of $\gamma$. 
We suppose that $\cG$ acts 
$\Gamma$-equivariantly on $X$, that is, $\sigma(g) \cdot \sigma(x) = \sigma (g\cdot x)$.

An element $g\in \cG$ is called a {\em cocycle} if $g\sigma(g) = 1$. We remark that since we 
use the Galois group $\Gamma$, this is equivalent to, but not the same as, the usual definition. The set of all cocycles is denoted by $Z^1(\cG,\sigma)$, or
just by $Z^1 \cG$ if $\sigma$ is clear. We define an equivalence relation
$\sim$ on $Z^1(\cG,\sigma)$ by $g_1 \sim g_2$ if there is 
an $h\in \cG$ with $h^{-1} g_1  \sigma(h)=g_2$. The set of equivalence
classes in $Z^1(\cG,\sigma)$ is called the first Galois cohomology set of $\cG$ and
denoted $H^1 (\cG,\sigma)$, or also just $H^1 \cG$ if $\sigma$ is clear.
For a cocycle $g\in \cG$ we denote its equivalence class in $H^1\cG$ by $[g]$.

By $\cG^\sigma$, $X^\sigma$ we denote the elements of $\cG$ and $X$ respectively that are 
fixed under $\sigma$. Let $x_0\in X^\sigma$ and let $Z_{\cG}(x_0) = \{ g\in \cG \mid 
g\cdot x_0=x_0\}$ be its stabilizer in $\cG$. We have a natural map $i_* : H^1 Z_{\cG}(x_0)\to 
H^1 \cG$, mapping the class of a cocycle in $H^1 Z_{\cG}(x_0)$ to its class in $H^1 \cG$.
By definition, its kernel is the set of all classes that are sent to the trivial class, i.e.,
$$\ker i_* = \{ [c]\in H^1 Z_{\cG}(x_0) \mid c \sim 1 \text{ in }
Z^1(\cG,\sigma)\}.$$
We now have the following theorem for whose proof we refer to 
\cite[Section I.5.4, Corollary 1 of Proposition 36]{serre_gal}. 

\begin{theorem}\label{thm:galco}
Let the notation be as above. Let $Y=\cG\cdot x_0$ be the orbit of $x_0$.  The orbits of 
$\cG^\sigma$ in $Y^\sigma$ are in bijection with $\ker i_*$. The bijection goes as follows:
let $[c]\in \ker i_*$, then there is a $g\in \cG$ with $g^{-1} \sigma(g) = c$ and $[c]$ 
corresponds to the $\cG^\sigma$-orbit of $g\cdot x_0$.
\end{theorem}

\subsection{The classification of real forms}

Let $N\subset \gl(n,\C)$ be a MANS given by a real basis. We want to classify, up to conjugation by $\GL(n,\R)$, all MANSs in $\gl(n,\R)$ that are 
$\GL(n,\C)$-conjugate to $N$. For this we use Theorem \ref{thm:galco}. In the notation of the previous subsection we let $\cG=\GL(n,\C)$ and
$\sigma(g) = \bar g$ for $g\in \GL(n,\C)$ (where $\bar g$ is given by complex conjugation of the matrix entries. We let $X$ be the set of subalgebras of 
$\gl(n,\C)$ and $\sigma(U) = \{ \bar u \mid u\in U\}$ for $U$ a subalgebra of $\gl(n,\C)$ and $\bar u$ is again given by entry-wise complex conjugation.
Then $\cG^\sigma = \GL(n,\R)$ and $X^\sigma$ consists of all subalgebras that have a real basis (i.e., whose elements are fixed under complex conjugation).
The stabilizer of $N$ is its normalizer,
$$Z_{\GL(n,\C)}(N) = \{ g\in \GL(n,\C) \mid gNg^{-1} = N\}.$$
We use the algorithms from \cite{borwdg} to compute the Galois cohomology for this group. These algorithms take as input a basis of the Lie algebra
of $Z_{\GL(n,\C)}(N)$ along with one element of each component of that group. The Lie algebra of $Z_{\GL(n,\C)}(N)$ is 
$$\mathfrak{z}_{\gl(n,\C)} = \{ x\in \gl(n,\C) \mid [x,N] \subset N\}$$
which can easily be determined by linear algebra computations. In order to determine the components let $u_1,\ldots,u_r$ denote a basis of $N$. Then 
$g\in \GL(n,\C)$ lies in $Z_{\GL(n,\C)}(N)$ if and only if there are $a_{ij}\in \C$ such that 
\begin{equation}\label{eq:eqns}
gu_i = (a_{i1} u_1+\cdots +a_{ir} u_r) g\text{ for } 1\leq i\leq r.
\end{equation}
Writing $g=(g_{kl})$ this is equivalent to a number of polynomial equations in the indeterminates $g_{kl}$, $a_{ij}$. We study these equations using 
Gr\"obner basis techniques and algorithms. In all cases of interest in this paper we were able to find the components this way. 

It is a well known fact that $H^1 \GL(n,\C)$ just consists of the trivial class. So, 
using the notation of Theorem \ref{thm:galco}, all $[c]\in H^1 Z_{\GL(n,\C)}(N)$ lie in $\ker i_*$. Secondly, 
\cite{borwdg} also has algorithms to compute $g\in \GL(n,\C)$ with $g^{-1}\bar g = c$. 

\begin{example}\label{ex:warmup}
By $e_{ij}$ we denote the $4\times 4$-matrix with 1 on position $(i,j)$ and zeros elsewhere. 
Let $N$ be the MANS in $\gl(4,\C)$ with the basis $u_1= e_{12}+e_{24}$, $u_2=e_{13}+e_{34}$, $u_3=e_{14}$. We consider the polynomial equations deriving from
\eqref{eq:eqns}. We add two indeterminates $d_1,d_2$ and the polynomials equations $d_1\det(g_{kl})-1=0$, $d_2\det(a_{ij}) -1=0$. This way we obtain a number of 
polynomial equations of the form $p=0$. We consider the ideal $I$ generated by all left hand sides $p$. We compute a Gr\"obner basis of $I$ with respect to the 
lexicographical order where the $g_{kl}$ are smallest. The Gr\"obner basis has 661 elements
and contains the following polynomials just in the $g_{kl}$:
\begin{align*}
& g_{11}g_{44} - g_{32}^2 - g_{33}^2,\\
& g_{22}^2 - g_{33}^2,\\
& g_{22}g_{23} + g_{32}g_{33},\\
& g_{22}g_{32} + g_{23}g_{33},\\
& g_{23}^2 - g_{32}^2,\\
& g_{21}, g_{31}, g_{41}, g_{42}, g_{43}.
\end{align*}
A small calculation shows that set of matrices satisfying these equations consists of two types:
$$\begin{pmatrix} g_{11} & g_{12} & g_{13} & g_{14}\\ 0 & g_{22} & g_{23} & g_{24}\\ 0 & -g_{23} & g_{22} & g_{24}\\ 0&0&0&g_{44}\end{pmatrix} \text{ and }
\begin{pmatrix} g_{11} & g_{12} & g_{13} & g_{14}\\ 0 & g_{22} & g_{23} & g_{24}\\ 0 & g_{23} & -g_{22} & g_{24}\\ 0&0&0&g_{44}\end{pmatrix}
\text{ with } g_{11} = g_{44}^{-1}(g_{22}^2+g_{23}^2), g_{22}^2+g_{23}^2\neq 0.$$
Conversely, a short calculation shows that all matrices of this form stabilize the algebra $N$. 
So the stabilizer has two components and an element of the non-identity component is $g_0 = \diag(1,1,-1,1)$. 
A computation with the implementation of the algorithms of \cite{borwdg} shows that $H^1 Z_{\GL(n,\C)}(N)$ has two elements, and that $g_0$ is a 
representative of the nontrivial class. We have $g_0 = h^{-1} \bar h$ with $h=\diag(1,1,i,1)$. We have
$hu_1h^{-1} = u_1$, $hu_2h^{-1} = -i(e_{13}-e_{34})$, $hu_3h^{-1} = e_{14}$. So we see that $hNh^{-1}$ has the real basis 
$u_1, e_{13}-e_{34}, u_3$. The conclusion is that $N$ has two real forms. 
\end{example}

\subsection{Results for \texorpdfstring{$n=5$}{n=5}}

Now we let $n=5$. According to \cite{suprunenko1968commutative} there are 18 classes of MANSs in $\gl(5,\C)$. They are listed in Theorem \ref{thm:unipotent-list}.
This theorem contains a list of unipotent groups. The corresponding MANSs are obtained by simply replacing the 1's on the diagonal by 0's. If in Theorem
\ref{thm:unipotent-list} the symbol for a MANS has an upper index then for the complex classification only upper index 1 has to be considered. The other upper indices 
then give the different real forms.
We have carried out the procedure outlined above for all cases. Table \ref{tab:galco} summarizes our findings. In this table the first column has the label of 
the MANS from Theorem \ref{thm:unipotent-list}. The second column gives the structure of the identity component of the stabilizer of $N$. Here $T_k$ denotes a torus of 
dimension $k$. Furthermore $A_1$ denotes a semisimple group whose root system is of type $A_1$. Finally, $U_k$ denotes a unipotent group of dimension $k$. 
The third and fourth column display, respectively,
the order of the component group and the order of the Galois cohomology set $H^1 Z_{\GL(5,\C)}(N)$. 

\begin{table}[htb]\caption{Stabilizers and Galois cohomology of the MANSs in $\gl(5,\C)$.}\label{tab:galco} 
\begin{center}
\begin{tabular}{|c|r|c|c|}
\hline  
$N$ & structure of $Z=Z_{\GL(5,\C)}(N)^\circ$ & $K$ & $|H^1 Z|$ \\
\hline
$N_{21}$ &  $\GL(4,\C)\times T_1 \ltimes U_4$ & 1 & 1 \\
$N_{22}$ &  $\GL(4,\C)\times T_1 \ltimes U_4$ & 1 & 1 \\
$N_{23}$ &  $\GL(3,\C)\times \GL(2,\C)\ltimes U_6$ & 1 & 1\\
$N_{24}$ &  $\GL(3,\C)\times \GL(2,\C)\ltimes U_6$ & 1 & 1\\
$N_{35}$ &  $T_4\ltimes U_{10}$ & 1 & 1 \\
$N_{36}$ &  $T_3\ltimes U_9$ & 1 & 1 \\
$N_{37}$ &  $T_3\ltimes U_9$ & 1 & 1 \\
$N_{38}^1$ & $T_3\ltimes U_8$ & $C_2$ & 2 \\
$N_{39}^1$ & $T_3\ltimes U_8$ & $C_2$ & 2 \\
$N_{310}$ & $\GL(2,\C)\times T_2 \ltimes U_9$ & 1 & 1\\
$N_{311}$ & $\GL(2,\C)\times T_2 \ltimes U_9$ & 1 & 1\\
$N_{312}^1$ & $T_4\ltimes U_9$ & $C_2$ & 2 \\
$N_{313}^1$ & $T_4\ltimes U_9$ & $C_2$ & 2 \\
$N_{314}^1$ & $A_1T_2 \ltimes U_7$ & 1 & 2 \\
$N_{415}$ & $T_3\ltimes U_8$ & 1 & 1 \\
$N_{416}$ & $T_3\ltimes U_8$ & 1 & 1 \\
$N_{417}$ & $T_2\ltimes U_8$ & 1 & 1 \\
$N_{518}$ &  $T_2\ltimes U_7$ & 1 & 1\\
\hline
\end{tabular}
\end{center}
\end{table}

We give some details for the cases with nontrivial Galois cohomology. 

Let $N=N_{38}^1$. Then the component group of $Z_{\GL(5,\C)}(N)$ is generated by $d=\diag(1,-1,1,-1,1)$. We have that $H^1 Z_{\GL(5,\C)}(N)$ consists of the
trivial class and of the class of $d$. Let $h=\diag(1,i,1,i,1)$. Then $h^{-1}\bar h = d$. We have that $hNh^{-1}$ consists of 
$$\begin{pmatrix} 0 & 0 & b & -a & d \\ 0 & 0 & a & b & c \\
  0 & 0 & 0 & 0 & b\\ 0 & 0 & 0 & 0 & a \\ 0 & 0 & 0 & 0 & 0\end{pmatrix}.$$
The algebra $N_{39}^1$ is conjugate to the transpose of $N_{38}^1$. So the second real form of $N_{39}^1$ is conjugate to the transpose of $N_{38}^2$.

Let $N=N_{312}^1$. Then the component group of $Z_{\GL(5,\C)}(N)$ is generated by $d=\diag(1,1,1,-1,1)$. We have that $H^1 Z_{\GL(5,\C)}(N)$ consists of the
trivial class and of the class of $d$. Let $h=\diag(1,1,1,i,1)$. Then $h^{-1}\bar h = d$. We have that $hNh^{-1}$ consists of 
 $$\begin{pmatrix} 0 & 0 & 0 & 0 & d\\ 0&0&a&b&c\\ 0&0&0&0&a\\
    0&0&0&0&-b\\ 0&0&0&0&0 \end{pmatrix}.$$
   The algebra $N_{313}^1$ is conjugate to the transpose of $N_{312}^1$. So the second real form of $N_{313}^1$ is conjugate to the transpose of $N_{312}^2$.

Let $N=N_{314}^1$. Then  $Z_{\GL(5,\C)}(N)$  is connected.  The set  $H^1 Z_{\GL(5,\C)}(N)$ consists of the
trivial class and of the class of $d = \diag(1,1,-1,-1,1)$. Let $h=\diag(1,1,i,i,1)$. Then $h^{-1}\bar h = d$. We have that $hNh^{-1}$ consists of 
 $$\begin{pmatrix} 0&a&b&c&d\\ 0&0&0&0&a\\ 0&0&0&0&-b\\ 0&0&0&0&-c\\
    0&0&0&0&0\end{pmatrix}.$$

\section{Rank 2 structures on \texorpdfstring{$\T^3$}{T3}}\label{appendix:T3}

In this appendix we apply the methods of Section \ref{sec:bricks}, with $d=3$, to the rank $2$ projective
structures on $\T^3$, and work through the exercises of \cite{Ben00} Section~4.6. Throughout,
$N\subset\mathrm{SL}(4,\R)$ is one of the groups $N_{i_1,\dots,i_k}$ of Proposition \ref{prop:open-orbits}, with
upper triangular blocks, and $B\subset N$ is the $2$-dimensional abelian group of a rank
$2$ structure, with Lie algebra $\mathfrak b$. Points of $\Sph^3$ are written $[x:y:z:w]$. As
in Section \ref{sec:rank3}, we determine the normal forms of $\mathfrak b$ and freeness from the $3\times3$
minors of the matrix $(b_1p,b_2p,p)$, find an invariant function, locate the faces with
Lemma \ref{lem:limits} and Corollary \ref{cor:parity}, and read off the bricks and words from Lemmas \ref{lem:bricks} and \ref{lem:constant-faces}.
Conversely, every admissible periodic word is realized, with holonomy group
$\langle\Gamma_0,g\rangle$, by Lemma \ref{lem:realize}.
 
Up to conjugation $N$ is one of $N_4$, $N_{3,1}$, $N_{2,2}$, $N_{2,1,1}$, $N_{1,1,1,1}$. The
unipotent case $N_4$ is \cite{Ben00} Exercise~4.6(b): rank $2$ structures exist, $V\cong \T^3$,
and $\Gamma_0$ is a lattice in one of two explicit $2$-dimensional subgroups of $N_4$. The
diagonal case is treated exactly as in Section \ref{sec:diag}: a $2$-dimensional subalgebra of the
traceless diagonal matrices is the kernel of a nonzero $\kappa\in\R^4$ with
$\sum_i\kappa_i=0$; the function $\sum_i\kappa_i\log x_i$ is an invariant function; $B$ is free
on the face $\{x_i=0\}$ if and only if $\kappa_i\neq0$; and the conclusions of
Proposition \ref{prop:diagonal} hold verbatim, with $n_++n_-\le4$. The three remaining groups are the
subject of \cite{Ben00} Exercises~4.13--4.15. Compared with the statements there, we find:
\begin{itemize}
  \item for $N_{3,1}$ there are rank $2$ projective structures, none of which is affine
  (Proposition~\ref{prop:B31});
  \item for $N_{2,1,1}$, which faces lie at the same end depends on both generators of
  $\mathfrak b$, through the sign of the number $K$ of Section \ref{sec:N211};
  \item for $N_{2,1,1}$, when $\lambda$ lies between $\mu$ and $\nu$ (including the
  boundary cases), only the brick $\mathcal O\cup \mathcal O_2^\pm$ occurs in a torus (Lemma \ref{lem:constant-faces}).
\end{itemize}

\begin{remark}[Structures occurring for two groups $N$]\label{rem:overlaps-B}
As in \Cref{sec:overlaps}, a structure on $\T^3$ occurs for two non-conjugate groups
$N$ and $N'$ exactly when, after conjugation, $\mathfrak b\subset\mathfrak n\cap\mathfrak n'$ and the
holonomy lies in $N'$; this is detected by merging blocks of $N$ on which $\mathfrak b$ has the
same weight, or splitting blocks on which $\mathfrak b$ acts decomposably. Running through
Propositions~\ref{prop:B31}, \ref{prop:B22} and~\ref{prop:B211}, and the unipotent and
diagonal cases, this happens only in two families:
{\sloppy
\begin{enumerate}
  \item[(I)] $N_{2,2}$, case (a) with $a=0$, and $N_{2,1,1}$, case (a) with $\mu=\nu$. Here
  $\mathfrak b=\langle\diag(1,1,-1,-1),\ E_{12}\rangle$ and the developing image is
  $\{y>0,\ (z,w)\neq(0,0)\}$. There is one brick for each group: for $N_{2,2}$ it is
  $\mathcal O\cup\mathcal O_4^+\cup\mathcal O_4^-$, with chain of period $2$; for $N_{2,1,1}$ it is
  $\mathcal O\cup\mathcal O_3\cup\mathcal O_4$, with chain of period $4$.
  \item[(II)] $N_{2,1,1}$ with $\mathfrak b=\mathfrak a$, and $N_{1,1,1,1}$ with $\kappa$
  proportional to $(1,-1,0,0)$. Here $\mathfrak b=\{\diag(s,s,t,u):2s+t+u=0\}$ and the
  developing image is $\{z>0,\ w>0,\ (x,y)\neq(0,0)\}$. Again there is one brick for
  each group: $\mathcal O\cup\mathcal O_2^+\cup\mathcal O_2^-$ with period $2$ for $N_{2,1,1}$,
  and one with period $4$ for $N_{1,1,1,1}$.
\end{enumerate}}
In all other cases the weights of $\mathfrak b$ on the blocks of $N$ are distinct and
$\mathfrak b$ acts indecomposably on each block. For $N_{3,1}$ this holds because $\mathfrak b$
contains a regular nilpotent element of the $3$-block; for $N_{2,2}$ in case (b), because
$X=E_{12}+bE_{34}$ with $b\neq0$; for $N_{2,1,1}$ with $\mathfrak b\neq\mathfrak a$, because
$\mathfrak b$ contains an element with nilpotent part $E_{12}$, and the weights of the two
singleton blocks coincide only if $\mu=\nu$ and $b=c$, which forces $K=0$, the case of
family~(I); and for the two
subgroups of $N_4$ in \cite{Ben00} Exercise~4.6(b), because the commutant of $\mathfrak b$ in
$\mathfrak{gl}(4,\R)$ contains no nontrivial idempotent.
 
In both families, the torus is a single-brick, affine structure, and as in
\Cref{sec:overlaps} it occurs for both groups only when its gluing element lies
in $N\cap N'$, for instance when $g\in B$; its finite covers, given by the words $w^m$ with
$g^m$, then occur for both groups as well. In particular, every rank $2$ structure on
$\T^3$ with at least two bricks determines $N$ up to conjugation.
\end{remark}
 
\subsection{\texorpdfstring{$N_{3,1}$}{N31} (Exercise 4.13)}
 
Here $\mathfrak n=\R D\oplus\mathfrak u$ with $D=\diag(1,1,1,-3)$ and
$\mathfrak u=\langle E_{12},E_{13},E_{23}\rangle$, a Heisenberg algebra with center $E_{13}$.
The open orbit $\mathcal O=\{z>0,\ w>0\}$ has the faces
\[
  \mathcal O_3^\pm=\{z=0,\ \pm y>0,\ w>0\},\qquad \mathcal O_4=\{w=0,\ z>0\}.
\]
 
\emph{Normal form.} If $\mathfrak b\subset\mathfrak u$, then $\mathfrak b$ contains the center $E_{13}$,
which fixes $\mathcal O_3^\pm$ pointwise, so only $\mathcal O_4$ is available and there is no brick. Hence
$\mathfrak b$ is spanned by $D+u_1$ and an element $Y\in\mathfrak u$. An element with a nonzero
$D$-coefficient has no eigenray on $\mathcal O$ or on its faces, since it would have the two
eigenvalues $1$ and $-3$ there. On $\mathcal O_3^\pm$ we have $Yp=Y_{12}\,y\,e_1$, so $B$ is free on
$\mathcal O_3^\pm$ if and only if $Y_{12}\neq0$. Normalizing $Y_{12}=1$ and using the abelian
condition, we obtain
\[
  X=D+bE_{13},\qquad Y=E_{12}+eE_{13}+fE_{23}.
\]
The minors include $-4fwz^2$ on $\mathcal O$, $-4wy^2$ on $\mathcal O_3^\pm$ and $bfz^3$ on $\mathcal O_4$. So $B$ is
free on $\mathcal O$ if and only if $f\neq0$, always free on $\mathcal O_3^\pm$, and free on $\mathcal O_4$ if and only
if $bf\neq0$.
 
\emph{Invariant function.} The function
\[
  f_{3,1}=\frac{y^2}{z^2}+2e\frac yz-2f\frac xz-\frac12bf\log\frac wz
\]
is $B$-invariant. For $f\neq0$ it is linear in $x$, so its fibers are graphs, and it is an
invariant function (Lemma \ref{lem:invariant-function}).
 
\emph{Faces.} The faces $\mathcal O_3^\pm$ have pole type with $k=2$ and $Q|_{\{z=0\}}=y^2$, so both
lie at $+\infty$ and are constant. The face $\mathcal O_4$ has logarithmic type and lies at
$\sgn(bf)\cdot\infty$, and it is constant.
 
\begin{proposition}\label{prop:B31}
A rank $2$ structure on $\T^3$ with respect to $N_{3,1}$ exists if and only if, in the normal
form above, $f\neq0$ and $bf<0$; conjugating by a diagonal and then by an upper
triangular matrix, one may take $f=1$, $b=-1$, $e=0$. There are two bricks,
$\mathcal O\cup \mathcal O_3^+\cup \mathcal O_4$ and $\mathcal O\cup \mathcal O_3^-\cup \mathcal O_4$, up to reflection. The chain alternates
between a face $\mathcal O_3^\pm$ and the face $\mathcal O_4$, with a free choice of sign at each crossing of
$\{z=0\}$, so there are infinitely many structures whose universal covers are pairwise non-isomorphic by isomorphisms normalizing $N$ (Proposition \ref{prop:words-distinct}). None of these structures is
affine.
\end{proposition}
 
\begin{proof}
A brick needs two faces at opposite ends (Lemma \ref{lem:bricks}(1)), which forces $\mathcal O_4$ to be free
and at $-\infty$, i.e.\ $bf<0$. Since all faces are constant, the chain alternates between
the two ends (proof of Lemma \ref{lem:constant-faces}), which gives the description of the words; they are
realized by Lemma \ref{lem:realize}. The normalization is as follows: a diagonal conjugation multiplies
$bf$ by a positive number and can be used to make $f=1$ and $b=-1$, and conjugation by
$\exp(sE_{23})\in N$ replaces $e$ by $e-s$ and fixes $X$. Finally, the only $B$-invariant
hyperplanes are $\{z=0\}$ and $\{w=0\}$, and every word crosses both; since the Zariski closure of the lattice $\Gamma_0$ contains $B$, no hyperplane invariant under the holonomy misses the developing
image, so no structure is affine (Remark \ref{rem:affine}).
\end{proof}
 
\subsection{\texorpdfstring{$N_{2,2}$}{N22} (Exercise 4.14)}
 
Here $\mathfrak n=\langle D_0,E_{12},E_{34}\rangle$ is abelian, with $D_0=\diag(1,1,-1,-1)$. The
open orbit $\mathcal O=\{y>0,\ w>0\}$ has the faces
\[
  \mathcal O_2^\pm=\{y=0,\ \pm x>0,\ w>0\},\qquad \mathcal O_4^\pm=\{w=0,\ \pm z>0,\ y>0\}.
\]
 
\emph{Normal forms.} The algebra $\mathfrak b$ is not $\langle E_{12},E_{34}\rangle$, since
$E_{12}$ fixes $\mathcal O_2^\pm$ and $E_{34}$ fixes $\mathcal O_4^\pm$ pointwise. So
$\mathfrak b\cap\langle E_{12},E_{34}\rangle$ is spanned by one element $pE_{12}+qE_{34}$, and
exchanging the two blocks if necessary we obtain:
\[
  \text{(a)}\ X=E_{12},\quad Y=D_0+aE_{34};\qquad
  \text{(b)}\ X=E_{12}+bE_{34}\ (b\neq0),\quad Y=D_0+aE_{34}.
\]
In case (a), $B$ is free on $\mathcal O$ and on $\mathcal O_4^\pm$ and fixes $\mathcal O_2^\pm$ pointwise. In case (b),
it is free on $\mathcal O$, $\mathcal O_2^\pm$ and $\mathcal O_4^\pm$: the minors include $2wy^2$ on $\mathcal O$, $-2bw^2x$ on
$\mathcal O_2^\pm$ and $2y^2z$ on $\mathcal O_4^\pm$.
 
\emph{Invariant function.} The function
\[
  f_{2,2}=-b\frac xy+\frac zw+\frac a2\log\frac wy
\]
is $B$-invariant, with $b=0$ in case (a); it is linear in $z$ with coefficient $1/w$, so it
is an invariant function.
 
\emph{Faces.} All faces in use have pole type with $k=1$: $Q=-bx$ at $\mathcal O_2^\pm$ and $Q=z$ at
$\mathcal O_4^\pm$. So each pair lies at opposite ends, and all faces are alternating. In $\mathcal O$, the
face $\mathcal O_4^+$ lies at $+\infty$, and $\mathcal O_2^+$ lies at $-\sgn(b)\cdot\infty$.
 
\begin{proposition}\label{prop:B22}
Let $V=\T^3$ carry a rank $2$ structure with respect to $N_{2,2}$.
\begin{itemize}
  \item In case (a) there is one brick, $D=\mathcal O\cup \mathcal O_4^+\cup \mathcal O_4^-$; the words are its powers,
  the chain alternates between $\mathcal O$ and $R_4 \mathcal O$, and the structure is affine, since it avoids
  $\{y=0\}$.
  \item In case (b), say with $b<0$ (the case $b>0$ is obtained by exchanging $\mathcal O_2^+$ and
  $\mathcal O_2^-$), there are four bricks,
  \[
    A=\mathcal O\cup \mathcal O_2^+\cup \mathcal O_2^-,\quad B'=\mathcal O\cup \mathcal O_4^+\cup \mathcal O_2^-,\quad C=\mathcal O\cup \mathcal O_2^+\cup \mathcal O_4^-,\quad
    D=\mathcal O\cup \mathcal O_4^+\cup \mathcal O_4^-.
  \]
  Since all faces are alternating, the chain enters every copy at the same end, and at
  every step it can cross either $\{y=0\}$ or $\{w=0\}$. So there are infinitely many structures whose universal covers are pairwise non-isomorphic by isomorphisms normalizing $N$ (Proposition \ref{prop:words-distinct}). The affine ones are the powers of $A$, which avoid
  $\{w=0\}$, and the powers of $D$, which avoid $\{y=0\}$.
\end{itemize}
Conversely, all of these words are realized by Lemma \ref{lem:realize}.
\end{proposition}
 
\subsection{\texorpdfstring{$N_{2,1,1}$}{N211} (Exercise 4.15)} \label{sec:N211}
 
Here $\mathfrak n=\mathfrak a\oplus\langle E_{12}\rangle$ is abelian, where $\mathfrak a$ consists of the
traceless matrices $\diag(\lambda,\lambda,\mu,\nu)$. The open orbit $\mathcal O=\{y>0,\ z>0,\ w>0\}$
has the faces
\[
  \mathcal O_2^\pm=\{y=0,\ \pm x>0,\ z,w>0\},\qquad \mathcal O_3=\{z=0,\ y,w>0\},\qquad
  \mathcal O_4=\{w=0,\ y,z>0\}.
\]
 
\emph{Normal forms.} If $\mathfrak b=\mathfrak a$, then $x/y$ is an invariant function, $B$ is free
on $\mathcal O$ and on $\mathcal O_2^\pm$ but not on $\mathcal O_3$ or $\mathcal O_4$, and the only brick is $\mathcal O\cup \mathcal O_2^\pm$; the
structures are affine. Otherwise $\mathfrak b\cap\mathfrak a$ is spanned by one element
$X=\diag(\lambda,\lambda,\mu,\nu)\neq0$, and $\mathfrak b$ is spanned by $X$ and
$Y=\diag(a,a,b,c)+E_{12}$, with $2\lambda+\mu+\nu=2a+b+c=0$. Put
\[
  K=(a-c)(\lambda-\mu)-(a-b)(\lambda-\nu).
\]
Then $K=0$ if and only if $E_{12}\in\mathfrak b$, since by the trace conditions $K=0$ exactly
when the diagonal part of $Y$ is proportional to $X$.
 
\emph{Freeness.} $B$ is free on $\mathcal O$. It is free on $\mathcal O_2^\pm$ if and only if $K\neq0$, on $\mathcal O_3$
if and only if $\lambda\neq\nu$, and on $\mathcal O_4$ if and only if $\lambda\neq\mu$; the relevant
minors are $wxzK$, $-wy^2(\lambda-\nu)$ and $-y^2z(\lambda-\mu)$, up to sign.
 
\emph{Invariant function.} The function
\[
  f_{2,1,1}=K\frac xy+(\lambda-\mu)\log w-(\lambda-\nu)\log z+(\mu-\nu)\log y
\]
is $B$-invariant. If $K\neq0$ it is linear in $x$; if $K=0$ it is a nonzero linear function of
$\log(z/y)$ and $\log(w/y)$, since $X\neq0$. In both cases its fibers are connected, and it
is an invariant function.
 
\emph{Faces.} For $K\neq0$, the faces $\mathcal O_2^\pm$ have pole type with $k=1$ and $Q=Kx$, so they
lie at opposite ends ($\mathcal O_2^+$ at $\sgn(K)\cdot\infty$) and are alternating. The faces $\mathcal O_3$
and $\mathcal O_4$ have logarithmic type and are constant; $\mathcal O_3$ lies at $\sgn(\lambda-\nu)\cdot\infty$
and $\mathcal O_4$ at $\sgn(\mu-\lambda)\cdot\infty$. Hence, when both are free, $\mathcal O_3$ and $\mathcal O_4$ lie at
the same end if and only if $\lambda$ lies strictly between $\mu$ and $\nu$.
 
\begin{proposition}\label{prop:B211}
Let $V=\T^3$ carry a rank $2$ structure with respect to $N_{2,1,1}$, with $\mathfrak b\neq\mathfrak a$.
\begin{itemize}
  \item[(a)] If $E_{12}\in\mathfrak b$ ($K=0$), a structure exists if and only if
  $\lambda\neq\mu,\nu$ and $\lambda$ does not lie between $\mu$ and $\nu$. There is one brick
  $\mathcal O\cup \mathcal O_3\cup \mathcal O_4$, the chain has period $4$, and the structure is affine, since it avoids
  $\{y=0\}$.
  \item[(b)] If $K\neq0$ and $\lambda=\mu$, $\lambda=\nu$, or $\lambda$ lies strictly between $\mu$
  and $\nu$, then the free constant faces lie at one end, so by Lemma \ref{lem:constant-faces} they are never
  crossed. The only brick is $\mathcal O\cup \mathcal O_2^\pm$, the words are its powers, and the structure is
  affine.
  \item[(c)] If $K\neq0$, $\lambda\neq\mu,\nu$ and $\lambda$ does not lie between $\mu$ and $\nu$,
  then $\mathcal O_3$ and $\mathcal O_4$ lie at opposite ends. The bricks over $\mathcal O$ are the four sets $\mathcal O\cup F_+\cup
  F_-$, where $F_+$ and $F_-$ are faces at $+\infty$ and $-\infty$ chosen among
  $\mathcal O_2^\pm,\mathcal O_3,\mathcal O_4$; so there are infinitely many structures whose universal covers are pairwise non-isomorphic by isomorphisms normalizing $N$ (Proposition \ref{prop:words-distinct}). The affine ones
  are the powers of $\mathcal O\cup \mathcal O_2^\pm$, which avoid $\{z=0\}$ and $\{w=0\}$, and the powers of
  $\mathcal O\cup \mathcal O_3\cup \mathcal O_4$, which avoid $\{y=0\}$.
\end{itemize}
Conversely, all of these words are realized by Lemma \ref{lem:realize}.
\end{proposition}
 
\begin{proof}
(a) Only $\mathcal O_3$ and $\mathcal O_4$ are available, and they must lie at opposite ends. (b) If
$\lambda=\mu$ only $\mathcal O_3$ is free, if $\lambda=\nu$ only $\mathcal O_4$ is free, and if $\lambda$ lies strictly
between $\mu$ and $\nu$ both are free but lie at the same end. (c) The constant faces lie at
both ends, so they can be crossed infinitely often. A word avoiding $\mathcal O_4$ can cross the
constant face $\mathcal O_3$ at most once, by the argument in the proof of Lemma \ref{lem:constant-faces}, so by
Lemma \ref{lem:bricks}(3) it does not cross $\mathcal O_3$ at all; similarly for a word avoiding $\mathcal O_3$. A word
avoiding $\mathcal O_2^\pm$ alternates between $\mathcal O_3$ and $\mathcal O_4$. This gives the description of the
affine structures, by Remark \ref{rem:affine}.
\end{proof}
 
\begin{example}\label{ex:N211}
Take $X=\diag(1,1,-1,-1)$ and $Y=\diag(1,1,0,-2)+E_{12}$. Then $\mu=\nu=-1$, so $\lambda=1$ does not
lie between $\mu$ and $\nu$, and $K=4$. The faces $\mathcal O_2^+$ and $\mathcal O_3$ lie at $+\infty$, and $\mathcal O_2^-$
and $\mathcal O_4$ at $-\infty$. With $Y=\diag(1,1,-2,0)+E_{12}$ instead, $K=-4$, and $\mathcal O_2^+$ and $\mathcal O_4$
lie at the same end. So which faces lie at the same end depends on $Y$, not only on $X$.
\end{example}

\section{Maximal abelian subgroups of \texorpdfstring{$\SL(5,\R)$}{SL5R} which are not unipotent}\label{appendix:nonunipotent}

Recall by Lemma \ref{lem:real-blocks} that maximal abelian subalgebras in $\gl(5,R)$ are sums of real blocks. 
The maximal abelian subalgebras of $\mathfrak{sl}(4,\R)$ are
listed in \cite{WZ84}; the list is corrected and completed in \cite{diatta2022classification} Section 5.2, 
who correct the representative of item 8 and add the missing decomposable item $\mathfrak j_2\oplus\mathfrak j_2$. The first
appearance of a group with more than one real form is in \cite{WZ84}, given without proof,
and Willem de Graaf gives a proof using Galois cohomology in Example \ref{ex:warmup}. In this
appendix we list the resulting groups in $\SL(5,\R)$ and determine which of them
contain the group $B$ of a rank 4 projective structure.
 
\subsection{ Block algebras and sums}
 
A \emph{block algebra} of size $m$ is one of the abelian subalgebras
$\bb\subset\gl(m,\R)$ in the table below. Within a block, $E_{ij}$ denotes the
elementary matrix. For blocks of size $3$ and $4$ we use the numbering of \cite{WZ84}, which
matches that of \cite{diatta2022classification} (where $\w j$ is called $\mathcal Y_j$): $\w1,\dots,\w8$ are the
maximal abelian subalgebras of $\mathfrak{sl}(4,\R)$ numbered $1$--$8$ there, with the same
representatives, and $\w9,\w{10},\w{11}$ are the $3\times3$ blocks of the decomposable
algebras numbered $9$--$11$. The column $(\lambda\mu\nu)$ gives the Kravchuk signature
of a nilpotent block, as in \cite{WZ84}. The column \cite{suprunenko1968commutative} gives the complex class in the
lists of \cite{suprunenko1968commutative} Chapter 3, Section 3, whose algebras are lower triangular; $N_i^T$ ($N_i$ in the list for $n=m$)
denotes the transpose, and $\sim$ means conjugate.
 
\begin{center}
\renewcommand{\arraystretch}{1.2}
\begin{tabular}{cllll}
size & name & basis & $(\lambda\mu\nu)$ & [ST03] \\ \hline
$1$ & $\bd$ & --- & & \\
$2$ & $\bj_2$ & $E_{12}$ & (101) & $N_1^T$ \\
$2$ & $\bc$ & $E_{12}-E_{21}$ & & not nilpotent \\ \hline
$3$ & $\w9$ & $E_{12},\ E_{13}$ & (102) & $N_1^T$ \\
$3$ & $\w{10}$ & $E_{13},\ E_{23}$ & (201) & $N_2^T$ \\
$3$ & $\w{11}$ & $E_{12}+E_{23},\ E_{13}$ & (111) & $N_3^T$ \\ \hline
$4$ & $\w1$ & $E_{13},\ E_{14},\ E_{23},\ E_{24}$ & (202) & $N_2^T$ \\
$4$ & $\w2$ & $E_{12},\ E_{13},\ E_{14}$ & (103) & $N_1^T$ \\
$4$ & $\w3$ & $E_{14},\ E_{24},\ E_{34}$ & (301) & $N_3^T$ \\
$4$ & $\w4$ & $E_{13}+E_{34},\ E_{14},\ E_{24}$ & (211) & $\sim N_5^T$ \\
$4$ & $\w5$ & $E_{12}+E_{24},\ E_{13},\ E_{14}$ & (112) & $\sim N_4^T$ \\
$4$ & $\w6^{\varepsilon}$ & $E_{12}+E_{24},\ E_{13}+\varepsilon E_{34},\ E_{14}$,
  \ $\varepsilon=\pm1$ & (121) & $\sim N_6^T$ \\
$4$ & $\w7$ & $E_{12}+E_{23}+E_{34},\ E_{13}+E_{24},\ E_{14}$ & (121) & $N_7^T$ \\
$4$ & $\w8$ & $X_0,\ N,\ X_0N$ & & not nilpotent
\end{tabular}
\end{center}
 
Here $X_0=E_{12}-E_{21}+E_{34}-E_{43}$, $N=E_{13}+E_{24}$ and $X_0N=E_{14}-E_{23}$; this
is the representative given by \cite{diatta2022classification} Section 5.2, item 8, which corrects the
non-commutative formula printed in \cite{WZ84}. The algebra $\w8$ is the realification of a
$2\times2$ complex Jordan block: $X_0^2=-I$ is a complex structure and $N$ is a
complex-linear nilpotent. All blocks except $\bc$ and $\w8$ consist of nilpotent
matrices.

The three algebras of size~$3$ are all the real forms of the corresponding complex
classes: their normalizers in $\GL(3,\C)$ are the stabilizer of the line $\C e_1$ (for
$\w9$), the stabilizer of the plane $\langle e_1,e_2\rangle$ (for $\w{10}$), and the
connected split solvable group $T'U_3$ with
$T'=\{\diag(t_1,t_2,t_3):t_1t_3=t_2^2\}$ (for $\w{11}$), all of which have trivial first
Galois cohomology, so by Theorem \ref{thm:galco} each class has a unique real form.
 
Two remarks explain the pattern of the table. The anti-transpose
$X\mapsto w_0X^{T}w_0$, where $w_0$ is the antidiagonal permutation matrix, exchanges
$\w2\leftrightarrow\w3$, $\w4\leftrightarrow\w5$ and $\w9\leftrightarrow\w{10}$, and
preserves the conjugacy classes of the other blocks. The algebras $\w6^{+}$ and
$\w6^{-}$ are the two real forms of the single complex class $N_6$, as in Remark \ref{rem:real-forms},

and Willem de Graaf gives a proof using Galois Cohomology in Example \ref{ex:warmup}.

\begin{definition}
Let $n_1+\dots+n_k=5$ and let $\bb_i\subset\gl(n_i,\R)$ be block algebras of sizes
$n_i$. Let $\mathfrak a$ be the space of trace-zero diagonal matrices that are scalar on
each block,
\[
  \mathfrak a=\{\diag(t_1I_{n_1},\dots,t_kI_{n_k})\ :\ n_1t_1+\dots+n_kt_k=0\}.
\]
We write $A(\bb_1,\dots,\bb_k)\subset\SL(5,\R)$ for the connected abelian group with Lie
algebra $\mathfrak a\oplus\bb_1\oplus\dots\oplus\bb_k$, the $\bb_i$ placed as diagonal
blocks in this order.
\end{definition}

By Lemma~\ref{lem:real-blocks}, every non-unipotent maximal abelian subalgebra of $\mathfrak{sl}(5,\R)$ is
conjugate to $\mathfrak a\oplus\mathfrak b_1\oplus\dots\oplus\mathfrak b_k$ with $k\ge2$
blocks from the table above. Conversely, each such sum is maximal abelian: $\mathfrak a$
contains an element with distinct scalars on the blocks, so the centralizer is block
diagonal, and within each block the centralizer of $\R I_{n_i}\oplus\mathfrak b_i$ is
itself. Enumerating the multisets of blocks whose sizes sum to $5$ gives the $24$ groups
of Table \ref{table:nonunip}.
 
For example, a rotation block contributes
$e^{s}\left(\begin{smallmatrix}\cos\theta & \sin\theta\\ -\sin\theta&\cos\theta
\end{smallmatrix}\right)$, and
\[
  A(\w7,\bd)=\left\{\begin{pmatrix}
  a&b&c&d&0\\ 0&a&b&c&0\\ 0&0&a&b&0\\ 0&0&0&a&0\\ 0&0&0&0&a^{-4}
  \end{pmatrix} : a>0\right\}.
\]
Permuting the blocks gives a conjugate group, so only the multiset of blocks matters;
we list larger blocks first. The groups of this form are listed in Table \ref{table:nonunip} organized by block size. %
All groups have dimension $4$ except $A(\w1,\bd)$, which has
dimension $5$.

\subsection{ A criterion for rank 4 holonomy}

For the non-unipotent groups the rank computation splits over blocks. Call a block algebra
$\bb\subset\gl(m,\R)$ \emph{good} if $\R I_m\oplus\bb$ has an open orbit in $\R^m$, that
is, if $(\R I_m\oplus\bb)\,q=\R^m$ (compare to the projective version in Lemma \ref{lem:C-criterion}).
 
\begin{proposition}\label{prop:C-blocks}
The group $A(\bb_1,\dots,\bb_k)$ has an open orbit in $\Sph^4$ if and only if every
block $\bb_i$ is good.
\end{proposition}
 
\begin{proof}
Let $\h=\mathfrak a\oplus\bb_1\oplus\dots\oplus\bb_k$ and $B=\R I_5\oplus\h$. Since
$\R I_5\oplus\mathfrak a$ is the space of all diagonal matrices that are scalar on each
block,
\[
  B=\bigoplus_{i=1}^k\bigl(\R I_{n_i}\oplus\bb_i\bigr),
\]
acting block by block on $\R^5=\R^{n_1}\oplus\dots\oplus\R^{n_k}$. For
$p=(p_1,\dots,p_k)$ we have $\h p+\R p=Bp=\bigoplus_i(\R I_{n_i}\oplus\bb_i)p_i$, so
$r(p)=5$ if and only if $(\R I_{n_i}\oplus\bb_i)p_i=\R^{n_i}$ for every $i$. By
Lemma~\ref{lem:C-criterion}, $A(\bb_1,\dots,\bb_k)$ has an open orbit if and only if such
$p_i$ exist for every $i$.
\end{proof}
 
\begin{corollary}\label{cor:C-good}
The blocks $\bd$, $\bj_2$, $\bc$, $\w{10}$, $\w{11}$, $\w3$, $\w4$, $\w6^{\pm}$, $\w7$,
$\w8$ are good, and the blocks $\w9$, $\w1$, $\w2$, $\w5$ are bad.
\end{corollary}
 
\begin{proof}
For each good block of size $m$, the last basis vector $q=e_m$ works: the basis
elements of the block send $e_m$ to vectors which, together with $e_m$, span $\R^m$.
For example, for $\w7$ we have $Je_4=e_3$, $J^2e_4=e_2$, $J^3e_4=e_1$; for $\w8$ we have
$X_0e_4=e_3$, $Ne_4=e_2$, $X_0Ne_4=e_1$; for $\bc$ we have $(E_{12}-E_{21})e_2=e_1$.
The other good blocks are checked in the same way.
 
For the bad blocks, let $q=(q_1,\dots,q_m)$. The basis elements of $\w9$ and $\w2$ take
values in $\R e_1$, so $\dim(\R I\oplus\bb)q\le2<m$. The basis elements of $\w5$ send $q$
to $q_2e_1+q_4e_2$, $q_3e_1$, $q_4e_1$, and those of $\w1$ take values in
$\R e_1+\R e_2$; in both cases $\dim(\R I\oplus\bb)q\le3<4$.
\end{proof}
 
\begin{example}[A group with no open orbit]\label{ex:C-no}
For $A(\w2,\bd)$, the Lie algebra is spanned by $D=\diag(1,1,1,1,-4)$, $E_{12}$,
$E_{13}$, $E_{14}$. For every $p=(v,w,x,y,z)$,
\[
  E_{12}p=w\,e_1,\qquad E_{13}p=x\,e_1,\qquad E_{14}p=y\,e_1
\]
are proportional, so $\delta\equiv0$. This group contains no group $B$ of a rank 4 structure.
\end{example}
 
\begin{example}[A group with open orbits, in full]\label{ex:C-full}
Consider $H=A(\w7,\bd)$, displayed above. Its Lie algebra has basis
\[
  D=\diag(1,1,1,1,-4),\qquad J=E_{12}+E_{23}+E_{34},\qquad J^2=E_{13}+E_{24},\qquad
  J^3=E_{14}.
\]
For $p=(v,w,x,y,z)$ we have
\[
  Dp=(v,w,x,y,-4z),\quad Jp=(w,x,y,0,0),\quad J^2p=(x,y,0,0,0),\quad J^3p=(y,0,0,0,0).
\]
Subtracting $p$ from the first column does not change $\delta$, and replaces $Dp$ by
$-5z\,e_5$. Expanding along that column,
\[
  \delta(p)=\det\begin{pmatrix}
  0&w&x&y&v\\ 0&x&y&0&w\\ 0&y&0&0&x\\ 0&0&0&0&y\\ -5z&0&0&0&z
  \end{pmatrix}
  =-5z\det\begin{pmatrix}
  w&x&y&v\\ x&y&0&w\\ y&0&0&x\\ 0&0&0&y
  \end{pmatrix}
  =-5z\cdot y\cdot(-y^3)=5\,y^4z .
\]
So $\delta\not\equiv0$, and by Lemma~\ref{lem:C-criterion} the open orbits of $H$ are
the four components
\[
  \Omega_{\pm\pm}=\{[v:w:x:y:z]\in\Sph^4:\ \pm y>0,\ \pm z>0\}.
\]
 
The structure can be written down explicitly. Take $\Omega=\Omega_{++}$. The group
$H\cong\R^4$ is simply connected, $\Omega$ is diffeomorphic to $\R^4$, and the orbit map
$H\to\Omega$ is an equivariant local diffeomorphism between homogeneous spaces, hence a
covering map; since $\Omega$ is simply connected it is a diffeomorphism. So $H$ acts
simply transitively on $\Omega$, and for any lattice $\Gamma\cong\Z^4$ in $H$ the
quotient $\Gamma\backslash\Omega\cong\Gamma\backslash H$ is a torus $\T^4$, and the left-invariant structure defined by $(H,[e_4+e_5],\Gamma)$ is a rank 4 structure on it with holonomy $\Gamma$. 
 
In the affine chart $z=1$, the domain $\Omega$ is the half-space $\{y>0\}\subset\R^4$,
so these structures are affine but not complete.
\end{example}

By Proposition~\ref{prop:C-blocks} and Corollary~\ref{cor:C-good}, a group has an open
orbit exactly when none of its blocks is $\w1,\w2,\w5,\w9$. To describe the open orbits
we also list $\delta$, given up to a nonzero constant in the coordinates
$p=(v,w,x,y,z)$ adapted to the block order; by the proof of
Proposition~\ref{prop:C-blocks} for the groups other than $A(\mathfrak w_1,\mathfrak d)$ it is a product of one factor for each block. We leave
the computations to the reader but present the summary in Table \ref{table:nonunip}.

\begin{table}[htb]
\caption{The non-unipotent maximal abelian subgroups and their open orbits}
\label{table:nonunip}
\begin{center}
\renewcommand{\arraystretch}{1.15}
\begin{tabular}{llll}
Block sizes & Group & $\delta(p)$ & Open orbit \\ \hline
$(4,1)$ & $A(\w1,\bd)$ & $\dim H=5$, $r\le4$ & no \\
        & $A(\w2,\bd)$ & $0$ & no \\
        & $A(\w3,\bd)$ & $y^4z$ & yes \\
        & $A(\w4,\bd)$ & $y^4z$ & yes \\
        & $A(\w5,\bd)$ & $0$ & no \\
        & $A(\w6^{+},\bd)$ & $y^4z$ & yes \\
        & $A(\w6^{-},\bd)$ & $y^4z$ & yes \\
        & $A(\w7,\bd)$ & $y^4z$ & yes \\
        & $A(\w8,\bd)$ & $z(x^2+y^2)^2$ & yes \\ \hline
$(3,2)$ & $A(\w9,\bj_2)$ & $0$ & no \\
        & $A(\w{10},\bj_2)$ & $x^3z^2$ & yes \\
        & $A(\w{11},\bj_2)$ & $x^3z^2$ & yes \\
        & $A(\w9,\bc)$ & $0$ & no \\
        & $A(\w{10},\bc)$ & $x^3(y^2+z^2)$ & yes \\
        & $A(\w{11},\bc)$ & $x^3(y^2+z^2)$ & yes \\ \hline
$(3,1,1)$ & $A(\w9,\bd,\bd)$ & $0$ & no \\
          & $A(\w{10},\bd,\bd)$ & $x^3yz$ & yes \\
          & $A(\w{11},\bd,\bd)$ & $x^3yz$ & yes \\ \hline
$(2,2,1)$ & $A(\bj_2,\bj_2,\bd)$ & $w^2y^2z$ & yes \\
          & $A(\bj_2,\bc,\bd)$ & $w^2z(x^2+y^2)$ & yes \\
          & $A(\bc,\bc,\bd)$ & $z(v^2+w^2)(x^2+y^2)$ & yes \\ \hline
$(2,1,1,1)$ & $A(\bj_2,\bd,\bd,\bd)$ & $w^2xyz$ & yes \\
            & $A(\bc,\bd,\bd,\bd)$ & $xyz(v^2+w^2)$ & yes \\ \hline
$(1,1,1,1,1)$ & $A(\bd,\bd,\bd,\bd,\bd)$ & $vwxyz$ & yes
\end{tabular}
\end{center}
\end{table}

In every positive case the open orbits are the components of $\{\delta\neq0\}$ in
$\Sph^4$. A factor such as $x^2+y^2$ comes from a rotation block or from the complex
structure of $\w8$, whose orbits wind around the subspace where the corresponding
coordinates vanish.

As for the two half-plane structures of Example \ref{ex:T2}, some of these groups can be dehomogenized in different ways, giving several affine structures.

\begin{remark}[Comparison with \cite{diatta2022classification}]
A block $\bb$ is good exactly when $\R I_m\oplus\bb$ has an open orbit in $\R^m$. By
\cite{diatta2022classification} Theorem 1 and Section 3.2 , this holds exactly when the $2m$-dimensional Lie
algebra $(\R I_m\oplus\bb^{T})\ltimes\R^m$ is Frobenius.  Indeed, $\mathbb RI\oplus\mathfrak b$ has an open orbit on $\mathbb R^m$ if and only if $\mathbb RI\oplus\mathfrak b^T$ has an open orbit on $(\mathbb R^m)^*$ under the dual action, since $(\xi X)^T=X^T\xi^T$. As $\bb^T$ is conjugate to the
anti-transpose of $\bb$, Corollary \ref{cor:C-good} can be read off from the lists in
\cite{diatta2022classification}, Section 5.2 : among $\mathcal Y_2,\dots,\mathcal Y_8$ (the algebra
$\mathcal Y_1=\w1$ has the wrong dimension for their theorem) exactly $\mathcal Y_2$,
$\mathcal Y_5$, $\mathcal Y_{6,\pm1}$, $\mathcal Y_7$, $\mathcal Y_8$ give Frobenius
algebras, and their anti-transposes are $\w3,\w4,\w6^{\pm},\w7,\w8$; and among the
$3\times3$ algebras, $L_{2,5}=\w9$ and $L_{2,6}=\w{11}$ give Frobenius algebras while
$L_{2,4}=\w{10}$ does not, whose anti-transposes are $\w{10}$, $\w{11}$, $\w9$. This
gives an independent confirmation of the table.
\end{remark}

\end{document}